\documentclass[12pt]{article}
\usepackage[latin1]{inputenc}
\usepackage{color}
\usepackage{dsfont}
\usepackage{tikz}
\usetikzlibrary{decorations.pathreplacing}
\usepackage{float}
\usepackage{multicol}
\usepackage[normalem]{ulem}
\usepackage{latexsym,epsfig}
\usepackage{amssymb,amsmath,amsfonts,amscd}
\usepackage{exscale}
\usepackage{ifthen}
\usepackage{longtable}
\usepackage{subfigure}
\usepackage{cite}
\usepackage{lscape}
\usepackage{todonotes}
\usepackage{enumitem}
\usepackage{hyperref}

 \newtheorem{theo}{\textbf{Theorem}}

\newtheorem{lemma}[theo]{\textbf{Lemma} }
\newtheorem{coro}[theo]{\textbf{Corollary} }

\newtheorem{prop}[theo]{\textbf{Proposition} }

\newenvironment{proof}[1][\unskip]{\par \noindent{\bf Proof #1. }}{\hfill $\Box$\par} 

\newcommand{\figu}[2]{\includegraphics[width=#1\textwidth]{#2}}

\newcommand{\E}{\mathbb{E}}
\renewcommand{\P}{\mathbb{P}}
\def \Q{\mathbb{Q}}

\def \N{\mathbb{N}}

\def \Z{\mathbb{Z}}
\def \R{\mathbb{R}}

\newcommand{\vv}{\mathbf{v}}

\def \cL{\mathcal{L}} 
\def \cI{\mathcal{I}} 
\def \cC{\mathcal{C}}

\def \cL{\mathcal{L}}

\def \cZ{\mathcal{Z}}

\def \fz{\mathfrak{z}}

\newcommand{\PRd}{\mathbf{P}(\R^d)} 

\renewcommand{\Pr}[1]{P_{#1}}

\def \AG{\mathcal{A}_1(\R^d)}

\newcommand{\vvt}{{\widetilde{\mathbf{v}}}}

\newcommand{\Aff}{{\mathrm{Aff}}}
\newcommand{\AffR}{{\mathrm{Aff}(\R)}}
\newcommand{\AffRd}{{\mathrm{Aff}(\R^d)}}
\newcommand{\GL}{\operatorname{GL}}
\newcommand{\supp}{{\mathrm{supp}}}

\newcommand{\sil}{{\sigma_\ell}}
\newcommand{\sir}{{\sigma_r}}

\newcommand{\ev}{\textrm{ev}}
\title{Local contraction of stochastic affine recursions of $\R^d$}
\date{}
\author{
Richard Aoun\thanks{
Laboratoire d'Analyse et de Math\'ematiques Appliqu\'ees (LAMA), UMR 8050, 
Universit\'e Gustave Eiffel, Universit\'e Paris-Est Cr\'eteil, CNRS, France. 
Email: \texttt{richard.aoun@univ-eiffel.fr}
}
\and
Sara Brofferio\thanks{
Laboratoire d'Analyse et de Math\'ematiques Appliqu\'ees (LAMA), UMR 8050, 
Universit\'e Paris-Est Cr\'eteil, Universit\'e Gustave Eiffel, CNRS, France. 
Email: \texttt{sara.brofferio@u-pec.fr}
}
\and
Marc Peign\'e\thanks{
Institut Denis Poisson (IDP), UMR 7013, 
Universit\'e de Tours, Universit\'e d'Orl\'eans, CNRS, France. 
Email: \texttt{marc.peigne@univ-tours.fr}
}
}

\date{}
\begin{document}
\maketitle

 \begin{center}
\textit{In   loving memory of    Martine Babillot.}    
 \end{center}
\begin{abstract}
We consider affine stochastic recursions $X^x_{n+1}=A_{n+1}X^x_n+B_{n+1}$ in $\R^d$,
where $(B_n,A_n)$ are i.i.d.\ in $\Aff(\R^d)$ with strongly irreducible and proximal
linear part. In the critical regime where the top Lyapunov exponent vanishes, we show
that the recursion is locally contractive in every direction,
yielding uniqueness, up to a positive multiplicative constant, of the $\mu$-invariant non-zero Radon measure. 
This extends to    every dimension $d\geq 2$ a result of
Babillot-Bougerol-\'Elie  (1997) in dimension one.
The proof is based on the study of the asymptotic behaviour of the ``joint process'' $(X^x_n,X^y_n)$,
using geometric and measure-theoretic properties of the diagonal affine action.
A key ingredient is an abstract criterion for the total dissipativity of semigroup actions, which we apply to this process.
\end{abstract}
 \bigskip
{\small 

}

  \tableofcontents

\section{Introduction}

\subsection{Main results}\label{sec.1.1}

 Let $d\geq 1$ and endow $\R^d$ with the Euclidean norm $\vert \cdot \vert$. Denote by $\GL_d(\R)$   the group of invertible $d\times d$ real matrices 
and by $\AffRd = \R^d \rtimes \GL_d(\R)$ the group of invertible affine transformations of $\R^d$. Each element of $\AffRd$ is   denoted by $\gamma=(B,A)$ with $B\in \R^d$ and $A\in \GL_d(\R)$.  The group $\AffRd$ acts    on $\R^d$ by $\gamma=(B,A): x\mapsto Ax+B$.

Let $(\gamma_n = (B_n, A_n))_{n \geq 1}$ be a sequence of 
i.i.d.\ $\AffRd$-valued random variables defined on a probability space $(\Omega, \mathcal{F}, \P)$ and with common distribution $\mu$. 
In this   paper, we study  the stochastic affine recursion on $\R^d$ defined by 
\begin{equation}\label{eq-rec.affine}
 X_{n+1}^x := \gamma_{n+1} X_n^x = A_{n+1} X_n^x + B_{n+1}, 
\qquad X_0^x = x \in \R^d.   
\end{equation}
This Markov process on $\R^d$ coincides with   
the orbit $(L_n x)$ of $x\in\R^d$ under the left random walk 
$L_n := \gamma_n \cdots \gamma_1$ on the group  $\AffRd$. 
The central question of the present paper is how  the process  $(X_n^x)$ is sensitive 
 to its the starting point $x$. Since for any $x,y\in \R^d$,  
\[
X_n^x - X_n^y = A_n\cdots A_1(x-y), 
\]
the study of $(X_n^x)_n$  relies crucially on the understanding of the left random walk $(A_n\cdots A_1)_n$ on $\GL_d(\R)$, which has been 
the subject of extensive research since the 1960s. When $\E(\ln^+ \vert A_1 \vert )<+\infty$, we denote by $
\lambda_1 = \lim_{n\to\infty} \frac{1}{n}\E(\ln | A_n\cdots A_1 |)$ the \emph{top Lyapunov exponent} of this product.

\paragraph{Local contraction for the affine recursion.} When $\lambda_1 < 0$, the  sequence $(A_n\cdots A_1)$
tends almost surely to zero. Hence the affine recursion is \emph{globally 
	contracting}:  
$|X_n^x - X_n^y| \to 0$ a.s, for every $x,y\in \R^d$.

The case 
 $\lambda_1 = 0$ is of a fundamentally different nature: the norm of the 
 linear part   oscillates between $0$ and $+\infty$ under generic conditions. 
 Nevertheless, despite the absence of global contraction, a landmark result of Babillot--Bougerol--\'Elie \cite{BBE} shows that 
 a subtle notion of contraction exists in dimension $d=1$.   Under 
  non-degeneracy and moment conditions,  they showed   in 
dimension $d=1$   that,     for every 
$x\in\R$ and every compact interval $I\subset\R$,
\begin{equation}\label{eq.local.contraction.1d}
	\forall\, y\in\R,\;
	|X_n^x - X_n^y|\,\mathds{1}_{X_n^x \in I}
	\underset{n\to+\infty}{\longrightarrow} 0\mbox{ a.s.}
\end{equation}
In other words, whenever the process visits a bounded region, 
trajectories starting from different points must eventually be close.
This phenomenon is known as \emph{local contraction}.

The goal  of this paper is to extend this local contraction phenomenon to any 
dimension $d\geq 2$. 
We actually show   a stronger result: the projection of $(X_n)$ on each one dimensional line is locally contractive. This is not a direct consequence of \cite{BBE} since  the projected process is not a one-dimensional affine recursion in general.

Before stating our first result, let us fix some notations. Denote by $\overline{\mu}$ the distribution of $A_1$ on $\GL_d(\R)$ and   by $\Gamma_{\overline{\mu}}$ the closed semigroup in $\GL_d(\R)$  generated by its support. We assume $\Gamma_{\overline{\mu}}$ to be \emph{strongly irreducible} (it does not preserve any finite union of non-zero proper subspaces of $\R^d$) and to be \emph{proximal} (it contains a matrix with a simple dominant eigenvalue). These are  natural non-degeneracy assumptions  in the study of product of random matrices since 
Furstenberg \cite{Furstenberg1963} and Guivarc'h-Raugi \cite{GR}. We also  
denote by $\Gamma_\mu$ the closed semigroup in $\AffRd$ generated by the support of $\mu$. 
\begin{theo}\label{theorem.main.contraction}
	Let $d\geq 2$. Assume  that $$
	\E\Bigl[\bigl(\ln^+|A_1| + \ln^+|A_1^{-1}|
	+ \ln^+|B_1|\bigr)^2\Bigr] < +\infty$$ and 
    $\lambda_1 = 0$.  Assume also that $\Gamma_{\overline\mu}$ is strongly 
	irreducible   and 
	proximal  and that 
	$\Gamma_\mu$ has no common fixed point in $\R^d$. 
	
	Then, for every linear form $\tilde{v}\in(\R^d)^*$ and every 
	compact interval $I\subset\R$,
	\begin{equation}\label{eq.local.contraction.y}
		\P\!\left(\forall\, x,y\in\R^d,\;
		\tilde{v}(X_n^x - X_n^y)\,\mathds{1}_{\tilde{v}X_n^x\in I}
		\underset{n\to+\infty}{\longrightarrow} 0\right) = 1.
	\end{equation}
	In particular, for every compact set $K\subset\R^d$:
	\begin{equation}\label{eq.local.contraction.compact}
		\P\!\left(\forall\, x,y\in\R^d,\;
		|X_n^x - X_n^y|\,\mathds{1}_{X_n^x\in K}
		\underset{n\to+\infty}{\longrightarrow} 0\right) = 1.
	\end{equation}
\end{theo}

\paragraph{ Uniqueness of the stationary Radon measure.} Local contraction is of interest primarily when 
the processes $(X_n^x)$, $x\in \R^d$ are topologically recurrent, 
i.e.~there  exists a compact set $K$ of $\R^d$ such that for every $x\in \R^d$, $\P$-almost surely,  the process $(X_n^x)$ returns to $K$  infinitely 
often. 
Topological recurrence of the affine recursion \eqref{eq-rec.affine} is established in  our 
previous work \cite{AounBrofferioPeigne} under slightly stronger moment conditions.
Standard results \cite[Theorem 5.1]{Lin} ensure the existence of a  invariant Radon  measure on $\R^d$ for the Markov chain $(X_n^x)$. The local contraction property of our main result Theorem~\ref{theorem.main.contraction} in this paper  yields its uniqueness, see for instance \cite[Theorem 2.13]{PeigneWoess}. 
\begin{coro}\label{coro.uniqueness}
	Under the assumptions of Theorem~\ref{theorem.main.contraction}, 
	assume further that $\overline\mu$ has an exponential moment, i.e.\
	$\E(|A_1|^\tau + |A_1^{-1}|^\tau) < +\infty$ for some $\tau > 0$ and that $\E((\log^+\vert B_1\vert)^\gamma)<+\infty$ for some $\gamma>2$.
	Then, up to a positive multiplicative constant,  there exists a unique non-zero invariant Radon measure for 
	$(X_n^x)$ on $\R^d$. This measure has infinite mass.
\end{coro}

\paragraph{The crossing theorem for the joint process.}
The natural object encoding the relative behavior of two trajectories 
is the joint process $(X_n^x, X_n^y)$, obtained by the diagonal action 
of $\AffRd$ on $\R^d\times\R^d$. The key geometric phenomenon behind 
Theorem~\ref{theorem.main.contraction} is in fact   more 
precise  than \eqref{eq.local.contraction.y}: the projected segment $[\tilde{v}X_n^x, \tilde{v}X_n^y]$ 
eventually stops ``crossing" any fixed interval $I$ in any direction $\tilde{v} \in (\R^d)^*$ (i.e. $[\tilde{v}X_n^x, \tilde{v}X_n^y]$ does not contain $I$ when $n$ is large enough).   Here and throughout the paper, for $s,t\in\R$ we denote by
$[s,t]:=[\min\{s,t\},\max\{s,t\}]$ the closed interval with endpoints $s$
and $t$, regardless of their order.
This is the content of the following result, which is the theorem we 
actually prove. \begin{figure}
    \centering
    \figu{0.3}{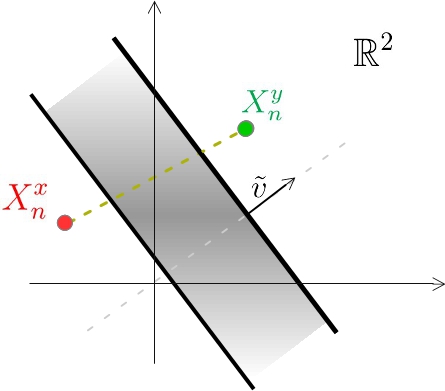} \qquad \qquad\figu{0.3}{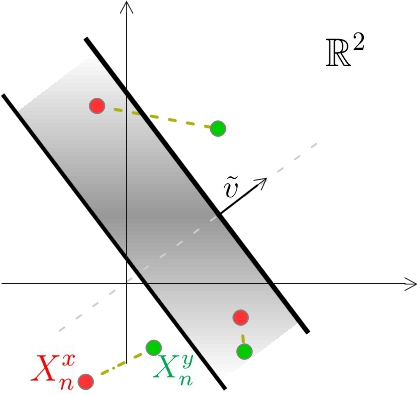}  
    \caption{Crossing / not crossing}
    \label{fig:Cross}
\end{figure}
\begin{theo}\label{theorem.main.crossing} 
	Under the assumptions of Theorem~\ref{theorem.main.contraction},
	for every non-zero $\tilde{v} \in (\R^d)^*$,
	\begin{equation}\label{eq.crossing}
		\P\Bigl(\forall\, x,y \in \R^d,\; \forall\, I \subset \R 
		,\;
		[ \tilde{v} X_n^x,  \tilde{v} X_n^y] \supset I
		\ \text{only finitely often}\Bigr) = 1,
	\end{equation}
	where $I$ ranges over all compact non-degenerate intervals of $\R$.
\end{theo}

\subsection{Context and prior results}

The study of contraction properties of   stochastic dynamical systems/iterated function systems as defined for instance in \cite{PeigneWoess}   has a rich history  at the interface of probability theory and dynamical systems.  In the setting of affine recursions on $\R^d$, 
these properties    depend on  the sign of the top Lyapunov 
exponent $\lambda_1$. The contracting regime $\lambda_1 < 0$ is by now well 
understood. The critical regime $\lambda_1 = 0$, 
which is the focus of the present paper, is more subtle: no stationary 
probability measure exists  under moment and non-degeneracy conditions, and the relevant notion of stability is local 
rather than global. We review the prior 
work  most directly related to our 
results, organized around these two regimes and the dimension of the ambient 
space.

\paragraph{The contracting regime and synchronization.}
When $\lambda_1 < 0$, the affine recursion is globally contracting and 
admits a unique invariant probability measure. Bougerol-Picard \cite{Bougerol-Picard} 
showed that $\lambda_1 < 0$ is also necessary for the existence of 
a finite $\mu$-invariant measure.
This global contraction phenomenon has a parallel history 
in the theory of random dynamical systems on other metric spaces, mainly compact ones.  
Under the name of \emph{synchronization},  this phenomenon   refers  to the convergence to $0$ 
of the distance between two trajectories starting from different initial 
conditions. In the linear setting,     going back to Furstenberg~\cite{Furstenberg1963} and 
Guivarc'h-Raugi~\cite{GR}, 
the action on  the projective space $\mathbf{P}(\R^d)$ of a left random walk on  $\GL_d(\R)$ is contracting under the strong irreducibility and 
proximality assumptions. More precisely,    the corresponding transition operator  on $\mathbf{P}(\R^d)$   has a 
spectral gap on a suitable functional space (see Le Page \cite{page}) and  it admits a unique invariant measure. On the circle $\mathbf{S}^1$, 
Antonov~\cite{Antonov} showed, under a minimality assumption, that random 
compositions of homeomorphisms of $\mathbf{S}^1$ synchronize; this is equivalent to 
the absence of a common invariant probability measure. This result was 
rediscovered   by Kleptsyn--Nalskii~\cite{KleptsynNalskii} under 
  a weaker assumption and extended 
by Malicet~\cite{Malicet}. Synchronization rates and limit laws for random dynamical systems on 
abstract compact metric spaces have been studied in \cite{GelfertSalcedo} 
and references therein. 
In the physics literature, a related but distinct 
phenomenon is studied for deterministic dynamical systems governed by 
differential equations: synchronization refers to the convergence of two 
trajectories coupled    by the same external forcing, and does not involve randomness. 
We refer to~\cite{PecoraCarroll} for a description of this phenomenon and 
further examples.

\paragraph{The critical regime in dimension $d=1$.}
Local contraction for one-dimensional affine recursions was established by 
Babillot-Bougerol-\'Elie~\cite{BBE} and subsequently extended in several 
directions, most remaining one-dimensional.
Brofferio~\cite{Brofferio-how} gave a new proof via a potential-theoretic 
approach, establishing local contraction as a consequence of the convergence 
of the right random walk to infinity; this approach directly inspired the 
present work. 

Deroin, Kleptsyn, Navas and Parwani~\cite{DKNP} established uniqueness of 
the stationary Radon measure for random walks on $\mathrm{Homeo}_+(\R)$ 
under symmetry and compact support assumptions, via a approach specific to 
the one-dimensional ordered setting, showing the robustness of the phenomenon 
beyond the affine case. In the same setting, Brofferio, Buraczewski 
and Szarek~\cite{BBS} subsequently established uniqueness of the stationary 
Radon measure under the more general conditions of recurrence, contraction 
and unbounded action, without symmetry assumptions.

Brofferio--Buraczewski~\cite{Brofferio-Buraczewski} established recurrence 
and local contraction for stochastic dynamical systems on the half line $\R_+$, noting 
in particular that local contraction follows rather directly for systems 
defined by not necessarily invertible maps, where a reflection at the 
origin provides a natural contraction mechanism.

\paragraph{Extensions and abstract frameworks.}
Cartwright, Kaimanovich 
and Woess~\cite{CKW} and Brofferio~\cite{brofferio.thesis} studied random 
walks on the affine group of homogeneous trees   in the critical regime and established local 
contraction. Such settings  carry a geometry that differs substantially from the real case as these spaces are totally disconnected and carry no natural order. 

Benda~\cite{benda} and Peign\'{e}-Woess~\cite{PeigneWoess,PeigneWoessII} 
developed a general framework for stochastic dynamical systems on metric 
spaces, studying the ergodic consequences of local contraction for systems 
generated by continuous maps.

\paragraph{The critical regime in dimension $d\geq 2$.}
For affine recursions on $\R^d$, Brofferio, Peign\'e and Pham~\cite{BrofferioPeignePham} 
treated the case of matrices with nonnegative entries and nonnegative 
translation vectors satisfying a strong uniform positivity condition. Their work was inspired by \cite{Brofferio-Buraczewski} and the   contraction of a cone in the sense of Birkhoff plays an essential role. 
Local contraction in the general setting of dimension $d \geq 2$ has remained 
open until now.

The present work also contributes to the broader program of classifying 
stationary Radon measures of infinite mass on homogeneous spaces $G/H$. 
In our setting, $G=\Aff(\R^d)$ is a non-reductive group acting on the 
algebraic homogeneous space $\R^d\simeq G/\GL_d(\R)$. 
Guivarc'h--Raja~\cite{GR2} developed an abstract ergodicity framework for 
random walks on homogeneous spaces; Bru\`{e}re~\cite{Bruere} established 
recurrence criteria when $G$ is semisimple and $H$ algebraic; 
B\'enard~\cite{Benard} classified stationary Radon measures on 
$\mathbb{T}^d\times\R\simeq G/\Gamma$ with $G$ non-reductive and $\Gamma$ 
a lattice. Closest to our setting, Benoist--Bru\`{e}re~\cite{Benoist-Bruere} 
characterized existence of stationary probability measures on all $k$-affine 
Grassmannians $\mathcal{A}_k(\R^d)$, generalizing 
Bougerol--Picard~\cite{Bougerol-Picard} to any $k$.
Later, Aoun--Sert \cite{aoun-sert1,aoun-sert2} classified  stationary probability measures on projective spaces,  a setting that encompasses both affine spaces and  affine Grassmannians.  Our companion 
paper~\cite{AounBrofferioPeigne} established  the existence of  a non-zero stationary 
Radon measure on $\R^d=\mathcal{A}_0(\R^d)$ (of infinite mass); the present 
work focuses on its uniqueness  up to positive multiplicative constant in the critical regime $\lambda_1=0$.
\subsection{Strategy of the proof and organisation of the paper}

The proof of Theorem~\ref{theorem.main.contraction} rests on two main ideas: 
a geometric reformulation of local contraction in terms of the diagonal 
process, and a potential-theoretic argument based on total dissipativity 
of the Markov operator.
\paragraph{A geometric point of view.} 
We study the joint process $(X_n^x, X_n^y)$ defined by the random diagonal action of 
$\AffRd$ on
\[
\cZ_d := (\R^d \times \R^d) \setminus \Delta,
\]
where $\Delta$ denotes the diagonal. As mentioned as the end of Section \ref{sec.1.1}, Theorem \ref{theorem.main.contraction} follows from Theorem \ref{theorem.main.crossing}. It therefore suffices to show that for every non-zero linear form 
$\tilde{v}\in(\R^d)^*$, 
the 
segment $[\tilde{v}X_n^x, \tilde{v}X_n^y]$ can contain $I$ only finitely 
many times;
 in other words, that the diagonal process eventually 
leaves each \emph{crossing set}
\begin{equation}\label{eq.def.crossing.set}
V_{\tilde{v},I} := \bigl\{(x,y)\in\cZ_d \mid [\tilde{v}x,\tilde{v}y]\supset I\bigr\}.
\end{equation}
These crossing sets are well-adapted to the affine structure and monotone 
under enlargement of segments: $[x,y]\subset[x',y']$ implies 
$(x,y)\in V_{\tilde{v},I}\Rightarrow(x',y')\in V_{\tilde{v},I}$. 
This monotonicity, combined with a countable intersection argument, will 
allow us to transfer the crossing property from almost every pair to every 
pair in the final step of the proof.

\paragraph{Counting entrances.}
Fix $\tilde{v}$ and $I$. Rather than estimating the total number of visits 
to $V_{\tilde{v},I}$, whose expected value diverges in our setting, we 
count the expected number of \emph{entrances} into $V_{\tilde{v},I}$, that is
\[
\sum_{n\geq 0}
\P\bigl((X_n^x,X_n^y)\notin V_{\tilde{v},I},\;
(X_{n+1}^x,X_{n+1}^y)\in V_{\tilde{v},I}\bigr).
\]
By Borel--Cantelli, finiteness of this series implies that the process 
re-enters $V_{\tilde{v},I}$ only finitely many times. 
Hence, after the last 
 entrance  time into $V_{\tilde v, I}$, the process either remains trapped in $V_{\tilde{v},I}$ forever, 
or it exits $V_{\tilde v, I}$ and never re-enters.
 Under our hypotheses, the first case is ruled out by the 
fact that $\liminf_{n\to\infty}|A_n\cdots A_1(x-y)|=0$ almost surely, forcing the projected segment 
$[\tilde{v}X_n^x,\tilde{v}X_n^y]$ to eventually be too short to contain $I$.

\paragraph{A potential-theoretic approach.}
Let $P_\mu$ be the Markov operator of the diagonal process $((X_n^x,X_n^y))$ on $\cZ_d$. Setting  
\[
\phi_{\tilde{v},I}(x,y) :=
\mathds{1}_{(x,y)\notin V_{\tilde{v},I}}\,
\P\bigl((X_1^x,X_1^y)\in V_{\tilde{v},I}\bigr),
\]
the expected number of entrances into $V_{\tilde v, I}$    is equal to   
$\sum_{n\geq 0} P_\mu^n\phi_{\tilde{v},I}(x,y)$  ,  that  is the potential of   $P_\mu$ evaluated at the function $\phi_{\tilde v, I}$. 
If $P_\mu$ is \emph{totally dissipative} on $L^1(\cZ_d,\sigma)$ for a 
suitable $\mu$-invariant $\sigma$-finite measure $\sigma$ on $\cZ_d$, 
and if $\phi_{\tilde{v},I}\in L^1(\cZ_d,\sigma)$, then this potential 
is finite $\sigma$-almost everywhere, giving the crossing 
property for $\sigma$-almost every starting pair $(x,y)$.
 The  three steps:  (a) construction of $\sigma$, (b) integrability of $\phi_{\tilde{v},I}$ 
and (c) total dissipativity of 
$P_\mu$, are the main 
building blocks of the paper. The extension from $\sigma$-almost every 
pair to every pair $(x,y)\in\cZ_d$ is then achieved via a countable 
intersection argument, using Oseledets' theorem and the monotonicity of the crossing sets 
$V_{\tilde{v},I}$ under enlargement of segments.

\paragraph{New contributions with respect to dimension $d=1$.}
The strategy of the proof is partly inspired by the potential-theoretic 
approach of Brofferio~\cite{Brofferio-how} in dimension $d=1$. In that 
setting, $\Aff(\R)$ acts simply and transitively on $\cZ_1$, allowing 
one to identify $\cZ_1\simeq\Aff(\R)$. As a consequence, the left Haar 
measure on $\Aff(\R)$ provides directly the $\mu$-invariant $\sigma$-finite 
measure $\sigma$ on $\cZ_1$.  Integrability of $\phi_{\tilde{v},I}$ is a simple computation, and the total dissipativity of $P_\mu$ follows 
from the non-unimodularity of $\Aff(\R)$ via the theorem of 
Guivarc'h-Keane-Roynette~\cite{Guivarch-Keane-Roynette}.

In dimension $d\geq 2$, each of these steps requires a genuinely new 
approach. 
\begin{enumerate}[label=(\alph*)]
	\item  It can be verified  that there is no $\Aff(\R^d)$-invariant non-zero Radon measure on $\cZ_d$.  We construct a $\mu$-invariant Radon measure $\sigma$ on $\cZ_d$  
	 as a fiber integral over the
 unique $\mu$-invariant probability measure on the   
	Grassmannian $\AG$ of affine lines of $\R^d$, using left Haar conditional measures on the fibers 
	$\cZ_1(\ell)\simeq\Aff(\R)$.
	\item The integrability of $\phi_{\tilde{v},I}$ involves a new transverse 
	component: the   \emph{height} of the stationary line 
	on $\AG$, that is its distance to the origin. We show that this component is log-integrable,  this has no analogue in dimension $d=1$.
	\item The total dissipativity can no longer be deduced from non-unimodularity 
	alone, since the action of $\Aff(\R^d)$ on $\cZ_d$ is 
    not free. We develop an abstract dissipativity criterion 
that extends  the  theorem of Guivarc'h-Keane-Roynette 
	 beyond the group framework to semigroup actions admitting a suitable 
	commuting transformation.
\end{enumerate}

Let us comment on the assumption $d\geq2$ in our main results. It is
only used to guarantee, together with {\bf I} and {\bf P}, that $\Gamma_\mu$
contains an element whose linear part is proximal with dominant eigenvalue of
modulus different from $1$ ; this is the
key input for the dissipativity criterion. In dimension one, {\bf I} and
{\bf P} are empty conditions, and one assumes instead directly that
$\P(|A_1|=1)<1$, as in \cite{BBE,Brofferio-how}.

\paragraph{Organisation of the paper.}
Section~\ref{sec.preliminaries}  gathers known results 
on random linear and affine actions on $\R^d$. It also  introduces the main 
geometric objects we use in the paper, namely the   space $\cZ_d$ and its fibration over the affine 
Grassmannian $\AG$.
Section~\ref{section.construction.sigma} constructs $\sigma$ and 
establishes the log-integrability of  the height of  the stationary line on $\AG$.
Section~\ref{sec.end.proof} proves integrability of $\phi_{\tilde{v},I}$, 
and derives Theorem~\ref{theorem.main.crossing} and 
Theorem~\ref{theorem.main.contraction} assuming that $P_\mu$ is totally dissipative.
Section~\ref{section.dissipative} establishes  the  total dissipativity of 
$P_\mu$ on $L^1(\cZ_d,\sigma)$: Section \ref{the.criterion} develops the
 abstract criterion 
(Theorem~\ref{theorem.main.dissipativity})  mentioned above, and 
Section~\ref{sec.affine.application} applies it to the diagonal action.

\section*{Acknowledgments}
The authors would like to warmly thank Philippe Bougerol, to whom this work owes a lot,  for all his encouragement during the long years. They also would like to thank Fran\c{c}ois Ledrappier and Wolfgang Woess for their interest through all these years and 
 Cagri Sert for helpful discussions.

\section*{AI disclosure}
This paper was born in  Halong Bay  in 2017 during a mini-course on stochastic dynamical systems. The mathematics and a first draft of the manuscript were  completed in March 2026 and obtained without any AI assistance. AI was then used recently to assist with language editing and  improve the clarity of some formulations. 
\section{Setting and preliminary results}
\label{sec.preliminaries}

Throughout the paper, $\mu$ is a probability measure on $\AffRd$ and
$(\gamma_k)_{k\in\Z}=\bigl((B_k,A_k)\bigr)_{k\in\Z}$ is a bi-infinite sequence
of i.i.d.\ $\AffRd$-valued random variables with distribution $\mu$, defined on
a probability space $(\Omega,\mathcal F,\P)$. We denote by $\overline\mu$ the
distribution of the linear part $A_k$; thus $(A_k)_{k\in\Z}$ is an i.i.d.\
sequence of $\GL_d(\R)$-valued random variables with distribution
$\overline\mu$.

The affine recursion~\eqref{eq-rec.affine},
$X_n^x=\gamma_n\cdots\gamma_1 x$, is obtained by iterating the random affine
maps $\gamma_n$ with positive indices $n\geq 1$. 
The  backward products
$A_0\cdots A_{-n+1}$ and $\gamma_0\cdots\gamma_{-n+1}$, for $n\geq 1$, are used to construct the stationary objects
we rely on:
the Furstenberg process in Section~\ref{projectiveaction} and the
stationary affine line in Section~\ref{sec.line}.
When needed, we work on the canonical space
$\Omega=(\AffRd)^{\Z}$ with $\P=\mu^{\otimes\Z}$, endowed with the shift
$\theta\bigl((\gamma_k)_{k\in\Z}\bigr):=(\gamma_{k+1})_{k\in\Z}$.
\subsection{Linear actions }

\subsubsection{Random linear action on $\R^d$}

We fix $d\geq 1$ and a Euclidean structure on $\mathbb R^d$, i.e.\ a scalar
product $\langle\cdot,\cdot\rangle$ and the associated norm $\vert\cdot\vert$.
Elements of the vector space $\R^d$ are denoted with letters $v, w$ and are
considered as column vectors. The space $\R^d$ is also endowed with an affine
structure; in this case its elements are denoted with letters $x, y, z$.

We write $(\R^d)^*$ for the dual space of $\R^d$ endowed with the corresponding
Euclidean structure. Its elements are denoted $\tilde{v},\tilde{w}$ and are
considered as row vectors. For short, when $\tilde v\in (\R^d)^*$ and
$w\in \R^d$, we write $\tilde v(w)=\tilde v w$.

The group $\GL_d(\R)$ of $d\times d$ invertible matrices acts naturally on the
left on $\R^d$ and on the right on $(\R^d)^*$; these actions satisfy
the associativity relation $(\tilde v A) w=\tilde v (Aw)$, which we simply
write $\tilde v A w$.
We endow $\GL_d(\R)$ with the operator norm
$|A|:=\sup\{|A v|\mid |v|=1\}=\sup\{|\tilde v A | \mid |\tilde v|=1\}$.

We also use the exterior product $\wedge^2\R^d$ of $\R^d$, endowed with the
Euclidean structure inherited from $\R^d$; we still denote its norm
by $|\cdot|$. Thus $|v\wedge w|=|v|\,|w|\,|\sin\angle(v,w)|$ is the area of the
parallelogram spanned by $v$ and $w$.

The linear action of an element $A\in \GL_d(\R)$ on $\R^d$ induces a
linear action on $\wedge^2\R^d$ defined by $(\wedge^2A)(v\wedge w)= Av\wedge Aw$
for any $v, w$ in $\R^d$. We denote by
$|\wedge^2 A|$ the corresponding operator norm.

For $p>0$, we will need the following classical
$\ln^p$-integrability assumption.

\vspace{3mm}

\noindent {\bf Hypothesis ${\rm \bf M}_p(\bar \mu)$}: {\it $\mathbb E[ (\ln \max\{
	\vert A_1\vert, \vert A_1^{-1}\vert\})^p]<+\infty$.}

\vspace{3mm}

Under hypothesis ${\rm M}_1(\bar \mu)$, the \emph{first Lyapunov exponent}
$\lambda_1=\lambda_1(\bar \mu)$ is well defined in $\R$ as
\begin{equation} \label{es-def-lambda1}
	\lambda_1:= \lim_{n \to +\infty} \frac{\E(\ln \vert A_n\cdots A_1\vert)}{n}.
\end{equation}
Since $(A_k)_{k\geq 1}$ and $(A_k)_{k\leq 0}$ have the same distribution,
we also have
\[
\lambda_1=\lim_{n \to +\infty} \frac{\E(\ln \vert A_0\cdots A_{-n+1}\vert)}{n}.
\]

By Kingman's subadditive ergodic theorem \cite{kingman} (the original
result is due to Furstenberg--Kesten \cite[Theorem 2]{FurstenbergKesten}), the
convergence \eqref{es-def-lambda1} also holds almost surely, for both the
forward and the backward products:
\begin{equation}\label{eq-lambda1-ps}
\lim_{n\to+\infty}\frac{\ln|A_n\cdots A_1|}{n}
=\lim_{n\to+\infty}\frac{\ln|A_0\cdots A_{-n+1}|}{n}
=\lambda_1 \qquad \P\text{-a.s.}
\end{equation}

Similarly, under ${\rm M}_1(\bar \mu)$, the \emph{second Lyapunov exponent} $\lambda_2$ is  defined   via the limit of the action of the  products $A_n \cdots A_1$ on the exterior product $\wedge^2 \R^d$. More precisely,
\begin{align}
\label{eq-def-lambda2}
    \lambda_1+\lambda_2
	&:= \lim_{n \to +\infty} \frac{\E(\ln \vert \wedge^2 (A_n\cdots A_1)\vert)}{n}
	\\
	 &=\nonumber 
	\lim_{n \to +\infty} \frac{\ln | \wedge^2(A_n\cdots A_1)|}{n}=\lim_{n \to +\infty} \frac{ \ln |\wedge^2(A_0\cdots A_{-(n-1)})|}{n} \quad \P{\rm -a.s.}    
\end{align}

 \subsubsection{Random projective action}
\label{projectiveaction}

Let us now describe the projective action of $\GL_d(\R)$.
The projective space $\PRd$ is the quotient space $(\R^d\setminus\{0\})/{\sim}$,
where $v \sim w$ if and only if $\R v = \R w$, for $v, w \in \R^d\setminus\{0\}$;
its elements are denoted $[v] := \R v$.
We endow $\PRd$ with the Fubini--Study metric
$\delta([v],[w]) := \frac{|v\wedge w|}{|v|\,|w|} = |\sin\angle(v,w)|$.
The group $\GL_d(\R)$ acts on $\PRd$ by
$$
\forall A\in \GL_d(\R),\ \forall [v]\in \PRd,
\qquad A \cdot [v] := [Av] = \R Av.
$$
For each $V_0:=[v]\in \PRd$, the Markov chain
$(V_n )_{n\geq 0}$ with
$V_n := A_n \cdots A_1\cdot V_0$ for $n\geq 1$,
is called the \emph{Furstenberg process}.
Since $\GL_d(\R)$ acts continuously on the compact space $\PRd$,
there exists a \emph{$\overline{\mu}$-invariant probability measure} $\nu$ on $\PRd$,
i.e.\ an invariant  probability measure for the Markov chain $(V_n)_{n\geq 0}$. In other words, for every bounded measurable function
$f\colon \PRd \to \R$,
$$
\overline{\mu} * \nu(f)
:= \int_{\GL_d(\R)}\int_{\PRd} f(A\cdot [v])\,d\overline\mu(A)\,d\nu([v])
= \nu(f).
$$
Equivalently,
$\E\bigl(f(A_1\cdot V_0)\bigr) = \E\bigl(f(V_0)\bigr)$,
where $V_0$ is a random element of $\PRd$ with distribution $\nu$,
independent of $(A_n)_{n\geq 1}$.

The following classical assumptions guarantee that the random action is contractive on $\PRd$.
\vspace{3mm}

\noindent{\bf Hypothesis I}: \textit{The semigroup $\Gamma_{\overline\mu}$ is
	strongly irreducible, i.e.\ no finite union of proper subspaces of $\R^d$ is
	invariant under $\Gamma_{\overline\mu}$.}

\vspace{3mm}

\noindent{\bf Hypothesis P}: \textit{The semigroup $\Gamma_{\overline\mu}$ is
	proximal, i.e.\ it contains a proximal matrix, that is, a matrix with a
	simple dominant eigenvalue.}

\vspace{3mm}

We recall some key properties of the Furstenberg process that hold
under hypotheses ${\bf M}_1(\overline\mu)$, {\bf I} and {\bf P}. 
These are consequences of landmark results of Guivarc'h--Raugi \cite{GR};
we refer to \cite[Chapter~III]{Bougerol-Lacroix} and
\cite[Chapters~III and~IV]{BQ-book} for further details.

The $\overline{\mu}$-invariant probability measure $\nu$ on $\PRd$ is unique
and \emph{non-degenerate} (i.e. $\nu$ assigns zero mass to every projective
hyperplane).

The first Lyapunov exponent is simple, i.e.\ $\lambda_1 > \lambda_2$,
and the Furstenberg process synchronizes at exponential speed:
for every $[v],[w]\in\PRd$,
$$
\lim_{n\to+\infty} \frac{1}{n}
\ln \delta\bigl(A_n\cdots A_1\cdot [v],\, A_n\cdots A_1\cdot [w]\bigr)
\leq  \lambda_2 - \lambda_1 < 0
\quad \text{a.s.}
$$

All cluster points of the normalized sequence
$(A_{0}\cdots A_{-n+1}/|A_{0}\cdots A_{-n+1}|)_{n\geq 1}$
are $\P$-a.s.\ rank-one matrices sharing the same random image $[\vv_0]$.
The distribution of the random direction $[\vv_0]$ is $\nu$.

For each $k\in\Z$, we set $[\vv_k] := [\vv_0] \circ \theta^k$, i.e.\ $[\vv_k]$
is the limiting direction of the shifted sequence
$(A_{k}\cdots A_{k-n}/|A_{k}\cdots A_{k-n}|)_{n>0}$.
The sequence $([\vv_k])_{k\in\Z}$ is a stationary bi-infinite sequence,
equivariant under the action of $(A_k)_{k\in\Z}$:
\begin{equation}\label{eq.equivarance.vk}
	[\vv_k] = A_k\cdot[\vv_{k-1}] \qquad \forall\, k\in\Z.
\end{equation}
Furstenberg's formula gives
\begin{equation}\label{eq-furstemberg-lambda1}
	\lambda_1
	= \E\ln\vert A_1\frac{\vv_0}{\vert\vv_0\vert}\vert
	= \iint \ln\frac{|Av|}{|v|}\,d\overline\mu(A)\,d\nu([v]).
\end{equation}

The dual projective space $\mathbf{P}((\R^d)^*)$ is the space of lines in
$(\R^d)^*$. 
The group $\GL_d(\R)$ acts on $\mathbf{P}((\R^d)^*)$ on the right: for $[\tilde v]\in \mathbf{P}((\R^d)^*)$ and $A\in \GL_d(\R)$, 
$$
[\tilde{v}]\cdot A := [\tilde{v}A].
$$
This dual right action is equivalent to the left action of $\GL_d(\R)$ on
$\PRd$ via transposition. In particular, Hypotheses {\bf I} and {\bf P} for
$\overline\mu$ imply the analogous hypotheses for the transpose  and  all the results above 
apply to this dual action.

Under hypotheses ${\bf M}_1(\overline\mu)$, {\bf I} and {\bf P},
all cluster points of the normalized sequence
$(A_n\cdots A_1/|A_n\cdots A_1|)_{n\geq 1}$
are $\P$-a.s.\ rank-one matrices sharing the same random kernel $\ker\vvt_1$.
The random direction $[\vvt_1]\in\mathbf{P}((\R^d)^*)$ is measurable with
respect to $(A_j)_{j\geq 1}$, hence independent of $(A_j)_{j\leq 0}$.
Its distribution $\tilde\nu$ is the unique right $\overline\mu$-invariant
probability measure on $\mathbf{P}((\R^d)^*)$, and it is non-degenerate.
In particular, $\vvt_1\vv_0\neq 0$ almost surely.

Under hypotheses {\bf I}, {\bf P} and ${\bf M}_1(\overline\mu)$,
Oseledets' theorem  gives, almost surely, for every
$v\in\R^d\setminus\{0\}$:

\begin{eqnarray}\label{eq-osseledats}
	\lim_{n\to+\infty}\frac{1}{n}\ln|A_n\cdots A_1 v| &=& \lambda_1
	\quad \text{if } \vvt_1 v \neq 0,
	\\\
	\limsup_{n\to+\infty}\frac{1}{n}\ln|A_n\cdots A_1 v| &\leq& \lambda_2
	\quad \text{if } \vvt_1 v = 0.
\end{eqnarray}

\subsection{Affine actions}

\subsubsection{Random affine recursion on $\R^d$}
The group $\Aff(\R^d)$ of invertible affine transformations $\gamma: x\mapsto A x+B$ of  $\R^d$, endowed with the composition operation,  
coincides with the semi-direct product $\R^d\rtimes \GL_d(\R)$. For now, for $\gamma\in \Aff(\R^d)$, we write $\gamma=(B,A)$ with $B\in \R^d$ and $A\in \GL_d(\R)$. The group law is $(B',A') (B,A)=(B'+A'B,A'A)$.

Let $\mu$ be a probability measure on $\AffRd$ and $(\gamma_n = (B_n, A_n))_{n \in \Z}$ a bi-infinite sequence of i.i.d.~random variables with common distribution $\mu$.
We consider the stochastic affine recursion  on $\R^d$ defined by
\[
X_{n+1}^x = \gamma_{n+1} X_n^x= A_{n+1} X_n^x + B_{n+1},\qquad X_0^x=x.
\]
For $p>0$ ,  we   introduce the following moment assumption that extends 
 $\ln^p$-integrability to $\AffRd$:

\vspace{3mm}

 \noindent {\bf Hypothesis  ${\rm \bf M}_p( \mu)$} : 
${\rm \bf M}_p( \overline \mu)$ {\it  holds and} 
 $\E((\ln^+\vert B_1\vert)^p)<+\infty$. 

 \vspace{3mm}

\subsubsection{Action on affine lines $\AG$}

We now consider the Grassmannian $\AG$ of affine lines in $\R^d$, that is,
\[
\AG := \{ \ell=x + \R v \mid x \in \R^d,\ v \in \R^d\setminus\{0\} \}.
\]
The group $\AffRd$ acts on $\AG$ by: for any $\gamma=(B,A) \in \AffRd$ and any  $\ell = x +\R v\in \AG$,
\[\gamma \cdot \ell = Ax + B + \R Av.\]

We assume that $\Gamma_{\overline \mu}$ satisfies assumptions \textbf{I} and \textbf{P} and that $\lambda_1=0$.  
By \cite{Bougerol-Picard} there is  no $\mu$-invariant probability measure on $\R^d$. Nevertheless, the fact that $\lambda_2<0$ gives sufficient contraction properties on $\AG$. 
Along with hypothesis ${\rm \bf M}_1( \mu)$, this  insures
the existence of a $\mu$-invariant probability measure $\eta$ on $\AG$, meaning that for any measurable function $f : \AG \to \R^+$, one has 
$$\mu * \eta(f) := \E\bigl(f(\gamma_1\cdot \ell_0)\bigr)
= \E\bigl(f(\ell_0)\bigr) = \eta(f),$$
where $\ell_0$ denotes a random affine line with law $\eta$, independent of the steps $\{(B_n, A_n)\}_{n \geq 1}$.

The existence of $\mu$-invariant probability measures on the  $k$-dimensional affine Grassmannian $\mathcal A_k(\mathbb R^d)$  was studied by Benoist--Bru\`{e}re~\cite{Benoist-Bruere}. Under compact support assumptions and hypotheses {\bf I} and {\bf P}, they showed that such a measure exists if and only if  the $k^{th}$ Lyapunov exponent of the sequence $(A_n)_n$  is negative. Thus, under this strong assumption on the support of $\mu$, the existence of $\eta$ follows as a particular case of \cite[Theorem 1.3]{Benoist-Bruere} for $k=1$. 

In   Section \ref{section.construction.sigma}, we provide an explicit construction of $\eta$ under first moment assumption only (Proposition~\ref{prop.line.integrable}). This construction  also allows us to derive moment estimates for $\eta$.

Note that the natural projection from $ \AG$ to $\PRd$, which sends $\ell = x + \R v$ to $[v]$, intertwines the corresponding actions of $\AffRd$ and $\GL_d(\R)$. Consequently, the pushforward of $\eta$ under this projection equals $\nu$.

\subsection{Diagonal affine action}
\label{sec.diag.action}
In order to study the synchronization between   two trajectories $(X_n^x)_n$ and $(X_n^y)_n$, 
it is natural to consider  the space $\cZ_d$ of pairs of 
distinct points of $\R^d$, defined by 
$$\cZ_d := \{ (x,y)\in \R^d \times \R^d \mid x\neq y\} = (\R^d \times \R^d)\setminus \Delta,$$
where $\Delta = \{(x,x) \mid x \in \R^d\}$ denotes the diagonal of $\R^d\times \R^d$. We then consider the  diagonal action   of  $\AffRd$ on $\R^d\times \R^d$, given  by : for all $\gamma \in \AffRd$  and all $x, y\in \R^d$,
\[
\gamma (x,y) := (\gamma x, \gamma y).
\]
The space $\cZ_d $ is $\Aff(\R^d)$-invariant and we restrict the action to it. 
An important feature of this space is that it fibers over the affine 
	Grassmannian   $\AG$ via the map $p$ that sends   each couple $(x,y)$ of distinct points of $
\R^d$ to the unique line passing through them,   \[
p : \left\{\begin{array}{clc}
\cZ_d   &\to &\AG
\\(x,y) &\mapsto& x + \R (y-x).
\end{array}\right.
\]
This is an $\Aff(\R^d)$-equivariant map. For each $\ell\in \AG$, the fiber $p^{-1}(\ell)$ above the line $\ell$, which we denote in the sequel by $\cZ_1(\ell)$, coincides with the set of pairs of distinct points on $\ell$: 
$$\mathcal{Z}_1(\ell):=p^{-1}(\ell)=(\ell\times \ell)\setminus \Delta \subset \cZ_d.$$
Once a base point $x_\ell$ and a direction vector $v_\ell$ of the line $\ell$ are chosen, $\R$ is 
naturally identified with $\ell$   via the map $\varphi_\ell:t\mapsto x_\ell+t v_\ell$. This allows to identify $\cZ_1=(\R\times \R)\setminus \Delta$ with $\cZ_1(\ell)$ via the bijection: 
\begin{equation}\label{phiell.initial}
\varphi_\ell : \left\{\begin{array}{clc}
\cZ_1    &\to &\cZ_1(\ell)
\\
(s,t) &\mapsto& (x_\ell+sv_\ell, x_\ell +t v_\ell).
\end{array}
\right.
\end{equation}
The fiber structure  of $\cZ_d$ with fibers
$\simeq \cZ_1$ is central to our 
	construction: it reduces the problem to the one-dimensional case, where 
	$\Aff(\R)$ acts simply transitively on $\cZ_1$ and the left Haar measure 
	provides a natural invariant Radon measure on each fiber $\cZ_1(\ell)$.

Finally, throughout the paper, we make use of a second system of coordinates on $\cZ_d$ which separates the affine action and the linear action. More precisely, $\cZ_d$ can be identified with $\R^d\times (\R^d\setminus \{0\})$ via 
the bijection 
\[
\begin{array}{clc}
\cZ_d
&\longrightarrow&
\R^d\times(\R^d\setminus\{0\})
\\
(x,y)
&\longmapsto&
(x,v):=(x,y-x)
\end{array}
\] We call $(x,y)$ the \emph{pair coordinates} and $(x,v)$ the
\emph{point--vector coordinates} on $\cZ_d$.
In  the point-vector coordinates, the action $\gamma=(B,A)\in \Aff(\R^d)$ on $\cZ_d$ is   given by: for $(x,v)\in \R^d\times (\R^d\setminus \{0\})$,
\begin{equation}\label{eq.action.Ab}
\gamma(x,v)=(Ax+B, Av). \end{equation}
 The projection $p$ introduced above then takes the form
\[
p(x,v)=x+\R v.
\]
In dimension one, these coordinates indentify $\cZ_1$ with $\R\times (\R\setminus\{0\})$ via the map 
$(s,t)\longmapsto (b,a):=(s,t-s)$.
Following our conventions, we use the letters $(s,t)$ for the pair
coordinates and $(b,a)$ for the point--vector coordinates on $\cZ_1$. 
In particular, for every $\ell=x_
\ell+\R v_\ell\in\AG$,
 the map $\varphi_\ell$
defined in \eqref{phiell.initial} takes the form 
\begin{equation}\label{phiell}
\varphi_\ell:
\left\{
\begin{array}{clc}
\R\times(\R\setminus\{0\})\simeq \cZ_1 
&\longrightarrow&
\cZ_1(\ell)
\\
(b,a)
&\longmapsto&
(x_\ell+bv_\ell,av_\ell).
\end{array}
\right.
\end{equation}
Depending on the context, we   work in $\cZ_d$ either in pair coordinates $(x,y)$ or in point-vector coordinates $(x,v)$,  in which case $\cZ_d$ is  identified     with $\R^d\times (\R^d\setminus \{0\})$.

\section{Construction of $\mu$-invariant measures}
\label{section.construction.sigma}
In dimension $d=1$,  the affine group $\Aff(\R)$ acts simply transitively on $\cZ_1$,  
so the  Haar measure on $\Aff(\R)$ provides a $\Aff(\R)$-invariant Radon measure on $\cZ_1$; this measure  is obviously $\mu$-invariant. In   dimension $d\geq 2$, the action of $\Aff(\R^d)$ on $\cZ_d$ is   transitive but is no more free; so that $\cZ_d$ is only an $\Aff(\R^d)$-homogeneous space. In fact, one can show that there is no $\Aff(\R^d)$-invariant non-zero Radon measure on $\cZ_d$.\footnote{The stabilizer of $(0,e_1)$ is the subgroup $H$ of affine transformations fixing the origin and with linear part of the form $A=\begin{bmatrix} 1 & \tilde w\\
0 & C\end{bmatrix}$, $\tilde w\in (\R^{d-1})^*$ and $C\in \GL_{d-1}(\R)$. 
Denote by $\Delta_G$ and $\Delta_H$ the respective modular functions of $G:=\Aff(\R^d)$ and $H=\textrm{Stab}_{(0,e_1)}$.  
A direct computation shows that $\Delta_G(B,A)=|\det A|^{-1}$ and that $\Delta_H(0,A)=|\det A|$. Hence 
  $(\Delta_G)_{|_H}\neq \Delta_H$. Using \cite[Theorem 2.51]{Folland}, it follows that no non-zero invariant Radon   measure on $G/H\simeq \cZ_d$ exists.}
To overcome this difficulty, we construct in this section a $\mu$-invariant Radon measure on $\cZ_d$. We proceed in two steps. We first construct a $\mu$-invariant probability measure $\eta$ on the affine Grassmannian $\AG$.  Then  we lift $\eta$ to a  Radon $\mu$-invariant measure $\sigma$ on $\cZ_d$ by disintegration along the fibers of the natural $\AffRd$-equivariant projection $p$ onto $\AG$.

\subsection{The $\mu$-invariant probability measure $\eta$ on $\AG$}
\label{sec.line}
For $[v] \in \PRd$  and $x\in \R^d$, we denote by
	$$
	\Pr{ [v]}x := x - \frac{\langle x,v\rangle}{|v|^2}v
	$$
	the orthogonal projection onto the hyperplane $(\R v)^\perp$ (observe that it only depends on the class $[v]\in\PRd$ and not on its representative). 
    
    In particular, for every affine line $\ell=x_\ell+\R v_\ell$, 
    the   projection 
    $P_{[v_\ell]}$ sends the whole line  $\ell$ 
   to a single point $z_\ell$ such that  $\{z_\ell\}=\ell\cap (\R v_\ell)^\perp $. 
 The norm $\vert z_\ell\vert$  is the height of $\ell$, i.e. its distance from the origin. Thus
	an affine line $\ell = x_\ell +\R v_\ell$ can be canonically identified  with the pair $([v_\ell],z_\ell)\in\PRd\times\R^d$, and, in this coordinate,  the action   on $\AG$ of $\gamma=(B,A)\in \Aff(\R^d)$ is the following  skew product over $\mathbf{P}(\R^d)$:
	\begin{equation}\label{eq-action-Aff-z}
		\gamma\cdot\ell=\big(A\cdot[v_\ell]\ ,\ (\Pr{A\cdot[v_\ell]}A)\ z_\ell+ (\Pr{A\cdot[v_\ell]}B)\big).
	\end{equation}	
    
	We are interested in the process $\ell_n=\gamma_{n}\cdot\ell_{n-1}$ defined by the action of the random affine transformations $\gamma_n=(B_n,A_n)\in\AffRd$ on $\AG$. In the coordinates $\ell_n=( [v_{\ell_n}],z_{\ell_n})=:([\vv_n],Z_n)$, equation~(\ref{eq-action-Aff-z}) 
	shows that  $([\vv_n])_n$ is precisely the Furstenberg 
	process on $\mathbf{P}(\R^d)$ introduced in Subsection~\ref{projectiveaction}. 
    Conditioning on $([\vv_n])_n$, the transverse component $Z_n$ evolves by a stationary (Markovian but not i.i.d.) affine recursion on $\R^d$ with linear part $\Pr{ [\vv_n]} A_n$ and translation part $ \Pr{[\vv_n]}B_n$. We will show that the top Lyapunov exponent of the linear part is $\lambda_2$, which is strictly negative by assumption.  Thus $Z_n$ satisfies a stationary contracting affine recursion and a $\mu$-invariant probability measure can be constructed by adapting a classical result (see e.g.~\cite{Brandt}, on contracting stationary one-dimensional affine recursions). 
	We have the following statement.
	
    \begin{theo}\label{thm.stationary.line}
		Assume  $\bf{I}$, $\bf{P}$, ${\bf{M}}_1(\mu)$ and moreover that        $ \lambda_2<0$.
		Let $[\vv_0]=\R \vv_0$ be the random line defined in Section \ref{projectiveaction}. Then the following almost sure limit holds in $\R^d$:
		\begin{equation}\label{eq.def.Z}
			Z_0 := \lim_{n\to +\infty}  \Pr{ [\vv_0]}(\gamma_0\cdots \gamma_{-n}(0))
			= \sum_{k=0}^{+\infty}  \Pr{ [\vv_0]}A_0\cdots A_{-k+1} B_{-k}.  
		\end{equation}
		and the distribution $\eta$ of the random affine line
		\[
		\ell_0 := Z_0 + \R \mathbf v_0
		\]
		is a $\mu$-invariant probability measure
        on $\AG$, independent of $(\gamma_n)_{n\geq 1}$.
	\end{theo}
\begin{proof}
	By linearity of the projection, for every $n\in \N$,
	\begin{eqnarray}\label{eq-Z.sum}
		Z_0^n&:=&\Pr{[\vv_0]} \big( \gamma_0 \cdots \gamma_{-n}(0) \big)\nonumber\\
		&= &\Pr{[\vv_0]} \left( \sum_{k=0}^{n} A_0 \cdots A_{-k+1} B_{-k} \right) \nonumber\\
		&= &\sum_{k=0}^{n} \Pr{[\vv_0]}  A_0 \cdots A_{-k+1} B_{-k} .    
	\end{eqnarray}
	Observe that $|\Pr{[v]}x|=|\sin\angle(x,v)||x|=|x\wedge v|/|v|$, thus the general term of the sum (\ref{eq-Z.sum}) may be bounded as follows: for  $k\in\N$, 
	\begin{align}
		\left\vert \Pr{[\vv_0]}(A_0 \cdots A_{{-k+1}} B_{-k}) \right \vert
		& =  \frac{\vert A_0\cdots A_{-k+1} B_{-k} \wedge \vv_0\vert}{|\vv_0|}\nonumber\\
		& = \frac{\vert\bigwedge^2 A_0\cdots A_{-k+1} (B_{-k}\wedge {\bf v}_{-k})\vert}{\vert A_0\cdots A_{-k+1} {\bf v}_{-k}\vert}\quad\textrm{by \eqref{eq.equivarance.vk}}\nonumber\\
		& \leq  \frac{\vert\bigwedge^2 A_0\cdots A_{-k+1}\vert \vert {\bf v}_{-k}\vert}{\vert A_0\cdots A_{-k+1} {\bf v}_{-k}\vert}
		\vert B_{-k}\vert. \label{eq.estime.series}
	\end{align}
	We claim that these terms decay exponentially as $\mathrm{e}^{\lambda_2|k|}$. First, under our assumptions, as recalled in    Section \ref{sec.preliminaries}, $\P$-a.s. 
	\begin{equation}\lim_{k\to +\infty}\label{eq.eq1}\frac{1}{k}\ln \vert \bigwedge^2 A_0\cdots A_{-k+1}\vert= \lambda_1+\lambda_2 \mbox{ and }\lim_{k\to +\infty} \frac{1}{k}\ln \vert B_{-k}\vert =0 ,\end{equation}
By equivariance (\ref{eq.equivarance.vk}), we can write the logarithm of the
denominator in \eqref{eq.estime.series} as the sum $\sum_{i=0}^{k-1}F\circ\theta^{-i}$,
with $\theta$ the ergodic Bernoulli shift on $\Omega$ and
$F:=\ln|A_0\,\vv_{-1}/|\vv_{-1}||$, which is integrable under
${\bf M}_1(\overline\mu)$. Thus, by Birkhoff's ergodic theorem and
Furstenberg's formula (\ref{eq-furstemberg-lambda1}),
	\begin{equation} \label{eq.eq2}
		\lim_{k \to +\infty}\frac{1}{k}
		\ln \vert A_0\cdots A_{-k+1} \frac{\vv_{-k}}{|\vv_{-k}|}\vert 
		= \lim_{k \to +\infty}\frac{1}{k} \sum_{i=0}^{k-1} \ln |A_{-i} \frac{{\vv}_{-i-1}}{|\vv_{-i-1}|}| 
		= \lambda_1.
\end{equation}
	Combining \eqref{eq.eq1} and \eqref{eq.eq2}, we obtain, $\P$-a.s.
	\begin{equation}\frac{1}{k}\ln\left( \frac{\vert \bigwedge^2 A_0\cdots A_{-k+1}\vert \vert {\bf v}_{-k}\vert}{\vert A_0\cdots A_{-k+1} {\bf v}_{-k}\vert }
		\vert B_{-k}\vert \right)\to \lambda_1+\lambda_2+0-\lambda_1=\lambda_2<0.\end{equation}
	Thus by \eqref{eq.estime.series}, the general term in \eqref{eq-Z.sum} decays exponentially fast almost surely, hence the series converges absolutely to some $Z_0\in (\R {\bf v}_0)^{\perp}$.

Let us show that the distribution $\eta$ of the affine line
$\ell_0 := Z_0 + \R\vv_0 \in \AG$ is $\mu$-invariant and, more precisely,
that we have the equivariance property: for $\P$-almost all
$\omega = (\gamma_n)_{n\in\Z}$,  
\begin{equation}\label{eq.equivariance.ell}
	\gamma_1\cdot\ell_0(\omega) = \ell_0\circ\theta(\omega).
\end{equation}
By~(\ref{eq-action-Aff-z}),
$$
\gamma_1\cdot\ell_0 = \bigl(A_1\cdot[\vv_0],\; \Pr{A_1\cdot[\vv_0]}(A_1 Z_0 + B_1)\bigr).
$$
The equivariance of the first component follows from the equivariance of
the projective action~\eqref{eq.equivarance.vk}, which ensures that
$A_1\cdot[\vv_0] = [\vv_1] = [\vv_0]\circ\theta$; this coincides with
the direction of $\ell_0\circ\theta$.

To show \eqref{eq.equivariance.ell}, it remains to prove that  
	$
	\Pr{[\vv_1]}(\gamma_1 Z_0)=Z_0\circ \theta.  
	$
	Observe that
	\begin{align*}
		\gamma_1 Z_0=B_1+A_1 Z_0= B_1+ \sum_{k=0}^{+\infty}A_1  \Pr{[\vv_0]}(A_0\cdots A_{-k+1} B_{-k})
	\end{align*}
	Therefore
	$$ \Pr{[\vv_1]}({\gamma_1 Z_0})=
	\Pr{[\vv_1]}(B_1)+\sum_{k=0}^{+\infty}( \Pr{[\vv_1]} A_1   \Pr{[\vv_0]})(A_0\cdots A_{-k+1} B_{-k}).$$
	Now, for all $[v]\in\PRd$ and  $A\in \GL_d(\R)$, we have
	$$ \Pr{A\cdot[v]} A= \Pr{A\cdot[v]}A\Pr{ [v]}.$$
	Consequently,  since $A_1\cdot [\vv_0]=[\vv_1]$, we obtain
	$$ \Pr{[\vv_1]}(\gamma_1 \,Z_0)= \Pr{[\vv_1]}(B_1)+\sum_{k=0}^{+\infty}{ \Pr{[\vv_1]} (A_1A_0\cdots  A_{-(k-1)} B_{-k})}=Z_0\circ \theta.$$
	This proves that $\gamma_1\cdot \ell_0=\ell_0 \circ\theta\ \P$-a.s., hence $\eta$ is a $\mu$-invariant probability measure on $\AG$.
\end{proof}

\subsection{Construction of the invariant measure $\sigma$ on $\cZ_d$ }\label{sec.invariant.fibers}
We work in point-vector coordinates $(x,v)$ on $\cZ_d\simeq \R^d\times (\R^d\setminus \{0\})$ introduced at the beginning of Section \ref{sec.diag.action}, in which the action is given by \eqref{eq.action.Ab}.
As observed, in dimension $d=1$, the product space $\cZ_1$ can be identified with the affine group $\Aff(\R)$, and thus the left Haar measure on the group gives a     Radon measure $m$  on $\cZ_1\cong \R\times (\R\setminus \{0\})$
\begin{equation}\label{eq.def.m}
	m(f):=\int_{\cZ_1} f(b,a) \frac{db\,da}{|a|^2}
\end{equation}
that is $\Aff(\R)$-invariant, i.e.\ for every $\tau=(\beta,\alpha)\in \Aff(\R)$,
\begin{equation*}
	\tau_*m(f)= \int_{\cZ_1} f(\alpha b+\beta, \alpha a)\, dm(b,a)=m(f).
\end{equation*}
\subsubsection{Measures on fibers $\cZ_1(\ell)$}

Let $\ell=x_\ell+\R v_\ell\in \AG$. For a given choice of coordinates, consider the identification of $\cZ_1$ with $\cZ_1(\ell)\subset \cZ_d$ via $\varphi_\ell$ 
(see  \eqref{phiell}). We define $\lambda^\ell:=(\varphi_\ell)_*m$ as the pushforward measure of $m$ to $\cZ_d$ by $\varphi_\ell$, more precisely, for any measurable function $f$ on $\cZ_d$,
\[
\lambda^\ell(f):=(\varphi_\ell)_*m(f)= \int_{\cZ_1} f(x_\ell+b v_\ell, a v_\ell)\, dm(b,a)
\]

We first establish key properties of the family $(\lambda^\ell)_{\ell\in\AG}$: well-definedness, equivariance under $\Aff(\R^d)$, and behaviour under linear projection.
For $\tilde v\in (\R^d)^*$ and $\ell\in \AG$ such that   $\tilde v v_\ell\neq 0$,
we still denote by $\tilde v$ the map  $\cZ_1(\ell)\to  \cZ_1, (x,v)\mapsto (\tilde v x, \tilde v v)$. 

\begin{prop}\label{prop.def.lambda}
The definition of $\lambda^\ell$ does not depend on the choice of the coordinates $x_\ell$ and $v_\ell$ on the line $\ell$.
Furthermore, $
\gamma_*\lambda^\ell=\lambda^{\gamma\cdot \ell}
$ for every $\gamma\in \Aff(\R^d)$
and $ \tilde v_*\lambda^\ell=m$ for any linear form $\tilde v\in (\R^d)^*$ such that $\tilde v v_\ell\neq 0$.
\end{prop}
\begin{proof}
Let $(x_\ell',v_\ell')$ be another choice of coordinates for the same line $\ell=x_\ell+\R v_\ell$ and let $\varphi_\ell'$ be the associated map.
Then, there exists $\tau=(\beta,\alpha)\in \Aff(\R)$ such that
\[
x_\ell'=x_\ell+\beta v_\ell,\qquad v_\ell'=\alpha v_\ell.
\]
Since $m$ is invariant under the action of $(\beta,\alpha)\in \Aff(\R)$:
\begin{align*}
    (\varphi_\ell')_*m(f)&=\int_{\cZ_1} f(x_\ell'+b v_\ell', a v_\ell')\, dm(b,a)\\
    &=\int_{\cZ_1} f(x_\ell+(\beta +\alpha b ) v_\ell, \alpha a v_\ell)\, dm(b,a)\\
    &=(\varphi_\ell)_*(\tau_*m)(f)=(\varphi_\ell)_*m(f)=\lambda^\ell(f).
\end{align*}
For $\gamma=(B,A)\in \Aff(\R^d)$, the image of $\ell$ is $\gamma\cdot \ell=\gamma x_\ell+\R (A v_\ell)$. By independence of the choice of coordinates, we may use $(\gamma x_\ell, A v_\ell)$ as point-vector coordinates on the line $\gamma\cdot \ell$.  Then, for $\varphi_\ell(b,a)\in \cZ_1(\ell)$ with affine coordinates $(b,a)$, we have
$$\gamma \varphi_\ell(b,a)=\gamma(x_\ell+b v_\ell,a v_\ell)=(\gamma x_\ell+ b A v_\ell,a Av_\ell)=\varphi_{\gamma \cdot \ell}(b,a).$$
Hence $\gamma \varphi_\ell=\varphi_{\gamma\cdot \ell}$ and
$\gamma_*\lambda^\ell=(\gamma\varphi_\ell)_*m  =(\varphi_{\gamma\cdot\ell})_*m=\lambda^{\gamma \cdot \ell}.$

Finally  let $\tilde v\in (\R^d)^*$ such that $\tilde v v_\ell\neq 0$. A direct computation shows that $$\tau:=\tilde v\circ \varphi_\ell: \cZ_1\to \cZ_1,\  (b,a)\mapsto((\tilde v v_\ell)b+(\tilde v x_\ell),( \tilde v v_\ell)a)
	.$$
is a one-dimensional affine map.	Since $m$ is $\Aff(\R)$-invariant, it follows that $\tilde v_\ast \lambda^\ell=(\tilde v \varphi_\ell)_\ast m=m$.
\end{proof}

\subsubsection{The $\mu$-invariant Radon measure $\sigma$ on  $\cZ_d$}
We define the measure $\sigma$ on $\cZ_d$ by
	\begin{equation}\label{eq.def.sigma}
	    \sigma(f) := \int_{\AG} \lambda^{\ell_0}(f)\, d\eta(\ell_0)=\int_{\cZ_1} \E(f(Z_0+b\vv_0,a\vv_0))\, dm(b,a),
	\end{equation}
	where $\ell_0=Z_0+\R\vv_0\sim\eta$ is the $\mu$-invariant probability measure on $\AG$ constructed in 
	Theorem~\ref{thm.stationary.line}.
	
	\begin{prop}\label{prop.invariant.sigma}
		The measure $\sigma$ is a positive and non-zero $\mu$-invariant Radon measure on $\cZ_d$.
	\end{prop}

\begin{proof}{}
For every $\gamma\in \Aff(\R^d)$,
by the equivariance property given in Proposition \ref{prop.def.lambda},
$$
\gamma_\ast \sigma=\int{\gamma_\ast  \lambda^{\ell_0}} d\eta(\ell_0)=\int{\lambda^{\gamma\cdot\ell_0} d\eta(\ell_0)} 
= \int{ \lambda^{\ell_0}} d(\gamma_\ast  \eta)(\ell_0).$$
Since $\eta$ is $\mu$-invariant (Theorem \ref{thm.stationary.line}), we deduce that
$$\int{\gamma_\ast \sigma \,d\mu(\gamma)}=\iint{\lambda^{\ell_0} d(\gamma_\ast \eta)(\ell_0) d\mu(\gamma)} =\int{\lambda^{\ell_0} d\eta(\ell_0)}=\sigma,$$ i.e. $\sigma$ is $\mu$-invariant on $\cZ_d$. The measure $\sigma$ is positive and non-zero.   
Furthermore, it is a Radon measure. Indeed, let $K\subset \cZ\simeq (\R^d)\times (\R^d\setminus\{0\})$ be a compact subset.
 There exists 
 $M>1$ such  that $K\subset K_M:=\{(x,v) \mid |x|\leq M, 1/M\leq |v|\leq M\}$. Fix 
 $\ell_0=Z_0+\R {\bf v}_0$  with 
 $Z_0\perp {\bf v}_0$ and $|{\bf v}_0|=1$. 
 Then 
 \begin{align*}
\lambda^{\ell_0}(K)\leq \lambda^{\ell_0}(K_M)&=m(\{(b,a) \mid |Z_0+ b {\bf v}_0|\leq M, 1/M\leq |a {\bf v}_0|\leq M\})\\
 &\leq m(\{(b,a) \mid |b|\leq M, 1/M\leq |a|\leq M\})\quad \textrm{ since $Z_0\perp {\bf v}_0$}\\
&=\left(\int_{-M}^M{db}\right)\left(
2 \int_{1/M}^M{\frac{da}{a^2}}\right)\leq 4M^2.
 \end{align*}
 Integrating with respect to the probability measure $\eta(d\ell_0)$, we get that $\sigma(K)<4M^2<+\infty$. This proves that $\sigma$ is a Radon measure on $\cZ_d$.  
 
 \end{proof}

 \subsection{Log-integrability of the height of the invariant line}\label{sec.integrable.line}
In this section,  we show that the distance of the origin to the invariant line $\ell_0=Z_0+\R {\bf v}_0$ constructed in Section \ref{sec.line} satisfies a logarithmic first moment condition, a crucial estimate needed in the next section. 
\begin{prop}\label{prop.line.integrable} 
 Under hypotheses ${\bf  M}_2(\mu)$, {\bf I} and {\bf P}  with  $\lambda_2<0$,  the random variable $Z_0$ defined by \eqref{eq.def.Z}  satisfies 
$
\E\big(\ln^+ \vert Z_0\vert\big) < +\infty$. 
\end{prop}

\begin{proof}
The proof of Theorem \ref{thm.stationary.line} relies on the fact that the
general term of the series defining $Z_0$ in \eqref{eq.def.Z} decays
exponentially, with rate $\lambda_2<0$. To prove that $\ln^+ |Z_0|$ is
integrable, we need the second moment assumption ${\bf M}_2(\mu)$, which
allows us to obtain an integrable bound, uniform in $n$.

First, we fix $0<\epsilon< -\lambda_2/3$ and set $K:=e^{-\lambda_2-3\epsilon}>1$. Observe that    
\begin{align*}
\vert Z_0\vert&\leq \sum_{n=0}^{+\infty}\left \vert P_{[\vv_0]}\left( A_0 \cdots A_{-n+1} B_{-n} \right)\right\vert\\
&= \sum_{n=0}^{+\infty} K^{-n} K^n\left\vert  P_{[\vv_0]}\left( A_0 \cdots A_{-n+1} B_{-n} \right)\right\vert\\
&\leq \frac{1}{1-K^{-1}}\sup_{n\geq 0} K^n\left\vert  P_{[\vv_0]}\left( A_0 \cdots A_{-n+1} B_{-n} \right)\right\vert.
\end{align*}
This yields 
\begin{equation}\label{eq.Z.sup}
  \ln\vert Z_0\vert\leq \ln \frac{1}{1-K^{-1}}+    \ln\left( \sup_{n\geq 0} K^n\left\vert  P_{[\vv_0]}\left( A_0 \cdots A_{-n+1} B_{-n} \right)\right\vert\right).  
\end{equation}
Fix an orthonormal basis $\{e_1,\ldots, e_d\}$ of $\R^d$ and write $B_{-n}=\sum_{i=1}^d{B_{-n,i} e_i}$ in this basis. 
Observe that $\vert P_{[\vv_0]}\,x\vert =|\sin\angle(x,\vv_0)||x| = \delta([x],[\vv_0])\,\vert x\vert$.
We decompose the term $\vert  P_{[\vv_0]}(A_0 \cdots A_{-(n-1)} B_{-n})|$  which appears in the right hand side   of \eqref{eq.Z.sup} as follows:  
\begin{align}\label{3factors}
 &\left\vert  P_{[\vv_0]}(A_0 \cdots A_{-n+1} B_{-n} ) \right\vert
 =  \left\vert  \sum_{i=1}^d B_{-n,i} \  P_{[\vv_0]}(A_0 \cdots A_{-n+1} e_i ) \right\vert  \nonumber \\
 &\leq \sum_{i=1}^d\vert B_{-n,i}\vert \ \delta(A_0\cdots A_{-n+1} \cdot[e_i],[\vv_0])\,\vert A_0\cdots A_{-n+1}  e_i\vert\
 \nonumber\\
  &\leq d  \left(\max_{1\leq i \leq d}\delta(A_0\cdots A_{-n+1}\cdot[e_i],[\vv_0])\right)\,\vert A_0\cdots A_{-n+1} \vert \,\vert B_{-n}\vert.
\end{align}
The first factor in \eqref{3factors} deals with the contraction properties of the products $A_0\cdots A_{-n+1}$   and decays exponentially as $\mathrm{e}^{-(\lambda_1-\lambda_2)n}$ (recall that $\lambda_1-\lambda_2>0$ under $\bf I$ and $\bf P$);  the  second one  behaves as $e^{\lambda_1n}$  while the last one is  a stationary factor related to the translation part of $\gamma_{-n}$. Let us  formalize these observations. We set 
 $ M_\delta =e^{\lambda_1-\lambda_2-\epsilon},  M_A=e^{-\lambda_1-\epsilon}$ and $ M_B=e^{-\epsilon}$  (observe that $M_\delta M_A M_B= K$) and 
 write
\begin{align} 
& K^n  P_{[\vv_0]}\left( A_0 \cdots A_{-(n-1)} B_{-n} \right)\vert \nonumber\\
&\leq d
\Big({M_\delta^n \max_{1\leq i \leq d}\delta(A_0\cdots A_{-n+1}\cdot[e_i],[\vv_0])}\Big) \Big({M_A^{n}\vert A_0\cdots A_{-n+1}\vert}\Big)\Big({M_B^{n} \vert B_{-n}\vert}\Big)
\label{kn.pro}\end{align} 
We deal with each factor separately.

\emph{Step 1: $\delta$-factor.} First, we show that
\begin{equation}\label{eq.sup.delta}
\E\left(\sup_{n\geq 0}{\ln^+\left(M_\delta^{n}\left(\max_{1\leq i \leq d}\delta(A_0\cdots A_{-n+1}\cdot[e_i],[\vv_0])\right)\right)}\right)<+\infty.    
\end{equation}

Set $\delta_n:=\left(\max_{1\leq i \leq d}\delta(A_0\cdots A_{-n+1}\cdot[e_i],[\vv_0])\right)$. Using the large deviation estimate \cite[Theorem B]{Peneau}, since $M_\delta < e^{\lambda_1-\lambda_2}$, there exists $C, \beta>0$ (depending on $M_\delta$) such that 
$\P(\delta_n>M_\delta^{-n})\leq C \mathrm{e}^{-\beta n}.$ 
 Hence 
\begin{align*}\E(\sup_{n\geq 0}{\ln^+(M_\delta^{n} \delta_n)})&\leq 
\sum_{n=0}^{+\infty}{\E(\ln^+(M_\delta^{n} \delta_n))}\nonumber\\
& =   \sum_{n=0}^{+\infty}\E(\ln^+(M_\delta^{n} \delta_n) \mathds{1}_{\delta_n\geq M_\delta^{-n}})\nonumber\\
& \leq   (\ln^+ M_\delta) \sum_{n=0}^{+\infty}{n\ \P(\delta_n>M_\delta^{-n})} \quad \quad \,\,\,\textrm{(since $\delta_n\leq 1$)}\nonumber\\
&\leq  (\ln^+{M_\delta})\sum_{n=0}^{+\infty}{n\ e^{-\beta n}}<+\infty. 
\end{align*}

\emph{Step 2 : $A$-factor.} The second factor $M_A^n \vert A_0\cdots A_{-n+1}\vert$ is split into   a main contribution  and a remainder we control jointly with the third factor. By definition of the Lyapunov exponent  and since $M_A< e^{-\lambda_1}$, there exists $n_0\in \N$ such that 
$$\frac{1}{n_0}\E(\ln \vert A_0\cdots A_{-n_0+1}\vert )< -\ln M_A.$$
We decompose the product $A_0\cdots A_{-n+1}$ into blocks of length $n_0$. More precisely, let $k(n)=\lfloor n/n_0\rfloor$ and $r(n)=n- k(n)n_0$.
For every $i\geq 1$, let 
$$Y_i:=M_A^{n_0}\vert A_{-(i-1)n_0}\cdots A_{-i n_0+1}\vert \mbox{ and }W_n:=M_A^{r(n)} \vert A_{-k(n) n_0}\cdots A_{-n+1}\vert.$$
By submultiplicativity of the matrix norm,
\begin{equation}\label{eq.decom}
M_A^n 
\vert A_0\cdots A_{-n+1}\vert\
 \leq \left(M_A^{n_0 k(n)} \vert A_0\cdots A_{-k(n) n_0+1}\vert\right) W_n
\leq\Big(Y_1\cdots Y_{k(n)} \Big) W_n.
\end{equation}
The sequence $(\ln Y_i)_{i\geq 1}$ is i.i.d. with moment of order 2; furthermore, by the choice of $n_0$ and $M_A$, the $\ln Y_i$ have a negative expectation.
 Using the main result of \cite{KieferWolfowitz}, we deduce that the positive part of its supremum is integrable, in other words that
\begin{equation}\label{eq.sup.norm}
\E\left(\left(\sup_{n\geq 0}{\ln(M_A^{n_0 k(n)}\vert A_0\cdots A_{-k(n)n_0+1}\vert)}\right)^+\right)\leq \E\left(\left(\sup_{k\geq 1}{\sum_{i=1}^k{\ln Y_i}}\right)^+\right)<\infty.    
\end{equation}
The remaining term $W_n$ is controlled hereafter.

\emph{Step 3: $B$-factor and reminder term.} Finally, we show that
\begin{equation}\label{eq.sup.Wb}
    \E(\sup_{n\geq 0}{\ln^+(M_B^{n} W_n \vert B_{-n}\vert )})<+\infty.
\end{equation}

Distinguishing according to the remainder $r(n)$ modulo $n_0$, we write  
\begin{align*}
\E(\sup_{n\geq 0}{\ln^+(M_B^n W_n \vert B_{-n}\vert)})&\leq \sum_{r=0}^{n_0-1}\sum_{k=0}^{+\infty}{\E(\ln^+(M_B^{n_0 k +r} X_k^{(r)}))},  \quad\end{align*}
with $X_k^{(r)}:=W_{n_0 k+r} \vert B_{-n_0 k -r}\vert$. For every $0\leq r \leq n_0-1$,  the variables $W_{n_0 k+r}$, $k\geq 0$,  are identically distributed. Furthermore,   the variables $B_{-n_0k -r}$ and $W_{n_0k+r}$   are independent for every $k\geq 0$. Therefore the random variables $X_k^{(r)}$, $k\geq 0$,  have the same distribution. Moreover  $\E((\ln^+ X_1^{(r)})^2)<+\infty$ thanks to assumption ${\bf M}_2(\mu)$ and the submultiplicativity of the matrix norm.
Consequently,  since $M_B<1$, 
\begin{align*}
 \E(\sup_{n\geq 0}{\ln^+(M_B^n W_n \vert B_{-n}\vert)})&\leq \sum_{r=0}^{n_0-1}\sum_{k=0}^{+\infty}{\E(\ln^+(M_B^{n_0 k +r} X_1^{(r)}))}\\
 &= \sum_{r=0}^{n_0-1}\sum_{k=0}^{+\infty} \E( \ln^+ (M_B^{n_0 k+r} X_1^{(r)})\mathds{1}_{M_B^{n_0 k+r} X_1^{(r)} > 1}))\\
& \leq  \sum_{r=0}^{n_0-1}\sum_{k=0}^{+\infty} \E\left(\ln^+(X_1^{(r)}) \mathds{1}_{
\ln^+(X_1^{(r)}) > (n_0 k +r) \ln (1/M_B)}\right),
\end{align*}
which is finite since  $\mathbb E((\ln^+(X_1^{(r)}))^2)<+\infty$ for every $r$.  
Combining  \eqref{kn.pro}, \eqref{eq.sup.delta}, 
\eqref{eq.decom}, \eqref{eq.sup.norm} and \eqref{eq.sup.Wb}, we obtain the log-integrability of the random variable 
$
\sup_{n\geq 0} K^n\left\vert P_{[\vv_0]}\left( A_0 \cdots A_{-n+1} B_{-n} \right)\right\vert $. Thus   $\E(\ln^+\vert Z_0 \vert)<+\infty$, by \eqref{eq.Z.sup}. 
\end{proof}

\section{Proof of the Main Theorems}\label{sec.end.proof}

In this section we prove Theorem~\ref{theorem.main.crossing}, which
asserts that the joint process $(X_n^x,X_n^y)$ eventually stops crossing
any fixed interval in any direction. As a straightforward
corollary, we deduce the local contraction result of
Theorem~\ref{theorem.main.contraction}.
Throughout this section, we assume the hypotheses of
Theorem~\ref{theorem.main.contraction} (in particular
$\lambda_2<\lambda_1=0$) and we use pair coordinates $(x, y)$ on $\cZ_d$. We denote by $\sigma$ the $\mu$-invariant Radon measure on $\cZ_d$ defined in \eqref{eq.def.sigma} and by $P_\mu$ the Markov operator
of the joint process $(X_n^x,X_n^y)$ on $\cZ_d$, that is, for any
Borel measurable function $f:\cZ_d\to\R_+$ and $\fz=(x,y)\in\cZ_d$,
\begin{equation}\label{def.Pmu}
P_\mu f(\fz)=\int_{\AffRd} f(\gamma \fz)\,d\mu(\gamma).
\end{equation}

The proof relies on the total dissipativity of $P_\mu$ on
$L^1(\cZ_d,\sigma)$, stated in Theorem~\ref{theorem.affine.dissipative};
its proof is postponed to Section~\ref{section.dissipative}.

\begin{figure}[h]
\centering
\figu{0.4}{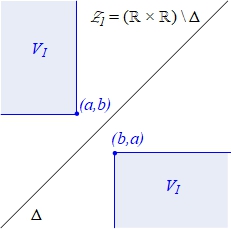}  
\caption{One-dimensional crossing set $V_I\subset \R\times \R$ for $I=[a,b]$.}
\label{fig.VI}
\end{figure}

Let $I = [a,b] \subset \R$ be a non-degenerated compact interval, then we define the 
\emph{one-dimensional crossing set} (see Figure~\ref{fig.VI}) as
\[
V_I := \bigl\{ (s,t) \in \cZ_1 \mid [s,t]  \supset I\bigr\}\subset \cZ_1 = (\R \times \R) \setminus \Delta,
\] 
that is, the set of pairs $(s,t)$ such that 
$\min\{s,t\} \leq a$ and $\max\{s,t\} \geq b$.

For $\tilde v\in(\R^d)^*\setminus\{0\}$, we still denote by $\tilde v$
the diagonal map $\cZ_d\to\R\times\R$, $(x,y)\mapsto(\tilde vx,\tilde vy)$.
The crossing set $V_{\tilde v,I}$ introduced in \eqref{eq.def.crossing.set}
is then  equal to $\tilde v^{-1}(V_I)$.

\subsection{Integrability of the one-step crossing function}\label{sec.cross.integrable}

\begin{prop}\label{prop.crossing.integrable} Assume {\bf I}, {\bf P}, ${\bf M}_2(\mu)$ and $\lambda_2<0$.
Let $\tilde v\in(\R^d)^*\setminus \{0\}$ and $I\subset\R$ be a compact non-degenerate
interval. Let $\gamma_1\sim\mu$ and define
$\phi_{\tilde v,I}:\cZ_d\to\R_+$ by
\[
\phi_{\tilde v,I}(x,y)
:=
\mathds{1}_{V_{\tilde v,I}^c}(x,y)\,
\P\bigl(\gamma_1 (x,y) \in V_{\tilde v,I}\bigr).
\]
Then $\phi_{\tilde v,I}\in L^1(\cZ_d,\sigma)$
\end{prop}
Notice that  $\phi_{\tilde v,I}(x,y)$ is the probability that the
joint process starting at $(x,y)$ enters $V_{\tilde v,I}$ at time $1$.
We start with a lemma in dimension $d=1$. The left Haar measure $m$ on $\Aff(\R)$ transported to $\cZ_1$
can be written, in pair coordinates, as
$dm(s,t)=\frac{ds\,dt}{(s-t)^2}$.
For $A\subset \R$,
denote $\mathrm{diam}(A):=\sup\{|s-t|; s,t\in A\}.$

\begin{lemma}\label{lemma.haar.crossing}
Let $I$ and $J$ be two compact non-degenerate intervals of $\R$. Then
\begin{equation}\label{eq.haar.intervals}
m(V^c_I\cap V_J)= 2\ln\frac{\mathrm{diam}(I\cup J)}{\mathrm{diam}(J)}.
\end{equation}
Furthermore, for $\tau=(\beta_\tau, \alpha_\tau)\in \Aff(\R)$,  we have
\begin{align}
m(\{(s,t)\in\cZ_1\mid [s,t]\not\supset I, \tau[s,t]\supset I\})
&=2\ln \frac{\mathrm{diam}(\tau I\cup I)}{\mathrm{diam}(I)}\label{eq.haar.equal} \\
& \leq C_I( 1+\ln^+|\alpha_\tau|+\ln^+|\beta_\tau|),
\label{eq.haar.affine}
\end{align}
where $C_I>0$ depends only on the interval $I$.
\end{lemma}

\begin{proof}{}
Set  $E=\{(s,t)\in \cZ_1 \mid  s<t, [s,t]\not\supset I, [s,t]\supset J\}$.  By symmetry, $m(  V^c_I\cap V_J)=2m(E)$.
Set $I=[a,b]$ and $J=[c,d]$, with $a<b$ and $c<d$. To prove \eqref{eq.haar.intervals}, we distinguish three cases according to the relative position of $I$ and $J$ (that is, of $V_I$ and $V_J$; see Figure \ref{fig:Lem-m-Cross}).
\begin{figure}
    \centering
    \figu{0.8}{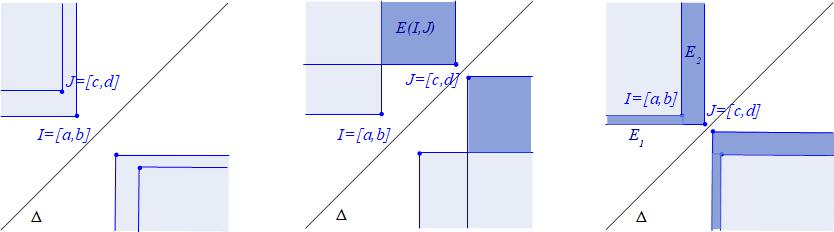}  
    \caption{The sets $V_I$, $V_J$ and $V^c_I\cap V_J$ in Cases 1, 2.ii and 3.}
    \label{fig:Lem-m-Cross}
\end{figure}

\textbf{Case 1. $I\subset J$, that is $(c, d)\in V_I$.}
In this case, $V_J\subset V_I$ and both side of \eqref{eq.haar.intervals} are equal to $0$. 

\textbf{Case 2. $I\not\subset J$ and $J\not\subset I$, that is $(c,d)\not\in V_I$ and $(a,b)\not\in V_J$.} We distinguish two symmetric situations:
\begin{enumerate}
    \item[2.i.]  Suppose $c<a$ and $d<b$. Hence $E=\{(s,t)\mid s\leq  c,\ d\leq  t<b\}$ so that 
$$m(E)=\int_{-\infty}^c\left(\int_d^b\frac{1}{(s-t)^2}dt\right)\,ds=\int_{-\infty}^c{\left(\frac{1}{d-s}-\frac{1}{b-s}\right) ds}=
\ln \frac{b-c}{d-c}.$$
Since $\mathrm{diam}(I\cup J)=b-c$, 
\eqref{eq.haar.intervals} holds. 
\item[2.ii.] Suppose $c>a$ and $d>b$. Hence
   $E=\{(s,t)\mid a<s\leq  c,\ t\geq  d\}$ so that 
$$m(E)=\int_a^c \left(\int_d^{+\infty}\frac{1}{(s-t)^2}dt\right) ds=\int_a^c{\frac{1}{d-s} ds}=\ln \frac{d-a}{d-c}.$$
Since $\mathrm{diam}(I\cup J)= d-a$, \eqref{eq.haar.intervals} holds. 
\end{enumerate}

\textbf{Case 3. $J\subset I$ that is $(a,b)\in V_J$.} In this case, $a\leq c<d\leq b$.
   We write $E=E_1\sqcup E_2$, where
   $E_1=\{(s,t)\mid s\leq  a,\ d\leq  t<b\}$ and
   $E_2=\{(s,t)\mid a<s\leq  c,\ t\geq  d\}$. 
A similar computation to that of Case 2 shows that 
$$m(E)=m(E_1)+m(E_2)
=\ln\frac{b-a}{d-a}+\ln\frac{d-a}{d-c}=\ln\frac{b-a}{d-c}.
$$
Since $I\cup J=I$, \eqref{eq.haar.intervals} holds. 
 Identity \eqref{eq.haar.intervals} is now proved.

To prove \eqref{eq.haar.equal} observe  that 
 $$A:=\{(s,t)\mid [s,t]\not\supset I, \tau[s,t]\supset I\}=  V^c_I \cap V_{\tau^{-1}I}=  V^c_I \cap \tau^{-1} V_{I}.$$ 
Hence,     since $m$ is invariant under the diagonal action of
   $\Aff(\R)$, applying \eqref{eq.haar.intervals} with $\tau I$ for $I$ and $I$ for $J$,  
 $$m(A)=m(\tau A)
 =m(  V^c_{\tau I} \cap V_{I})=2\ln \frac{\mathrm{diam}(\tau I\cup I)}{\mathrm{diam}(I)}.$$
 This proves \eqref{eq.haar.equal}.
Finally,
\begin{align*}
    \mathrm{diam}(I&\cup\, \tau I)=\max\{a,b,\tau a,\tau b\}-\min\{a,b,\tau a,\tau b\}\\
    &\leq 2 \max\{|a|,|b|,|\tau a|,|\tau b|\}= 2 \max\{|a|,|b|,|\alpha_\tau a+\beta_\tau|,|\alpha_\tau b+\beta_\tau|\}
\end{align*}
 Consequently
   $\mathrm{diam}(I\cup\, \tau I)\leq C'_I(1+|\alpha_\tau| + |\beta_\tau|)$
   where $C'_I:=2\max\{1,|a|,|b|\}$. 
Using the inequality $\ln^+(x+y)\leq \ln^+x+\ln^+y+\ln 2$ for $x,y>0$, 
\eqref{eq.haar.affine} follows. 
\end{proof}
  \medskip

\begin{proof}[of Proposition \ref{prop.crossing.integrable}]
By Tonelli's theorem,
\begin{align*}
&\int \phi_{\tilde v,I} \, d\sigma\\
&= \int_{\AffRd}\left(\int_{\AG} \left(\int_{\cZ_d} \mathds{1}_{ V_I^c }(\tilde v(x,y))\mathds{1}_{V_I}(\tilde v \gamma_1 (x,y)) d\lambda^{\ell_0}(x,y)\right)d\eta(\ell_0)\right)d\mu(\gamma_1)
\end{align*}
We first estimate
$\displaystyle \cI(\gamma_1,\ell_0):=\int_{\cZ_d} \mathds{1}_{ V_I^c }(\tilde v(x,y))\mathds{1}_{V_I}(\tilde v \gamma_1 (x,y)) d\lambda^{\ell_0}(x,y)$. By strong irreducibility of $\Gamma_{\overline{\mu}}$, the projection $\nu$ of $\eta$ on $\PRd$ is non-degenerate (see \cite[Chapter III, Proposition 2.3]{Bougerol-Lacroix}). Hence, by $\mu$-stationarity of $\nu$, for $\mu$-a.e.~$\gamma_1$ and $\eta$-a.e.~$\ell_0=Z_0+\R \vv_0$ one has $\tilde v \vv_0\neq 0$ and $\tilde{v} A_1 \vv_0\neq 0$. Fix such a line $\ell_0$, with $|\vv_0|=1$, and such an element $\gamma_1$.
Then the restrictions $\tilde v|_{\ell_0}$ and $\tilde v|_{\gamma_1\cdot \ell_0}$ are affine bijections onto $\R$. Let $\tau_{\ell_0,\gamma_1}\in\Aff(\R)$ be the unique affine map making the following diagram commute:
\[
\begin{array}{ccc}
\ell_0 & \xrightarrow{\ \gamma_1\ } & \gamma_1\cdot \ell_0 \\
{\scriptstyle \tilde v}\big\downarrow & & \big\downarrow{\scriptstyle \tilde v} \\
\R & \xrightarrow{\ \tau_{\ell_0,\gamma_1}\ } & \R
\end{array}
\]
that is, $\tau_{\ell_0,\gamma_1}$ is characterized by
$
\tilde v (\gamma_1 x)=\tau_{\ell_0,\gamma_1}(\tilde v x) \mbox{ for every }x\in \ell_0 .
$
A direct computation shows that $\tau_{\ell_0,\gamma_1}=(\beta_{\ell_0,\gamma_1},\alpha_{\ell_0,\gamma_1})\in \AffR$ with
$$\beta_{\ell_0,\gamma_1}:=\tilde v(A_1 Z_0+B_1)-\alpha_{\ell_0,\gamma_1}\, \tilde v Z_0, \quad \alpha_{\ell_0,\gamma_1}:=\frac{\tilde v A_1 \vv_0}{\tilde v
\vv_0}.$$
Therefore, since the support of $\lambda^{\ell_0}$ is $\cZ_1(\ell_0)=(\ell_0\times\ell_0)\setminus\Delta$ and since $\tilde v_\ast \lambda^{\ell_0}=m$ on $\cZ_1$ by Proposition \ref{prop.def.lambda}, it holds that
\begin{align*}
\cI(\gamma_1,\ell_0)&=\int_{\cZ_1(\ell_0)} \mathds{1}_{ V_I^c }(\tilde v(x,y))\mathds{1}_{V_I}( \tau_{\ell_0,\gamma_1}\tilde v (x,y)) d\lambda^{\ell_0}(x,y) \\
&=\int_{\cZ_1} \mathds{1}_{ V_I^c }((s,t))\mathds{1}_{V_I}( \tau_{\ell_0,\gamma_1}(s,t)) dm(s,t) \\
&= m\bigl(\left\{(s,t)\in \cZ_1;\ [s,t]\not\supset I,\ \tau_{\ell_0,\gamma_1}[s,t]\supset I \right\}\bigr)
\\
&\leq C_I (1+\ln^+|\alpha_{\ell_0,\gamma_1}| + \ln^+|\beta_{\ell_0,\gamma_1}|), \,  \textrm{by Lemma \ref{lemma.haar.crossing}.}
\end{align*}
Since
$\displaystyle
\vert \alpha_{\ell_0,\gamma_1}\vert\leq \frac{\vert \tilde v \vert \vert A_1\vert }{|\tilde v \vv_0|}$ and $\displaystyle
|\beta_{\ell_0,\gamma_1}|\leq \vert \tilde v \vert\left( \vert B_1\vert +\vert A_1\vert \vert Z_0\vert\Big(1+\frac{|\tilde v|}{|\tilde v \vv_0|}\Big)\right)$, we obtain
\begin{align*}
  &  \int{\phi_{\tilde v, I} d\sigma}
    = \int_{\AffRd}\int_{\AG} \cI(\gamma_1,\ell_0) d\eta(\ell_0) d\mu(\gamma_1)\\
&\leq C_{\tilde v,I}\Big(1+\E(\ln^+\vert A_1\vert)+ \E(\ln^+\vert B_1\vert)+\E(\ln^+\vert Z_0\vert) +  \int{\ln \frac{1}{|\tilde v \vv_0| } d\nu([{\vv_0}])}\Big).
\end{align*}
with $\E(\ln^+\vert Z_0\vert)<+\infty$ by Proposition \ref{prop.line.integrable}, and $\int \ln\frac{1}{|\tilde v \vv_0|}d\nu([\vv_0])<+\infty$ by \cite[Proposition 4.5]{BQCLT}, under ${\bf M}_2(\mu)$.
\end{proof}
\medskip 
\medskip

\subsection{Proof of Theorems~\ref{theorem.main.crossing} and~\ref{theorem.main.contraction}}\label{sec.proof.main}
In this subsection, we prove Theorems
\ref{theorem.main.crossing} and \ref{theorem.main.contraction}.
The fact that $\sigma$ is $\mu$-invariant is equivalent to
$\sigma P_\mu=\sigma$, that is
$\int{(P_\mu f)(\fz)\, d\sigma(\fz)}=\int{f(\fz)\, d\sigma(\fz)}$, which
implies that $P_\mu$ acts as a positive contraction on
$L^1(\cZ_d, \sigma)$.We assume Theorem~\ref{theorem.affine.dissipative}, proved in
Section~\ref{section.dissipative}, which states that $P_\mu$ is totally
dissipative on $L^1(\cZ_d,\sigma)$, i.e.
\begin{equation}\label{eq.total.before}
\forall f\in L^1(\cZ_d,\sigma), f\geq0, \quad
\sum_{n=0}^{+\infty}{P_\mu^n f}<+\infty \quad \textrm{$\sigma$-a.e.}
\end{equation}

\begin{proof}[of Theorem \ref{theorem.main.crossing}]
For $\tilde v\in (\R^d)^*\setminus\{0\}$, $(x,y)\in\R^d\times\R^d$
and $I\subset \R$ a compact non-degenerate
interval, set
$$\cC_{(x,y)}^I:=\{\omega\in \Omega \mid
[\tilde v X_n^x(\omega), \tilde vX_n^y(\omega)]\not\supset I\ \text{for all large } n \},$$
where $X_n^x(\omega)=\gamma_n\cdots \gamma_1(x)$. Since a single point
contains no non-degenerate interval, $\cC_{(x,x)}^I=\Omega$. We want to
show that
$\displaystyle \bigcap_{(x,y)\in\R^d\times\R^d}\bigcap_{I\subset\R}
\cC_{(x,y)}^I$ contains an event of probability one.

\noindent\emph{Step 1: for every interval $I\subset\R$, $\P(\cC_{(x,y)}^I)=1$ for $\sigma$-almost all $(x,y)\in\cZ_d$.}\\
By \eqref{eq.total.before} and the integrability of $\phi_{\tilde v,I}$ with respect to $\sigma$ (Proposition \ref{prop.crossing.integrable}),
the expected number of   entrances of the process $\tilde v (X_n^x, X_n^y)$   to $V_{I}$ is finite for $\sigma$-a.e.\ $(x,y)$, i.e.
\begin{align*}
    \sum_{n=0}^{+\infty}&\P([\tilde v X_n^x, \tilde v X_n^y]\not\supset I,   [\tilde v X_{n+1}^x, \tilde v X_{n+1}^y]\supset  I)\\
    &=
    \sum_{n=0}^{+\infty}\P(\tilde v[ X_n^x,   X_n^y]\not\supset I,  \tilde v \gamma_{n+1}[ X_{n}^x, X_{n}^y]\supset I)
    \\&=\sum_{n=0}^{+\infty} \E \left(\phi_{\tilde v,I}(X_{n}^x,X_{n}^y)\right)=\sum_{n=0}^{+\infty} P_\mu^n \phi_{\tilde v,I}(x,y) \ <+\infty .
\end{align*}
By Borel--Cantelli lemma, we get that for $\sigma$-a.e.~$(x,y)$, a.s., either for all large $n$, $[\tilde v X_n^x, \tilde v X_n^y]\supset I$ or, for all large $n$, $[\tilde v X_n^x, \tilde v X_n^y]\not\supset I$.

Under ${\bf M}_2(\overline{\mu})$, {\bf I} and {\bf P}, the process
$S_n:=\ln \vert A_n\cdots A_1\vert $ satisfies a CLT with positive variance
\cite{BQCLT}. Since $\lambda_1=0$, for every $M>0$,
$\P(S_n<-M \text{ i.o.})\geq\limsup_n\P(S_n<-M)=\frac12$, so that
$\P(\liminf_n S_n=-\infty)\geq\frac12$. Since the matrices $A_n$ are
invertible, this event does not depend on $A_1,\dots,A_k$ for any $k$; it
is therefore a tail event, and it has probability one by Kolmogorov's
$0$-$1$ law. Hence, almost surely,
$\liminf \vert (A_n\cdots A_1)(x-y)\vert =\liminf \vert X_n^x-X_n^y\vert=0$,
which implies that the interval $[\tilde v X_n^x,\tilde v X_n^y]$ has
arbitrarily small diameter infinitely often. Therefore the first
alternative cannot hold, so that
\[
\textrm{for $\sigma$-a.e.~$(x,y)$}\qquad \P(\cC_{(x,y)}^I)
=\P\big(I\not\subset [\tilde v X_n^x,\tilde v X_n^y]\ \text{for all large } n\big)=1.
\]

\medskip
\noindent\emph{Step 2: $\P(\bigcap_{I\subset\R} \cC_{(x,y)}^I)=1$ for $\sigma$-almost every pair $(x,y)\in\cZ_d$.}\\
First, observe that if $J\subset I$, then 
  $\cC_{(x,y)}^{J}\subset \cC_{(x,y)}^I$.
Since every compact non-degenerate interval contains a non-degenerate interval
$[a,b]$ with $a$ and $b$ rational, the above observation yields $\bigcap_{I\subset\R} \cC_{(x,y)}^I=\bigcap_{a\neq b\in\Q} \cC_{(x,y)}^{[a,b]}$. Thus by Step 1 and a countable intersection argument,  
for $\sigma$-a.e.~$(x,y)$
\begin{equation}
    \P(\bigcap_{I\subset\R} \cC_{(x,y)}^I)
=\P\big(\forall I\subset \R : I\not\subset \tilde v [X_n^x, X_n^y]\ \text{for all large } n\big)=1.\label{eq.crossing.strong2}
\end{equation}

\medskip \noindent\emph{Step 3: $\P\Bigl(\bigcap_{(x,y)\in \ell_0\times\ell_0}\bigcap_{I\subset\R} \cC_{(x,y)}^I\Bigr)=1$ for $\eta$-almost every affine line $\ell_0$.}\\
By the disintegration formula of $\sigma$ (Proposition~\ref{prop.invariant.sigma}) 
and Tonelli's theorem, for $\eta$-a.e.~$\ell_0\in \AG$, the set
\[
\mathcal{G}_{\ell_0} := \bigl\{(x,y)\in\cZ_1(\ell_0) \mid 
\P\bigl(\bigcap_{I\subset\R} \cC_{(x,y)}^I\bigr)=1\bigr\}
\]
has full $\lambda^{\ell_0}$-measure. Since $\lambda^{\ell_0}$ gives positive 
mass to every non-empty open subset of $\cZ_1(\ell_0)$, the set 
$\mathcal{G}_{\ell_0}$ is dense in $\cZ_1(\ell_0)$. As $\cZ_1(\ell_0)$ is 
separable, we may therefore fix a countable dense subset 
$D_{\ell_0}\subset \mathcal{G}_{\ell_0}$.

  By \eqref{eq.crossing.strong2},  for $\eta$-a.e.~$\ell_0$
\begin{align*}
 &\P\Bigl(\bigcap_{(x',y')\in D_{\ell_0}}\bigcap_{I\subset\R} \cC_{(x',y')}^I\Bigr)=1.
\end{align*}
On the other hand, for every $(x,y)\in \cZ_1(\ell_0)$, there exists $(x',y')\in D_{\ell_0}$ with $[x,y]\subset [x',y']$. 
Hence $
I\subset \tilde v[X_n^x,X_n^y]\ \Longrightarrow\ I\subset \tilde v[X_n^{x'},X_n^{y'}]
$ for any $n$, then
 $\cC_{(x',y')}^{I}\subset \cC_{(x,y)}^I$. This proves that for $\eta$-a.e.\ $\ell_0$,
\begin{equation}\label{eq.crossing.ae.line}\P\Bigl(\bigcap_{(x,y)\in\cZ_1(\ell_0)}\bigcap_{I\subset\R} \cC_{(x,y)}^I\Bigr)=1.\end{equation}
Since a single point contains no non-degenerate interval, this
remains true with $\cZ_1(\ell_0)$ replaced by $\ell_0\times\ell_0$.
 
\medskip
\noindent\emph{Step 4: $\P\Bigl(\bigcap_{(x,y)\in\R^d\times \R^d}\bigcap_{I\subset\R} \cC_{(x,y)}^I\Bigr)=1$.}\\
Fix a line $\ell_0=z_0+\R v_0$ satisfying \eqref{eq.crossing.ae.line}. 
We first prove that
there exists $\mathcal{L}_{v_0}\subset (\Aff(\R^d))^{\Z}$ of full 
$\P$-probability measure such that for every $\omega\in \mathcal{L}_{v_0}$ and  every $(x,y)\in \cZ_d$, there exist $\hat x=\hat x(\omega)\in \ell_0$ and $\hat y=\hat y(\omega)\in \ell_0$ such that 
\begin{equation}\label{eq.random.point.line}  A_n\cdots A_1(x-\hat x)\to 0 \qquad \textrm{and}\qquad A_n\cdots A_1(y-\hat y)\to 0.\end{equation}
Indeed, since $\lambda_1>\lambda_2$, Oseledets' theorem gives, for $\P$-almost every $\omega\in \Omega$, a random $\vvt_1 =\vvt_1 (\omega)\in {(\R^d)}^*\setminus \{0\}$ such that 
\begin{equation}\label{eq.oss2}
 \mbox{ for every $v\in \R^d$ with $\vvt_1 (\omega) v=0$,}    \  \lim_{n\to+\infty}\frac{1}{n}\ln\vert A_n\cdots A_1 v\vert\leq \lambda_2,
\end{equation}
hence $A_n\cdots A_1 v\to 0$, since $\lambda_2<0$ here.
Moreover, the distribution $\tilde \nu$ of $\vvt_1 $ is the unique right $\mu$-invariant probability measure on $\mathbf{P}((\R^d)^*)$ and it  is non-degenerate under hypothesis {\bf I} \cite{Furstenberg1963}. Therefore, the set $\cL_{v_0}$ of $\omega\in\Omega$ such that \eqref{eq.oss2} holds and  
$\vvt_1 (\omega)\, v_0\neq 0$ has $\P$-full measure. On this event, the line $\ell_0=z_0+\R v_0$ is not parallel to the hyperplane $\ker(\vvt_1 )$, hence it intersects every affine hyperplane $x+\ker(\vvt_1 )$ in a unique point $\hat{x}:= z_0+\frac{\vvt_1  (x-z_0)}{\vvt_1  v_0} v_0$.
\begin{figure}[h]
\centering
\figu{0.3}{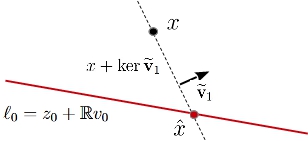}
\caption{Projection $\hat x$ of $x$ onto $\ell_0$ along $\ker\vvt_1(\omega)$.}
\label{fig.droite-proj}
\end{figure}
This implies that $\vvt_1 (x-\hat{x})=0$, yielding \eqref{eq.random.point.line}. The same result holds for $y$ and its projection $\hat y$.

\noindent We now conclude the proof of Step 4.
Let $\omega\in \cL_{v_0}\cap \Bigl(\bigcap_{(x,y)\in \ell_0\times\ell_0}\bigcap_{I\subset\R} \cC_{(x,y)}^I\Bigr)$.  Let $x,y\in \cZ_d$ and $I\subset \R$. We claim that $I\not\subset \tilde v[X_n^x(\omega),X_n^y(\omega)]$ for large $n$ (i.e. $\omega \in \cC_{(x,y)}^I$). 
Indeed take a smaller interval  $J\subset \operatorname{Int}(I)$ with non-empty interior.  Let $(\hat x,\hat y)\in \ell_0$ be as in \eqref{eq.random.point.line}. Using $X_n^x-X_n^{\hat x}=A_n\cdots A_1(x-\hat x)$, we have 
 \[
\tilde v(X_n^x - X_n^{\hat x}) \to 0 \qquad \text{and} \qquad \tilde v(X_n^y - X_n^{\hat y}) \to 0,
\]
so the endpoints of the projected segments $\tilde v[X_n^x,X_n^y]$ and $\tilde v[X_n^{\hat x},X_n^{\hat y}]$ become arbitrarily close. Since $J\subset \operatorname{Int}(I)$, it follows that for all sufficiently large $n$,
\[
I\subset \tilde v[X_n^x,X_n^y]
\ \Longrightarrow\
J\subset \tilde v[X_n^{\hat x},X_n^{\hat y}].
\]
By Step 3 applied to $(\hat x,\hat y)$ and $J$, the latter inclusion fails for all large $n$, hence so does the former.  In particular 
$$\cL_{v_0}\cap \Bigl(\bigcap_{(x,y)\in\ell_0\times\ell_0}\bigcap_{I\subset\R} \cC_{(x,y)}^I\Bigr)\subseteq \bigcap_{(x,y)\in\R^d\times \R^d}\bigcap_{I\subset\R} \cC_{(x,y)}^I$$
and since the left term has full measure this concludes the proof.
\end{proof}

We now complete the proof of Theorem~\ref{theorem.main.contraction} 
using a compactness argument. 
\begin{proof}[of Theorem~\ref{theorem.main.contraction}]
Fix $\tilde v \in (\R^d)^*\setminus\{0\}$ and a bounded interval $I=[a,b]\subset \R$. 
Let $\varepsilon>0$. 
By compactness of $I$, there exist finitely many subintervals 
$I_0,\ldots, I_{m-1}\subset \R$ such that
\begin{equation}\label{eq.cross}
s\in I,\ t\in \R,\ |s-t|\geq \varepsilon 
\,\Longrightarrow\,
\exists\,0\leq i \leq m-1\ \text{such that }  [s,t]\supset I_i.
\end{equation}
Indeed, choose for instance points
$
a_0<a=a_1<a_2<\cdots<a_{m-1}=b<a_m
$ 
with $a_{j+1}-a_j<\varepsilon/2$ for $0\leq j\leq m-1$, and set $I_j=[a_j,a_{j+1}]$. Observe that, if $s \in I$ then there exists $1\leq i\leq m-2$ such that $s\in I_i$; now, if $\vert s-t\vert \geq \varepsilon$, then either $I_{i-1}\subset [s, t]$ or $I_{i+1}\subset [s, t]$.
Theorem~\ref{theorem.main.crossing} applied to each $I_j$ implies that, for $\P$-almost every $\omega$ and any $x,y\in \R^d$, there exists $n_0=n_0(\omega, x,y)\in \N$ such that $
 [\tilde v X_n^x,\tilde v X_n^y]\not\supset I_i
$ for every $n\geq n_0$ and every $0\leq i\leq m-1$,
Applying \eqref{eq.cross} to $s=\tilde v(X_n^x)$ and $t=\tilde v(X_n^y)$, we deduce that almost surely, for any $x,y\in \R^d$ and every $n\geq n_0$,
\[
\tilde v(X_n^x)\in I \ \Longrightarrow\ |\tilde v(X_n^x - X_n^y)|<\varepsilon.
\]
Intersecting over $\varepsilon\in \Q_{>0}$, we obtain
$\tilde v(X_n^x - X_n^y)\,\mathds{1}_{\tilde v(X_n^x)\in I}\longrightarrow 0$ a.s., proving \eqref{eq.local.contraction.y}. Then   
\eqref{eq.local.contraction.compact}follows by applying \eqref{eq.local.contraction.y} to  a basis of $(\R^d)^*$.
\end{proof}

 \section{Total dissipativity of the joint process}\label{section.dissipative}

As announced, we prove in this section  the total dissipativity on $L^1(\cZ_d,\sigma)$ of the transition operator $P_\mu$ of the diagonal affine process:  Theorem
\ref{theorem.affine.dissipative}. 
 
In dimension $d=1$, this statement follows from Guivarc'h--Keane--Roynette's theorem  \cite{Guivarch-Keane-Roynette}  since $\cZ_1\simeq \Aff(\R)$ and $\Aff(\R)$ is a non-unimodular group. As explained in Section \ref{sec.diag.action}, when $d \geq 2$, the space $\cZ_d$ is only an $\Aff(\R^d)$-homogeneous space and \cite{Guivarch-Keane-Roynette}'s result cannot be applied.

The key point in Guivarc'h--Keane--Roynette's theorem is that the right action in $G$ commutes with the left action, hence with the left random walk; nevertheless, by non-unimodularity of $\Aff(\R)$,  it 
   does not preserve the left Haar measure (which is $\mu$-invariant). On $\cZ_d$, we replace the right action  with the   transformations $\Phi_\lambda(x,v)=(x,\lambda v)$ that commute with the left action, while not preserving the $\mu$-invariant measure $\sigma$  when $|\lambda|\neq1$.

Our proof uses many ideas of \cite{Guivarch-Keane-Roynette}, but in a
different framework. Our
criterion is purely measure-theoretic: it relies on the ergodic theory of
positive contractions of $L^1$ (namely the Hopf decomposition and the
Chacon--Ornstein ratio ergodic theorem \cite{Foguel,Garsia,Revuz}), rather than
on the group and its topological structure used in
\cite{Guivarch-Keane-Roynette}.

We first develop in Theorem~\ref{theorem.main.dissipativity} an abstract
framework underlying this argument, for semigroup actions admitting a
suitable commuting transformation. It recovers Guivarc'h--Keane--Roynette's
theorem as a special case (Corollary~\ref{cor.GKR}).

\medskip 

\subsection{An abstract dissipativity criterion}\label{the.criterion}
Throughout this subsection, $(\Gamma,\mu)$
is a measurable semigroup with $\mu$ a probability measure on $\Gamma$,
and $(\gamma_n)_{n\geq1}$ is a sequence of i.i.d.\ random variables with
law $\mu$. Let $\cZ $ be a measurable space carrying a measurable action of
$\Gamma$, and denote by
$P_\ell f(\fz)=\int f(\gamma\fz)\,d\mu(\gamma)$ the associated Markov
operator. Let $\sigma_\ell$ be a  non-zero $\sigma$-finite
measure on $\cZ$  $\mu$-invariant, i.e.\
$P_\ell$-invariant: $\int P_\ell f\,d\sigma_\ell=\int f\,d\sigma_\ell$ for
every non-negative measurable $f$. 
Then $P_\ell$ is a positive contraction on $L^1(\cZ ,\sigma_\ell)$.

 \medskip
\noindent\textbf{The companion space.}
In this section, we call \emph{observable} any non-negative measurable function
$\chi:\cZ\to\R_+$. The semigroup $\Gamma$ acts on the right on observables:
for an observable $\chi$ and $\gamma\in\Gamma$, set
\[
(\chi \gamma)(\fz):=\chi(\gamma \fz).
\]
Let $\cZ ^*$ be a measurable space of  observables  such that:
\begin{itemize}
\item $\cZ ^*\subset L^1(\cZ ,\sigma_\ell)$;
\item   
the map $\cZ\times\cZ^*\to\R_+$, $(\fz,\chi)\mapsto \chi(\fz)$, is
measurable with respect to the product $\sigma$-algebra;
\item $\cZ ^*$ is stable under the right action, and
the map $\cZ^*\times\Gamma\to\cZ^*$, $(\chi,\gamma)\mapsto \chi \gamma$,
is measurable with respect to the product $\sigma$-algebra.
\end{itemize}
The space $\cZ^*$ should be thought of as a family of integrable
observables stable under the dynamics; it provides a substitute for the
missing right action on $\cZ$.
This induces a Markov operator $P_r$ on $\cZ ^*$, acting on non-negative
measurable functions $F:\cZ^*\to\R_+$ by
\[
P_rF(\chi)=\int_\Gamma F(\chi \gamma)\,d\mu(\gamma).
\]
   For  any $\fz\in\cZ $, the \emph{evaluation map}
 $\ev_\fz:\cZ^*\to\R_+$, $\ev_\fz(\chi):=\chi(\fz)$, satisfies
   the duality relation
   \begin{equation}\label{eq.duality}
   P_\ell(\chi)(\fz)=P_r(\ev_\fz)(\chi),\qquad \chi\in\cZ^*,
   \end{equation}
which propagates to all iterates: for any $n \geq 1$, 
\begin{equation}\label{eq.duality.propagate}
P_\ell^n\chi(\fz)
=\E\bigl(\chi(\gamma_n\cdots \gamma_1 \fz)\bigr)
=\E\bigl((\chi \gamma_1\cdots \gamma_n)(\fz)\bigr)
=P_r^n(\ev_\fz)(\chi).
\end{equation}
This duality identity allow us to move from left to right random walks. 

\medskip

\noindent\textbf{Standing assumptions.}
We assume the following:
\begin{enumerate}
\item[(A1)] There exists a  non-zero $\sigma$-finite measure $\sigma_r$
on $\mathcal Z^*$ which is $P_r$-invariant, i.e.\
$\int P_rF\,d\sigma_r=\int F\,d\sigma_r$ for every non-negative measurable
$F$, and such that $L^1(\mathcal Z^*, \sigma_r)$ is separable.
\item[(A2)] For every $\fz\in\cZ $, the evaluation map $\ev_\fz$ belongs to
$L^1(\cZ ^*,\sigma_r)$.
\item[(A3)] (Weak irreducibility) For $(\sigma_\ell\otimes\sigma_r)$-a.e.\
$(\fz,\chi)$, $\sum_{n\ge 0} P_\ell^n \chi(\fz)>0$.
\end{enumerate}
By (A1), $P_r$ is a positive contraction on $L^1(\cZ ^*,\sigma_r)$.
Assumption (A3) ensures that almost every observable detects almost
every point along the $P_\ell$-orbit.

\medskip
\noindent\textbf{An intrinsic metric on $\cZ $.}
Assumption (A2) allows us to define a pseudo-metric on $\cZ $:
\begin{equation}\label{eq.dr}
d_r(\fz,\fz'):=\|\ev_{\fz}-\ev_{\fz'}\|_{L^1(\cZ ^*,\sigma_r)}
=\int_{\cZ ^*} |\chi(\fz)-\chi(\fz')|\,d\sigma_r(\chi).
\end{equation}
Since 
$L^1(\cZ ^*,\sigma_r)$ is separable, the pseudo-metric space $(\cZ ,d_r)$ is also
separable.

Let $\Phi:\cZ \to\cZ $ be measurable. We say that the action of
$(\Gamma,\mu)$ is \emph{$\Phi$-Lipschitz minimal} if, for every
$d_r$-Lipschitz map $f:\cZ \to\R$,
\begin{equation}\label{def.lipshitz}
f(\gamma \fz)=f(\fz)
\ \ \text{for }(\mu\otimes \sigma_\ell)\text{-a.e. }(\gamma,\fz)
\quad\Longrightarrow\quad
f\circ \Phi=f
\ \ \sigma_\ell\text{-a.e.}
\end{equation}
This weak form of minimality with respect to $d_r$ replaces the
irreducibility (adaptedness) assumption of Guivarc'h--Keane--Roynette's
theorem; it holds for instance when the orbits of the semigroup generated by
$\mathrm{supp}(\mu)$ are $d_r$-dense.

\medskip

\noindent\textbf{Conservative/dissipative decomposition.}
For a positive contraction $P$ on $L^1(\mathcal Y,\sigma)$, with
$\sigma$ a $\sigma$-finite measure on $\mathcal Y$, there exist measurable sets $C$ and $D$, with
$\mathcal Y=C\sqcup D$, such that for every $f\in L^1_+(\mathcal Y,\sigma)$,
\[
\sum_{n\ge 0} P^n f =
\begin{cases}
0 \text{ or } +\infty & \textrm{$\sigma$-a.e.~on $C$},\\
<+\infty & \textrm{$\sigma$-a.e.~on $D$}
\end{cases}
\]
(Hopf decomposition, see for instance \cite[Ch.~II]{Foguel}).
We say that $P$ is totally conservative if $\sigma(D)=0$, and totally
dissipative if $\sigma(C)=0$.

The main result of this section is the following. 
\begin{theo}\label{theorem.main.dissipativity}
Let $\Gamma$ act measurably on $\cZ $, let $\mu$ be a probability measure on $\Gamma$, and let $\sigma_\ell$ be a non-zero $\mu$-invariant $\sigma$-finite measure on $\cZ $. Assume that there exists a measurable space $\cZ ^*$ of observables satisfying (A1)--(A3).
Let $\Phi:\cZ \to\cZ $ be a measurable map commuting with the action:
\begin{equation}\label{eq.commute.phi}
\Phi(\gamma  \fz)=\gamma  \Phi(\fz)
\qquad (\mu\otimes\sigma_\ell)\text{-a.e.}
\end{equation}
Assume that the action is $\Phi$-Lipschitz minimal and that, for $\sigma_r$-a.e.\ $\chi$, the function $\fz\mapsto \chi(\Phi(\fz))$ belongs to $L^1(\cZ ,\sigma_\ell)$.

If $P_\ell$ is not totally dissipative, then for $\sigma_r$-a.e.\ $\chi$,
\[
\int_{\cZ }\chi(\Phi(\fz))\,d\sigma_\ell(\fz)
=
\int_{\cZ }\chi(\fz)\,d\sigma_\ell(\fz).
\]
\end{theo}
In other words, if $\Phi$ does not preserve $\sigma_\ell$ when tested
against observables, then $P_\ell$ is totally dissipative. This provides the
dissipativity criterion used in Corollary~\ref{cor.GKR} and in
Section~\ref{sec.affine.application}.

We first show shows that each of the Markov operators
$P_\ell$ and $P_r$ is either totally conservative or totally dissipative
(there is no partial conservativity), and that both have the same type. This can be seen as a generalization of  the classical fact that the
left and right random walks on a group are either both recurrent or both transient.
\begin{prop}\label{prop.dichotomy}
In the above setting, assume (A1)--(A3). Then the operators $P_\ell$ on $L^1(\cZ ,\sigma_\ell)$ and
$P_r$ on $L^1(\cZ ^*,\sigma_r)$ are either both totally conservative or both
totally dissipative.
\end{prop}

\begin{proof}
Let $\cZ =C_\ell\sqcup D_\ell$ and $\cZ ^*=C_r\sqcup D_r$
be the conservative-dissipative decompositions of the positive contractions
$P_\ell$ and $P_r$ on $L^1(\cZ ,\sigma_\ell)$ and $L^1(\cZ ^*,\sigma_r)$.

Since $\sil$ and $\sir$ are non-zero, it is enough to show that if $P_\ell$ has a dissipative part
(i.e.~$\sil(D_\ell)>0$), then $P_r$ is totally dissipative
(i.e.~$\sir(C_r)=0$), and symmetrically with the roles of $\ell$ and $r$
reversed. 

Assume $\sigma_\ell(D_\ell)>0$. By definition of $\cZ ^*$, every
$\chi\in\cZ ^*$ belongs to $L^1(\cZ ,\sigma_\ell)$. It follows from the
definition of the dissipative part $D_\ell$ that for every $\chi\in \cZ ^*$
and for $\sigma_\ell$-a.e.~$\fz\in D_\ell$,
\[
\sum_{n\geq 0}P_\ell^n\chi(\fz)<\infty.
\]
By (A3), for $(\sigma_\ell\otimes\sigma_r)$-a.e.\ $(\fz,\chi)$,
\[
0<\sum_{n\geq 0}P_\ell^n\chi(\fz).
\]
Using the iterated duality relation \eqref{eq.duality.propagate}
and Tonelli's theorem, we deduce that for $\sigma_r$-a.e.\ $\chi$ and
$\sigma_\ell$-a.e.\ $\fz\in D_\ell$,
\[
0<\sum_{n\geq 0}P_r^n(\ev_\fz)(\chi)<\infty.
\]
By Tonelli's theorem again, since $\sil(D_\ell)>0$, we deduce that there
exists $\fz\in D_\ell$ such that
\[
0<\sum_{n\geq 0}P_r^n(\ev_\fz)(\chi)<\infty
\qquad\textrm{for $\sigma_r$-a.e.~$\chi$}.
\]
By (A2), $\ev_\fz\in L^1(\cZ ^*,\sigma_r)$. Hence the Hopf decomposition for
$P_r$ implies that, for $\sir$-a.e.~$\chi$, $\chi\in D_r$. In other words
$\sir(D_r^c)=0$ as claimed.
The converse implication is proved in the same way, exchanging the roles of
$(\cZ ,\sigma_\ell,P_\ell)$ and $(\cZ ^*,\sigma_r,P_r)$ using the duality
relation.
\end{proof} 
We shall use the following fact: a non-negative superharmonic function
for a totally conservative Markov operator is constant along trajectories.
Such a fact is folklore in the theory of Markov chains  but we could not find a reference for the  purely measure theoretical statement we use.

\begin{lemma}\label{lem.subharmonic.constant}
Let $(\mathcal Y,\sigma)$ be a $\sigma$-finite measure space and let $P$ be a
Markov operator on $\mathcal Y$ such that $\sigma$ is $P$-invariant. Assume
that $P$ is totally conservative on $L^1(\mathcal Y,\sigma)$. Let
$h:\mathcal Y\to[0,+\infty)$ be measurable and finite $\sigma$-a.e. Then the
following are equivalent:
\begin{enumerate}
\item $Ph\le h$ $\sigma$-a.e.;
\item $Ph=h$ $\sigma$-a.e.;
\item for $\sigma$-a.e. $x$ and for $P(x,\cdot)$-a.e. $y$, one has $h(y)=h(x)$.
\end{enumerate}
\end{lemma}
\begin{proof}{}
The directions $3\Rightarrow 2\Rightarrow 1$ are immediate. It remains to show $1\Rightarrow 2$
and $2\Rightarrow 3$.

\medskip
\noindent
\textit{Step 1: $1\Rightarrow 2$.}
Let $g:=h-Ph\geq 0$. Since $Ph\leq h$, we have $P^kh\leq h<+\infty$
$\sigma$-a.e.\ for every $k$, so that, for every $n\in \N$,
$$\sum_{k=0}^{n-1}P^k g=h-P^n h.$$
Since $P$ is a positive operator and $h\geq 0$, we have $P^n h\geq 0$ so that
$$\forall n\in \N, 0 \leq \sum_{k=0}^{n-1}P^k g \leq h.$$
Hence, for $\sigma$-a.e.\ $y$,
$$0\leq \sum_{k=0}^{+\infty}P^k g(y)\leq h(y)<+\infty.$$
To use conservativity we truncate $g$.
Let $(A_m)$ be an increasing sequence of finite-measure sets with
$\bigcup_m A_m=\mathcal{Y}$ and define
$g_m:=\min\{g \mathds{1}_{A_{m}}, m\}$.
Then $0\leq g_{m}\leq g$ and
$\|g_{m}\|_{L^1(\sigma)}\leq m\, \sigma(A_{m})<+\infty$, so
$g_{m}\in L_+^1(\mathcal{Y},\sigma)$.
Since $P$ is totally conservative,
$$\sum_{k\geq 0}P^k g_m \in \{0,+\infty\}\quad \sigma\textrm{-a.e.}$$
But $\sum_{k\geq 0}
P^k g<+\infty$ $\sigma$-a.e. Hence for every
$m\in \N$, $\sum_{k\geq 0}P^k g_{m}=0$, so $g_m=0$
$\sigma$-a.e. Since $g_m\uparrow g$ $\sigma$-a.e., we obtain $g=0$
$\sigma$-a.e., i.e.\ $Ph=h$ $\sigma$-a.e.

\medskip
\noindent
\textit{Step 2: $2\Rightarrow 3$.}
Since $z\mapsto \sqrt{z}$ is concave, Jensen's inequality for the Markov
operator $P$ yields
\begin{equation}\label{eq.jensen}
h^{\frac{1}{2}}=(Ph)^{\frac{1}{2}}\geq Ph^{\frac{1}{2}}
\qquad \sigma\textrm{-a.e.}
\end{equation}
Hence $h^{\frac{1}{2}}$ is $P$-superharmonic, thus harmonic by the
implication $1\Rightarrow 2$.
Consequently equality holds in \eqref{eq.jensen}.
Since the function $z\mapsto \sqrt{z}$ is strictly concave, equality in
Jensen's inequality occurs only when $h(z)$ is constant for
$P(y,dz)$-almost every $z$. By harmonicity, this constant must be $h(y)$.
This finishes the proof.
\end{proof}

Before proving Theorem \ref{theorem.main.dissipativity} we recall the
Chacon--Ornstein ratio theorem in the form we will use (see for instance
\cite[Theorem 2.6.1]{Garsia}, \cite[Ch. III, Theorem D]{Foguel},
\cite[Theorem 3.3]{Revuz}).
Let $(\mathcal{Y},\sigma)$ be a $\sigma$-finite measure space and let
$P:L^1(\mathcal{Y},\sigma)\to L^1(\mathcal{Y},\sigma)$ be a totally
conservative positive contraction.
For $f,g\in L^1(\mathcal{Y},\sigma)$ with $\sum_{k=0}^{+\infty} P^k g > 0$
$\sigma$-a.e., define (for $n$ large enough)
$$
R_n(f/g)(y)
=\frac{\sum_{k=0}^{n-1}P^k f(y)}{\sum_{k=0}^{n-1}P^k g(y)} .
$$
Then $R_n(f/g)$ converges $\sigma$-a.e.~to a $\sigma$-a.e.~finite limit.
Denote the limit by $$R(f/g):=\lim_{n\to \infty}R_n(f/g).$$
Moreover the map $f\mapsto R(f/g)\,g$ is a positive linear contraction
of $L^1(\mathcal{Y},\sigma)$.
Finally
\begin{equation}\label{eq.mass.principle.total}
\int R(f/g)g\, d\sigma = \int f\, d\sigma .
\end{equation}

\begin{proof}[of Theorem \ref{theorem.main.dissipativity}]{}
Throughout this proof, given a measure $m$, we write ``$\forall_m x$''
for ``for $m$-almost every $x$''.
Assume that $P_\ell$ has a conservative part.
By Proposition \ref{prop.dichotomy}, both $P_\ell$ and $P_r$ are totally
conservative as operators on $L^1(\cZ ,\sil)$ and $L^1(\cZ ^*,\sir)$
respectively.

In all the proof, we fix $\fz_0\in \cZ $ such that for $\sir$-almost every
$\chi$, $\sum_{n=0}^{+\infty}{P_\ell^n\chi(\fz_0)}=+\infty$. This is possible by
 assumption (A3), total conservativity of
$P_\ell$ and Tonelli's theorem.

For every $\fz\in \cZ $,
one can then define for $\sir$-almost every $\chi$ and for every $n$ large
enough, say $n\geq n_0(\chi)$,
\begin{equation}\label{eq.expression}\mathcal{R}_n(\fz,\chi):=\frac{\sum_{k=0}^{n-1} P_r^k \ev_{\fz}(\chi)}
{\sum_{k=0}^{n-1}P_r^k\ev_{\fz_0}(\chi)}=
\frac{\sum_{k=0}^{n-1}{P_\ell^k \chi(\fz)}}{\sum_{k=0}^{n-1}{P_\ell^k \chi(\fz_0)}},\end{equation}
where the equality follows from duality \eqref{eq.duality.propagate}. The proof proceeds in four steps.
First, using \ Chacon--Ornstein
ratio theorem for $P_r$, we show that
$\mathcal{R}_n(\fz,\chi)$ converges to a non-degenerate limit $\mathcal{R}(\fz,\chi)$.
Second, we prove that the function $\mathcal{R}(\cdot, \chi)$ is constant
along the $P_\ell$-trajectories.
Then, the $\Phi$-Lipschitz minimality assumption implies that this function
is $\Phi$-invariant.
Finally, we conclude by using the commutation property of $\Phi$ and by
applying the Chacon--Ornstein theorem to $P_\ell$.

\medskip\noindent\textit{Step 1:} We show that the ratio limit exists
and is finite: 
\begin{equation}\label{eq.step1}
\forall \fz\in\cZ, \forall_\sir \chi, \quad
\mathcal{R}(\fz,\chi):=\lim_n \mathcal{R}_n(\fz,\chi)\in[0,+\infty),
\end{equation}
and   $\mathcal{R}(\fz,\xi)>0$ for $\sil\otimes\sir$-a.e.\ $(\fz,\chi)$.

By \eqref{eq.expression}, $\mathcal R_n(\fz,\cdot)=R_n(\ev_\fz/\ev_{\fz_0})$ is
the Chacon--Ornstein ratio of the pair $(\ev_\fz,\ev_{\fz_0})$ for the
operator $P_r$. Since for every $\fz\in \cZ $, $\ev_{\fz}$ and $\ev_{\fz_0}$
belong to $L^1(\cZ ^*,\sir)$ thanks to assumption (A2), and the denominator
is $>0$ for $\sir$-a.e.~$\chi$ by choice of $\fz_0$, the Chacon--Ornstein
theorem applied to the totally conservative operator $P_r$ yields
\eqref{eq.step1}, with $\mathcal R(\fz,\cdot)=R_r(\ev_\fz/\ev_{\fz_0})$.

For the positivity, note that by assumption~(A3), the
duality relation and Tonelli's theorem, for $\sigma_\ell$-a.e.~$\fz$ the
Chacon--Ornstein theorem applies to both pairs $(\ev_{\fz_0},\ev_{\fz})$ and
$(\ev_{\fz},\ev_{\fz_0})$, whose finite-stage ratios are reciprocal, hence so
are their limits. Since, by the Chacon--Ornstein theorem, these limits are
finite $\sigma_r$-a.e., neither can vanish, and therefore
$\mathcal{R}(\fz,\chi)>0$ $\sigma_\ell\otimes\sigma_r$-a.e.

\medskip\noindent \textit{Step 2:} We show that 
the functions $\fz\mapsto \mathcal{R}(\fz,\chi)$ are constant along the trajectories of the Markov chain, namely 
\begin{equation}\label{eq.alpha.constant}\forall_{\sir}\ \chi, \forall_\sil\ \fz, \forall_{P_\ell(\fz,\cdot)}\ \fz', \,\mathcal{R}(\fz,\chi)=\mathcal{R}(\fz',\chi).
\end{equation}
We first prove that for $\sir$-a.e. $\chi$, the function $\mathcal{R}(\cdot, \chi)$ is $P_\ell$-superharmonic. 
For $\chi\in \cZ ^*$, let
$$A_\chi:=\{\fz\in \cZ ; \mathcal{R}_n(\fz,\chi)\to \mathcal{R}(\fz,\chi)\}.$$
By \eqref{eq.step1} and Tonelli, for $\sir$-a.e.~$\chi$, $\sil(A_\chi^c)=0$. From now on fix such a generic $\chi\in \cZ ^*$ verifying in addition $\sum_{k=0}^{+\infty}{P_\ell^k \chi(\fz_0)}=+\infty$. 
Since $\sil$ is $P_\ell$-invariant and $\sil(A_\chi^c)=0$, 
we have
\begin{equation}\label{eq.abs.cont.almost}\forall_\sil \fz\in A_\chi, \,P_\ell(\fz,A_\chi)=1.\end{equation}
Therefore, for such a $\sil$-generic 
$\fz$, we have:  
\begin{eqnarray}
    P_\ell
   \mathcal{R}(\fz,\chi)& = & \int{\mathcal{R}(\fz',\chi)\, P_\ell(\fz,d\fz')}\nonumber\\
    & = & \int{\lim_n \mathcal{R}_n(\fz',\chi)\, P_\ell(\fz,d\fz')}\quad \textrm{(by \eqref{eq.step1} and \eqref{eq.abs.cont.almost})}\nonumber\\
    & \leq  & \liminf_n P_\ell \mathcal{R}_n(\cdot,\chi)(\fz)\quad \textrm{(Fatou)}\nonumber\\
    & = & \liminf_n\left[\mathcal{R}_n(\fz,\chi) - \frac{\chi(\fz)}{\sum_{k=0}^{n-1}{P_\ell^k \chi(\fz_0)}} + \frac{P_\ell^n \chi(\fz)}{\sum_{k=0}^{n-1}{P_\ell^k \chi(\fz_0)}}\right]\nonumber\\
    & = & \mathcal{R}(\fz,\chi) + \liminf_n \frac{P_\ell^n \chi(\fz)}{\sum_{k=0}^{n-1}{P_\ell^k \chi(\fz_0)}},\label{eq.harmonic.1}
\end{eqnarray}
It remains to check that 
\begin{equation}\label{eq.fatou2}\forall \sil \fz, \,\liminf_n \frac{P_\ell^n\chi(\fz)}{\sum_{k=0}^{n-1}{P_\ell^k \chi(\fz_0)}}=0.\end{equation}
By Fatou's lemma and the $P_\ell$-invariance of $\sil$, 
$$\int{\liminf_n P_\ell^n \chi(\fz) \, d\sil(\fz)}
\leq  \liminf_n\int{P_\ell^n \chi(\fz)\, d\sil(\fz)}
=  \int{\chi(\fz)\,d\sil(\fz)}<+\infty.
$$
Hence 
$$
\forall_\sil \fz, \,\liminf_n P_\ell^n \chi(\fz)
<+\infty.
$$
   Since $\sum_{k=0}^{+\infty}{P_\ell^k \chi(\fz_0)}=+\infty$,  the ratio
   in \eqref{eq.fatou2} tends to $0$ along any  subsequence on which
   $(P_\ell^n\chi(\fz))_n$ stays bounded, and \eqref{eq.fatou2} follows. Plugging \eqref{eq.fatou2} 
in \eqref{eq.harmonic.1} yields
$$P_\ell \mathcal{R}(\cdot,\chi) \leq \mathcal{R}(\cdot,\chi) \quad (\textrm{$\sil$-a.e.}).$$
Assertion \eqref{eq.alpha.constant} follows from Lemma \ref{lem.subharmonic.constant}. 

\medskip\noindent\textit{Step 3:} We show the $\Phi$-invariance of
$\mathcal{R}(\cdot,\chi)$, namely:
\begin{equation}\label{eq.alpha.commute}
\forall_\sil \fz, \forall_\sir \chi, \,\mathcal{R}(\fz,\chi)=\mathcal{R}(\Phi(\fz),\chi).\end{equation}
Assume first the following assumption which is stronger than (A3):
$\chi(\fz)>0$ for $\sil\otimes \sir$-a.e.~$(\fz,\chi)$.
 The Chacon--Ornstein theorem applied
to the totally conservative operator $P_r$ ensures that, for every
$\fz\in \cZ $, the function
$$L_{\fz}(\chi):=R_r(\ev_{\fz}/\ev_{\fz_0})(\chi)\,\ev_{\fz_0}(\chi)$$
belongs to $L^1(\cZ ^*,\sir)$.
The same theorem ensures that for every $f,g\in L^1(\cZ ^*,\sir)$ with $g>0$
$\sir$-a.e., the map $f\mapsto R_r(f/g) g$ is a contraction of
$L^1(\cZ ^*,\sir)$.
Applying this with $f_{\fz}=\ev_{\fz}$ and $g=\ev_{\fz_0}$ gives for every
$\fz_1,\fz_2\in \cZ $,
\begin{equation}\label{eq.isometric}\|L_{\fz_1}-L_{\fz_2}\|_{L^1(\cZ ^*,\sir)}\leq \|\ev_{\fz_1}-\ev_{\fz_2}\|_{L^1(\cZ ^*,\sir)}=d_r(\fz_1,\fz_2).\end{equation}
Fix $\fz'\in \cZ $, and denote
$$f_{\fz'}: \cZ \to \R, \fz\mapsto \|L_{\fz} - L_{\fz'}\|_{L^1(\cZ ^*,\sir)}.$$
By \eqref{eq.isometric} and the triangle inequality,
this function is $d_r$-Lipschitz.
Moreover, by \eqref{eq.alpha.constant}, we have that
$$\forall_\sil \fz, \forall_{P_\ell(\fz,\cdot)} \fz'', \,f_{\fz'}(\fz)=f_{\fz'}(\fz'').$$
By $\Phi$-Lipschitz minimality, this implies
\begin{equation}\label{eq.invariance.shift} \forall_\sil \fz, \, f_{\fz'}(\fz)=f_{\fz'}(\Phi(\fz)) \quad \textrm{i.e.}\quad \|L_{\fz}-L_{\fz'}\|_1=\|L_{\Phi(\fz)}-L_{\fz'}\|_1.\end{equation}
Since $(\cZ ,d_r)$ is separable, we choose
a countable dense subset $\{\fz_n; n\in \N\}$ in $(\cZ , d_r)$. Intersecting
the corresponding conull sets from \eqref{eq.invariance.shift} gives
\begin{equation}\label{eq.33}\forall_\sil \fz, \forall n\in \N, \|L_{\fz}-L_{\fz_n}\|_1=\|L_{\Phi(\fz)}-L_{\fz_n}\|_1.\end{equation}
Let $\fz\in \cZ $ be $\sil$-generic and take
$\fz_{n_k}\underset{k\to \infty}{\to} \fz$ in $(\cZ ,d_r)$. By
\eqref{eq.isometric}, $L_{\fz_{n_k}}\to L_{\fz}$ in $L^1(\cZ ^*,\sir)$,
hence letting $k\to \infty$ in \eqref{eq.33} shows that
$L_{\fz_{n_k}}\to L_{\Phi(\fz)}$ as well, which yields
$$\forall_\sil \fz, \forall_\sir \chi, \,L_{\fz}(\chi)=L_{\Phi(\fz)}(\chi).$$
Since, under the stronger assumption, we may also choose $\fz_0$ such
that $\ev_{\fz_0}(\chi)>0$ for $\sir$-a.e.~$\chi$, the definition of
$L_{\fz}(\chi)$ implies \eqref{eq.alpha.commute} under the stronger
assumption.

Now we show \eqref{eq.alpha.commute} under only (A3). Let $P:=P_r$ and fix
$m\in \N$ and denote by $P^m:=P_r^m$. Define
$f_{m,\fz}=\sum_{k=0}^{m-1}P^k\ev_{\fz}$, $g_m=\sum_{k=0}^{m-1}P^k \ev_{\fz_0}$.
Since $P$ is totally conservative, so is $P^m$.
Hence the Chacon--Ornstein theorem applies to the operator $P^m$
with the pair $(f_{m,\fz},g_m)$. Observe that
$R_r(\ev_{\fz}/\ev_{\fz_0})=R_r^m(f_{m,\fz}/g_m)$, where $R_r^m$ is the
Chacon--Ornstein ratio associated to the Markov operator $P^m$.
Define $L_{\fz}^{(m)}=R_r^m(f_{m,\fz}/g_m) \,g_m$.
Repeating the same proof as above, with \eqref{eq.isometric} replaced by
$\|L^{(m)}_{\fz_1}-L^{(m)}_{\fz_2}\|_1\leq m\,d_r(\fz_1,\fz_2)$, shows that for
$\sil\otimes \sir$-a.e.~$(\fz,\chi)$, $L_{\fz}^{(m)}(\chi)=L_{\Phi(\fz)}^{(m)}(\chi)$.
Hence, for every $m\in \N$ and for $\sil\otimes \sir$-a.e.~$(\fz,\chi)$, if
$g_{m}(\chi)>0$ then $\mathcal{R}(\fz,\chi)=\mathcal{R}(\Phi(\fz),\chi)$.
Intersecting these countable conull sets and using the choice of $\fz_0$ (which ensures
$\bigcup_m\{g_m>0\}=\{\sum_kP_r^k\ev_{\fz_0}>0\}$ has full $\sir$-measure),
we conclude that \eqref{eq.alpha.commute} holds for
$\sigma_\ell\otimes\sigma_r$-a.e.~$(\fz,\chi)$.

\medskip\noindent\textit{Step 4: End of the proof.}
Combining \eqref{eq.step1},   the
positivity of $\mathcal R$ (Step 1) and \eqref{eq.alpha.commute}, we obtain
by Tonelli's theorem for $\sir$-a.e.~$\chi$,
$$
\forall_\sil \fz, \, \lim_n \frac{\sum_{k=0}^{n-1} P_\ell^k \chi(\fz)}{\sum_{k=0}^{n-1} P_\ell^k \chi(\fz_0)}
=\lim_n \frac{\sum_{k=0}^{n-1} P_\ell^k \chi(\Phi(\fz))}{\sum_{k=0}^{n-1} P_\ell^k \chi(\fz_0)}>0.
$$
(Here we use that \eqref{eq.step1} holds for every $\fz\in\cZ$, in
particular at the point $\Phi(\fz)$.)
Dividing the two limits yields for $\sir$-a.e.~$\chi$,
\begin{equation}\label{eq.ratio.1}\forall_{\sil} \fz, \,\lim_n
\frac{\sum_{k=0}^{n-1} P_\ell^k \chi(\Phi(\fz))}{\sum_{k=0}^{n-1} P_\ell^k \chi(\fz)}=1.\end{equation}
On the other hand, since  $\Phi$ commutes with the $\mu$-action \eqref{eq.commute.phi},  and since $\sil$ is $\mu$-invariant, 
we have
$\Phi(\gamma_k\cdots\gamma_1\fz)=\gamma_k\cdots\gamma_1\Phi(\fz)$ almost surely
for $\sil$-a.e.\ $\fz$ and every $k$. Hence,   for
$\sil$-a.e.\ $\fz$ and every $n\in \N$,
$$
\frac{\sum_{k=0}^{n-1} P_\ell^k \chi(\Phi(\fz))}{\sum_{k=0}^{n-1} P_\ell^k \chi(\fz)}=
\frac{\sum_{k=0}^{n-1} P_\ell^k (\chi \circ \Phi)(\fz)}{\sum_{k=0}^{n-1} P_\ell^k \chi(\fz)},
$$
so the ratio above is the Chacon--Ornstein ratio associated with the
functions $f_\chi:=\chi \circ \Phi$ and $g_\chi:=\chi$.
By assumption and the definition of $\cZ ^*$, both functions belong to
$L^1(\cZ ,\sil)$, and by assumption~(A3), for $\sir$-a.e.~$\chi$,
$\sum_{k=0}^{+\infty}{P_\ell^k g_\chi}>0$ $\sil$-a.e.
Applying the Chacon--Ornstein theorem to the totally conservative operator
$P_\ell$ yields for $\sir$-a.e.~$\chi$,
\begin{equation}\label{eq.ratio.2}\forall_\sil \fz,\,\lim_n
\frac{\sum_{k=0}^{n-1} P_\ell^k (\chi \circ \Phi)(\fz)}{\sum_{k=0}^{n-1} P_\ell^k \chi(\fz)}=R_\ell(\chi\circ \Phi\,\big/\,\chi)(\fz).\end{equation}
Comparing \eqref{eq.ratio.1} and \eqref{eq.ratio.2}, we obtain for
$\sigma_r$-a.e.~$\chi$ that\xout{,}
$$R_\ell(\chi\circ \Phi\,\big/\,\chi)(\fz)=1\qquad \textrm{$\sil$-a.e.}\xout{.}$$
Finally by \eqref{eq.mass.principle.total} we deduce that
for $\sir$-a.e.~$\chi$,
$$\int{\chi(\Phi(\fz))\,d\sil(\fz)}=\int{R_\ell(\chi\circ \Phi/\chi)(\fz)\,\chi(\fz)\,d\sil(\fz)}=\int{\chi(\fz)\, d\sil(\fz)}.$$
\end{proof}

\subsubsection{The group case}\label{subsec.group.case}
We now show how the abstract framework applies to random walks on a
locally compact group; this construction will serve as base model for the
affine setting of Section~\ref{sec.affine.application}.

Let $G$ be a locally compact second countable group, with left and right
Haar measures $m_\ell$ and $m_r$, acting on $\cZ:=G$ by left
multiplication, and let $\mu$ be a probability measure on $G$. The left Haar
measure $\sigma_\ell:=m_\ell$ is $\mu$-invariant.  When $G=\Aff(\R)$, the measure $m_\ell$ is the measure
$m$ of Section~\ref{sec.invariant.fibers}.

\emph{Companion space.} Fix $\xi\in L^1(G,m_\ell)\cap L^1(G,m_r)$
with $\xi>0$. We let $G$ act on observables on the right by composition i.e. 
$(\xi \cdot h)(\fz)=\xi(h\fz)$ for $h\in G$ and $\fz\in \cZ=G$ and we set $\cZ^*:=\{\xi_h := \xi \cdot h \mid h \in G\}$. 
 In other words, $\mathcal Z^*$ is the image of $G$ by the map $h\mapsto \xi_h$; note that $\cZ^*\subset L^1(\cZ,\sil)$ by left-invariance of $m_\ell=\sil$. Moreover,  $\cZ^*$ is stable under this right action as $\xi_h\cdot \gamma=\xi_{h\gamma}$ so that 
the right action 
of $G$ on $\cZ^*$ is induced by the right
multiplication in $G$.  From now, we endow $\cZ^*$ with the pushforward Borel $\sigma$-algebra on $G$ by the map $h\mapsto \xi_h$. Denote by $P_r$ the Markov kernel on $\cZ^*$ induced by the right random walk on $G$ i.e. $P_r F(\xi_h):=\E(F(\xi_h \cdot \gamma_1))$.

  Denote by  $\sigma_r$   the image of the right Haar measure $m_r$
under $h\mapsto\xi_h$. We now check Assumptions (A1)--(A3). The measure $\sir$ is $P_r$-invariant by right-invariance of
$m_r$. 
 This right invariance    implies also that     
$\|\ev_g\|_{L^1(\cZ^*,\sigma_r)}=m_r(\xi)<+\infty$ for every $g\in G$; showing (A2).  Since $\ev_g>0$ on $\cZ^*$ (since $\xi>0$ on $G$), and $\ev_g\in L^1(\cZ^*,\sir)$ as we just showed, 
this   implies that $\sir$ is $\sigma$-finite.
Now $L^1(G,m_r)$ is separable (since $G$ is second countable and locally compact) and the linear map $F\mapsto (h\mapsto F(\xi_h))$ is  an isometry from  $L^1(\cZ^*,\sir)$ to $L^1(G,m_r)$. 
Hence $L^1(\cZ^*,\sir)$ is separable. This finishes checking assumption (A1).
Finally (A3) holds since $\xi>0$.

Hence, in the group case,  both spaces $\cZ$ and $\cZ^*$ can heuristically  be identified with
the group $G$, the left Haar measure plays the role of $\sigma_\ell$, and the
right Haar measure the role of $\sigma_r$.

\emph{The intrinsic metric.} In this setting, \eqref{eq.dr} reads
\begin{equation}\label{eq.dr.group}
d_r(g,g')=\int_G|\xi(hg)-\xi(hg')|\,dm_r(h),
\qquad g,g'\in G.
\end{equation}
By right-invariance of $m_r$, $d_r$ is left-invariant, i.e.\
$d_r(hg,hg')=d_r(g,g')$, bounded by $2m_r(\xi)$, and
$d_r(g,g')=\|R_{g'^{-1}g}\xi-\xi\|_{L^1(m_r)}$, where $R_gf(h):=f(hg)$.
Since right translations are strongly continuous on $L^1(G,m_r)$
\cite[Prop.~2.42]{Folland}, $d_r$ is continuous for the topology of $G$:
\begin{equation}\label{eq.dr.group.continuity}
g_n\to g \ \text{in } G \quad\Longrightarrow\quad d_r(g_n,g)\to 0 .
\end{equation}
In particular, $d_r$-Lipschitz functions are continuous on $G$.

We now recover Guivarc'h--Keane--Roynette's theorem.

\begin{coro}[Guivarc'h--Keane--Roynette's theorem]\label{cor.GKR}
Let $G$ be a locally compact second countable non-unimodular group and let
$\mu$ be a probability measure on $G$ such that the closed subgroup $G_\mu$
generated by $\mathrm{supp}(\mu)$ equals $G$. Then the operator
$P_\mu f(g)=\int f(\gamma g)\,d\mu(\gamma)$ is totally dissipative on
$L^1(G,m_\ell)$.
\end{coro}

\begin{proof}
We apply Theorem~\ref{theorem.main.dissipativity} with $(\cZ,\sil)=(G,m_\ell)$ and $(\cZ^*,\sir)$ as defined above. 
For $g_0\in G$, the right multiplication $\Phi_{g_0}(g)=gg_0$ commutes with
the left action by associativity.
It remains to check the $\Phi_{g_0}$-Lipschitz minimality. Let $f$ be $d_r$-Lipschitz
with $f(\gamma g)=f(g)$ for $(\mu\otimes m_\ell)$-a.e.~$(\gamma,g)$. By
\eqref{eq.dr.group.continuity}, $f$ is continuous, which implies
$f(\gamma g)=f(g)$ for all $(\gamma,g)\in\mathrm{supp}(\mu)\times G$. Hence
$H:=\{\gamma\mid f(\gamma g)=f(g)\ \forall g\}$ is a closed subgroup
containing $\mathrm{supp}(\mu)$, so $H=G_\mu=G$ and $f$ is constant.

Finally, denoting by $\Delta$ the modular function of $G$, the Jacobian of
$\Phi_{g_0}$ with respect to $m_\ell$ equals $\Delta(g_0)^{-1}$, that
is, $\int\xi_h\circ\Phi_{g_0}\,dm_\ell=\Delta(g_0)^{-1}\int\xi_h\,dm_\ell$
for every $h\in G$; in particular $\xi_h\circ\Phi_{g_0}\in L^1(m_\ell)$.
If $P_\mu$ were not totally dissipative, Theorem~\ref{theorem.main.dissipativity}
would give $\int\xi_h\circ\Phi_{g_0}\,dm_\ell=\int\xi_h\,dm_\ell$
for $m_r$-a.e.\ $h$. Since $\int\xi_h\,dm_\ell>0$, this forces
$\Delta(g_0)=1$ for all $g_0\in G$, contradicting non-unimodularity.
\end{proof}

\subsection{Application to the affine diagonal action}
\label{sec.affine.application}
We now apply Theorem \ref{theorem.main.dissipativity} to show the
following. 
\begin{theo}\label{theorem.affine.dissipative}
Let $d\geq 2$.  Let $\mu$ be a probability measure on $\Aff(\R^d)$ satisfying assumptions
{\bf I}, {\bf P}, 
${\bf M}_1(\mu)$ and $\lambda_2<0$. Assume that
$\mu(\{\gamma\in \Aff(\R^d) \mid \gamma x=x\})<1$ for any $x\in \R^d$. Let
$\sigma$ be the $\mu$-invariant 
$\sigma$-finite measure on $\cZ_d$
of Proposition \ref{prop.invariant.sigma}. Then the Markov operator
$P_\mu f(\fz) = \int f(\gamma\fz)\,d\mu(\gamma)$ is totally dissipative on
$L^1(\cZ_d ,\sigma)$.
\end{theo}
The key observation is that the scaling map
\begin{equation}\label{eq.philambda}
\Phi_\lambda: \cZ_d \to \cZ_d, (x,v)\mapsto (x,\lambda v), \qquad \lambda\in\R^*,
\end{equation}
commutes with the diagonal action of $\AffRd$.
Moreover it preserves each fiber $\mathcal{Z}_1(\ell)\simeq \mathrm{Aff}(\R)$
and acts on it by right multiplication by $(0,\lambda)$. Since the
conditional measures of $\sigma$ are left Haar measures on fibers
(Proposition~\ref{prop.invariant.sigma}), it therefore  
 rescales $\sigma$ by the modular factor  
 $|\lambda|$ 
(see \eqref{eq.modular.jacobian}).
If $P_\mu$ were not totally dissipative, then
Theorem~\ref{theorem.main.dissipativity} would force $\Phi_\lambda$ to
preserve $\sigma$ when tested against suitable observables, yielding
$|\lambda|=1$, a
contradiction for $|\lambda| \neq 1$.

In this subsection, we verify the assumptions of
Theorem \ref{theorem.main.dissipativity}. We first construct the companion
space $(\cZ_d ^*,\sigma_r)$ and verify assumptions (A1)--(A3). Then, we
establish the $\Phi_\lambda$-Lipschitz minimality for the metric
$d_r$; this is the main geometric input.

In the sequel, we set $P_\ell:=P_\mu$  and
$\sil:=\sigma$, the $P_\mu$-invariant measure defined in ~\eqref{eq.def.sigma}. As $\sigma$ is a non-zero
Radon measure on the $\sigma$-compact space $\cZ_d$, it is
$\sigma$-finite, and we are in the setting of
Theorem \ref{theorem.main.dissipativity}.

\subsubsection{The companion space of $(\cZ_d, \sigma)$}
\label{subsec.companion}

\medskip

\noindent\textbf{The companion space $\cZ_d^*$.}

Set
$
\Theta := (\R^d)^* \times \R
$
and fix a measurable function $\xi:\R\times\R\to \R_+$ such that:
\begin{enumerate}
\item[i.] $\xi(b,a) > 0$ for $(b,a)\in\AffR$ and $\xi(b,0)=0$;
\item [ii.] $\xi$ is integrable on $\Aff(\R)$ with respect to both the left Haar measure $m=da\,db/|a|^2$ and the right Haar measure $m_r= da\,db/|a|$.
\end{enumerate}
For $\theta = (\tilde{v},b) \in \Theta$, define
\begin{equation}\label{eq.action.trial}
\xi_\theta:\cZ_d\to \R_+ \quad \fz=(x,v) \mapsto \xi_\theta(x,v):= \xi\bigl(\tilde{v}x+b,\,\tilde{v}v\bigr).
\end{equation}
Set $\cZ_d^* := \{\xi_\theta \mid  \theta\in\Theta\}$, endowed
with the $\sigma$-algebra   given by the
pushforward of the Borel structure of $\Theta$ by the map
$\theta \mapsto \xi_\theta$.

\medskip

\noindent\textbf{The right action of $\AffRd$ on $\cZ_d^*$.}

Recall that $\AffRd$ acts on the right on observables by
$(\chi\gamma)(\fz):=\chi(\gamma\fz)$.
We check that $\cZ_d^*$ is invariant for this action.
Indeed, for any
$\gamma=(B_\gamma,A_\gamma)\in\AffRd$, $\theta=(\tilde v, b)\in \Theta$
and $\fz=(x,v)\in\cZ_d$,
\begin{align*}
\xi_\theta(\gamma\fz)
= \xi\bigl(\tilde{v}(A_\gamma x+B_\gamma)+b,\,\tilde{v}(A_\gamma v)\bigr) 
=\xi\bigl((\tilde{v}A_\gamma) x + (\tilde{v}B_\gamma + b),\,
(\tilde{v}A_\gamma) v\bigr).
\end{align*}
Hence,
$\xi_\theta\gamma= \xi_{\theta\cdot \gamma}$, with
$\theta\cdot \gamma:=(\tilde v A_\gamma, \tilde v B_\gamma +b)\in \Theta$.
This defines a right action of $\AffRd$ on $\Theta$, and $\cZ_d^*$ is
stable under the right action.

In particular, the right $\mu$-random walk $(\gamma_1\cdots \gamma_n)$ on
$\Aff(\R^d)$ induces a Markov chain on $\cZ_d^*$ whose Markov operator is
given for every measurable function $F:\cZ_d^*\to \R_+$ and
$\theta\in \Theta$,
$$P_r F(\xi_\theta)=\int{F(\xi_\theta \gamma) d\mu(\gamma)}.$$


\medskip

\noindent\textbf{The measure $\sigma_r$ on $\cZ_d^*$.}
Recall (Subsection~\ref{projectiveaction}) that $\tilde\nu$ denotes the
unique right $\overline\mu$-invariant probability measure on
$\mathbf{P}((\R^d)^*)$. Moreover it is non-degenerate.
We lift this measure on $(\R^d)^*\setminus \{0\}$ by setting:
\[
\int_{(\R^d)^*\setminus\{0\}} f(\tilde{v})\,d\tilde{\frak{m}}(\tilde{v})
:= \int_{\mathbf{P}((\R^d)^*)} \int_{\R^*} f(a\tilde{v})\,\frac{da}{|a|}
\,d\tilde\nu([\tilde{v}]).
\]
The measure $\tilde{\frak{m}}$ is well-defined since
$da/|a|$ is the Haar measure on $\R^*$. Moreover one can
easily check that $\tilde{\frak{m}}$ is right $\overline{\mu}$-invariant.

We equip $\cZ_d^*$ with the pushforward measure of the
measure $d\tilde{\frak{m}}(\tilde{v})\,db$ on $\Theta$ by the map
$\theta \mapsto \xi_\theta$, that is
for any Borel function $F:\Theta\to \R_+$,
\begin{align}\label{eq.sigmar}
\int_{\cZ_d^*}{F d\sigma_r}
&:= \int_{(\R^d)^*}{\int_\R F(\xi_{(\tilde{v},b)})\, 
d\tilde{\frak{m}}(\tilde v) db}\notag \\
&=
\int_{\P((\R^d)^*)} \int_{\Aff(\R)} F(\xi_{(a\tilde{v},b)})\, 
d\tilde{\nu}([\tilde{v}])\,dm_r(b,a),
\end{align}
where $dm_r(b,a)=dbda/|a|$  is the right Haar measure on $\Aff(\R)$.
This measure is $P_r$-invariant. Indeed,
\begin{align*}
\int_{\cZ_d^*} &P_rF\,d\sigma_r
= \int_{\AffRd}\int_\Theta F(\xi_\theta \gamma)\,d\sigma_r(\xi_\theta)\,d\mu(\gamma)\\
&= \int_{\AffRd}\int_{(\R^d)^*\times\R} 
F\bigl(\xi_{(\tilde{v}A_\gamma,\,\tilde{v}B_\gamma+b)}\bigr)
\,d\tilde{\frak{m}}(\tilde{v})\,db\,d\mu(\gamma)\\
&= \int_{\AffRd}\int_{(\R^d)^*\times\R} 
F\bigl(\xi_{(\tilde{v}A_\gamma,\,b)}\bigr)
\,d\tilde{\frak{m}}(\tilde{v})\,db\,d\mu(\gamma)
\quad
\ \textrm{(invariance of $db$)}\\
&= \int_{\AffRd}\int_{(\R^d)^*\times\R} 
F\bigl(\xi_{(\tilde{v},\,b)}\bigr)
\,d\tilde{\frak{m}}(\tilde{v})\,db\,d\mu(\gamma) \quad 
\ \textrm{(right\ $\overline\mu$-stationarity of $\tilde{\frak{m}}$)}\\
&= \int_\Theta F(\xi_\theta)\,d\sigma_r(\xi_\theta),\quad \textrm{ so\  $\sigma_r$ is $P_r$-invariant.}
\end{align*}


\medskip

\noindent\textbf{Integrability of observables.}

We verify that $\xi_\theta \in L^1(\cZ_d,\sigma)$ for every 
$\theta=(\tilde{v},b)\in\Theta$.
Using the disintegration $\sigma = \int \lambda^\ell\,d\eta(\ell)$ 
(Proposition~\ref{prop.invariant.sigma}), it suffices to compute the 
integral of $\xi_\theta$ over each fiber $\cZ_1(\ell)$ for 
$\eta$-a.e.\ $\ell\in\AG$.

Fix $\ell = x_\ell +\R v_\ell  \in\AG$ with $\tilde{v}v_\ell \neq 0$ (which holds 
for $\nu$-a.e.\ direction $[v_\ell ]$ by non-degeneracy of $\nu$, hence for 
$\eta$-a.e.\ $\ell$). Under the identification 
$\cZ_1(\ell)\simeq\Aff(\R)$ via $(b',a')\mapsto (x_\ell +b'v_\ell ,\,a'v_\ell )$, 
and using $dm(b',a')=db'\,da'/|a'|^2$, we compute:
\begin{align*}
\int_{\cZ_1(\ell)} \xi_\theta(x,v)\,d\lambda^\ell(x,v)
&= \int_{\R\times\R^*} 
\xi\bigl(\tilde{v}(x_\ell +b'v_\ell )+b,\,\tilde{v}(a'v_\ell )\bigr)
\,dm(b',a')\\
&= \int_{\R\times\R^*} 
\xi\bigl((\tilde v x_\ell + b, \tilde v v_\ell) (b',a')\bigr)
\,dm(b',a').
\end{align*}
Since $m$ is left-invariant 
on $\Aff(\R)$, we obtain
\[
\int_{\cZ_1(\ell)} \xi_\theta(x,v)\,d\lambda^\ell(x,v)
= \|\xi\|_{L^1(\Aff(\R),m)}<+\infty.
\]

\medskip

\noindent\textbf{Verification of (A1)--(A3).}

\begin{prop}\label{prop.verification.A}
The companion space $\cZ_d^*$ with measure $\sigma_r$ satisfies 
assumptions \emph{{\rm (A1)--(A3)}} of Theorem~\ref{theorem.main.dissipativity}.
\end{prop}

\begin{proof}
\medskip
\noindent\emph{Assumption {\rm (A2}): $\ev_\fz\in L^1(\cZ_d^*,\sigma_r)$ for 
every $\fz = (x,v)\in\cZ_d$.}
 By \eqref{eq.sigmar},  
 right-invariance 
of $m_r$ and using  $\tilde{v}v\neq 0$ a.s.\ (by non-degeneracy of $\tilde\nu$),
\begin{align*}
\|\ev_\fz\|_{L^1(\cZ_d^*,\sigma_r)}
&= \int \vert \chi(\fz)\vert \,d\sigma_r(\chi)\\
&= \int \xi(a\tilde{v}x+b,\,a\tilde{v}v)\,dm_r(b,a)\,d\tilde{\nu}([\tilde{v}])\\
&= \iint \xi\bigl((b,a)(\tilde{v}x,\tilde{v}v)\bigr)
\,dm_r(b,a)\,d\tilde\nu([\tilde{v}])\\
&= \|\xi\|_{L^1(\Aff(\R),m_r)} < +\infty.   
\end{align*}
\medskip

 \emph{Assumption {\rm(A1):} $\sigma_r$ is $P_r$-invariant and 
$(\cZ_d^*,\sigma_r)$ is $L^1(\cZ^*,\sir)$ is separable.}
The $P_r$-invariance of $\sigma_r$ is established above. Fix $\fz=(x,v)\in \cZ$. 
We have that  $\ev_{\fz}>0$, $\sir$-almost everywhere. Indeed, since $\xi>0$ everywhere on $\Aff(\R)$, 
the zero set of $\ev_{\fz}$ is 
$\ev_{\fz}^{-1}(\{0\})=\{\xi_{(\tilde v, b)} \mid \tilde v v=0\}$,  which has zero $\sir$-measure by non-degeneracy of $\tilde \nu$. 
Since $\ev_z\in L^1(\cZ^*,\sir)$ as shown above,   we deduce that  
  $\sir$ is $\sigma$-finite. Finally, 
  $L^1(\cZ^*,\sir)$ 
  embeds isometrically into $L^1(\Theta,d\tilde{\frak{m}}(\tilde{v})\,db)$ via the map 
  $F\mapsto (\theta \mapsto F(\xi_\theta))$ and the space $L^1(\Theta,d\tilde{\frak{m}}(\tilde{v})\,db)$ is separable (as $\Theta$ is second countable and $d\tilde{\frak{m}}(\tilde{v})\,db$ is $\sigma$-finite). Hence $L^1(\cZ^*,\sir)$ is separable.

\emph{Assumption {\rm (A3)}: Weak irreducibility.}
Since $\xi > 0$ everywhere on $\Aff(\R)$, and since $\tilde{v}v\neq 0$ 
for $\tilde{\frak m}\otimes\nu$-a.e.\ $(\tilde{v},[v])$ by non-degeneracy of 
$\tilde\nu$ and $\nu$, we have $\xi_\theta(\fz) > 0$ for 
$(\sigma_r\otimes\sigma_\ell)$-a.e.\ $(\theta,\fz)$. Therefore,
\[
\sum_{n\geq 0} P_\ell^n \xi_\theta(\fz) \geq \xi_\theta(\fz) > 0
\qquad (\sigma_r\otimes\sigma_\ell)\text{-a.e.},
\]
which gives (A3). 
\end{proof}

\medskip
\noindent\textbf{Properties of the intrinsic metric $d_r$.}

The pseudo-metric $d_r$ on $\cZ_d$ defined in \eqref{eq.dr} takes the
explicit form: for $\fz=(x,v),\fz'=(x',v')\in\cZ_d$,
\begin{equation}\label{eq.dr.affine}
d_r(\fz,\fz')
= \int_{(\R^d)^*}\int_{\R}
\left|\xi\bigl(\tilde{v}x+b,\,\tilde{v}v\bigr)
-\xi\bigl(\tilde{v}x'+b,\,\tilde{v}v'\bigr)\right|
\,
d\tilde{\frak{m}}(\tilde v) db.
\end{equation}
Observe that, in dimension $d=1$, this gives the pseudo-metric $d_r^{\AffR}$ on
$\cZ_1\cong\AffR$
 defined in the group case (see
\eqref{eq.dr.group}) for $G=\Aff(\R)$, that is
for $g_1=(b_1,a_1),g_2=(b_2,a_2)\in\AffR$,
\begin{align*}
d_r^\AffR\bigl((b_1,a_1),(b_2,a_2)\bigr)
&:= \int_{\Aff(\R)}{|\xi(hg_1)-\xi(hg_2)|\, dm_r(h)}\\
&= \int_{\AffR} \left|\xi(a\,b_1+b, a\,a_1)
- \xi(a\,b_2+b,a\,a_2)\right| \frac{db\,da}{|a|}.
\end{align*}
Hence
$d_r^\AffR$ is left-invariant, bounded by $2\|\xi\|_{L^1(\AffR,m_r)}$ and,  by \eqref{eq.dr.group.continuity}, is 
continuous for the Euclidean topology of $\AffR$.

We then can deduce the following result:
\begin{lemma}[Properties of $d_r$]\label{lem-prop-dr}
The pseudo-metric $d_r$ on $\cZ_d$ can be written as
\begin{equation}\label{eq.dr.vs.d1}
d_r(\fz,\fz') = \int_{\mathbf{P}((\R^d)^*)}
d_r^\AffR\bigl((\tilde{v}x,\tilde{v}v),\,(\tilde{v}x',\tilde{v}v')\bigr)
\,d\tilde\nu([\tilde{v}]).
\end{equation}
Consequently, $d_r$ is bounded on $\cZ_d$, continuous for the Euclidean
topology on $\cZ_d\times \cZ_d$ and invariant under the action of
$\R^d\rtimes\R^*I < \AffRd$: for any $(B,\lambda I)\in\R^d\rtimes\R^*I$ and
$\fz,\fz'\in\cZ_d$,
\begin{equation}\label{eq.dr.aff.inv}
d_r\bigl((B,\lambda I)\fz,\,(B,\lambda I)\fz'\bigr)
= d_r(\fz,\fz').
\end{equation}
\end{lemma}

\begin{proof}
The representation \eqref{eq.dr.vs.d1} follows directly from the definition
of $d_r$ and $d_r^{\AffR}$ and the measure $\tilde{\frak{m}}$.  Note that,
for fixed $\fz,\fz'$, one has $\tilde vv\neq0$ and $\tilde vv'\neq0$ for
$\tilde\nu$-a.e.\ $[\tilde v]$ by non-degeneracy of $\tilde\nu$, and that the
integrand does not depend on the representative $\tilde v$ of $[\tilde v]$,
since replacing $\tilde v$ by $a\tilde v$ amounts to a left multiplication by
$(0,a)$ in $\AffR$.
Since $d_r^{\AffR}$ is bounded by $2\|\xi\|_{L^1(\AffR,m_r)}$ and
$\tilde{\nu}$ is a probability measure, it follows that $d_r$ is bounded on
$\cZ_d$ by $2\|\xi\|_{L^1(\AffR, m_r)}$.
We now show continuity of $d_r$. Let $\fz_n=(x_n,v_n)\to\fz_0=(x_0,v_0)$ in
$\cZ_d$ with respect to the Euclidean topology. For $\tilde\nu$-a.e.\
$[\tilde{v}]$, non-degeneracy of $\tilde\nu$ ensures $\tilde{v}v_0\neq 0$, so
that $(\tilde{v}x_n,\tilde{v}v_n)\to(\tilde{v}x_0,\tilde{v}v_0)$ in $\AffR$ for
the Euclidean topology. By the continuity property \eqref{eq.dr.group.continuity} of $d_r^\AffR$ , we have
\[
d_r^\AffR\bigl((\tilde{v}x_n,\tilde{v}v_n),\,
(\tilde{v}x_0,\tilde{v}v_0)\bigr)\to 0
\qquad \tilde\nu\text{-a.e.}
\]
Using \eqref{eq.dr.vs.d1} and the fact that $d_r^{\AffR}$ is bounded, the
dominated convergence theorem yields $d_r(\fz_n,\fz_0)\to 0$, which easily
implies the continuity of $d_r$ as a function on $\cZ_d\times \cZ_d$.
Finally, we show invariance of $d_r$ under the action of similitudes. For
$(B_0,\lambda I)\in\R^d\rtimes\R^*I$ and $\fz=(x,v)$, we have
$(B_0,\lambda I)\fz = (\lambda x+B_0,\lambda v)$,
$\tilde{v}(\lambda x+B_0) = \lambda\tilde{v}x+\tilde{v}B_0$ and
$\tilde{v}(\lambda v) = \lambda\tilde{v}v$. By \eqref{eq.dr.vs.d1},
\begin{align*}
d_r\bigl((B_0,\lambda I)\fz&,(B_0,\lambda I)\fz'\bigr)\\
&= \int_{\mathbf{P}((\R^d)^*)}
d_r^\AffR\bigl((\lambda\tilde{v}x+\tilde{v}B_0,\lambda\tilde{v}v),\,
(\lambda\tilde{v}x'+\tilde{v}B_0,\lambda\tilde{v}v')\bigr)
\,d\tilde\nu([\tilde{v}])\\
&= \int_{\mathbf{P}((\R^d)^*)}
d_r^\AffR\bigl((\tilde{v}B_0,\lambda)(\tilde{v}x,\tilde{v}v),\,(\tilde{v}B_0,\lambda)(\tilde{v}x',\tilde{v}v')\bigr)
\,d\tilde\nu([\tilde{v}])\\
&= \int_{\mathbf{P}((\R^d)^*)}
d_r^\AffR\bigl((\tilde{v}x,\tilde{v}v),\,(\tilde{v}x',\tilde{v}v')\bigr)
\,d\tilde\nu([\tilde{v}])\\
&= d_r(\fz,\fz'),
\end{align*}
where the third equality uses the left-invariance of $d_r^\AffR$.
\end{proof}

\subsubsection{$\Phi$-Lipschitz minimality of the diagonal action on $\cZ_d$}
Let $\gamma_\star=(B_\star,A_\star)\in \Gamma_\mu$ be such that $A_\star$ is proximal with dominant
eigenvalue $\lambda_\star$ satisfying $|\lambda_\star|\neq 1$.  For $d\geq2$,
such an element exists under assumptions {\bf I} and {\bf P}, as we show in
Lemma~\ref{lemma.proximal.modulus}.
We establish here the $\Phi_{\lambda_\star}$-Lipschitz minimality assumption of
Theorem~\ref{theorem.main.dissipativity}, where
$\Phi_{\lambda_\star}$ is the function introduced in \eqref{eq.philambda}.
\begin{prop}\label{prop.rigid.measure}
Let $d_r$ be the pseudo metric defined in \eqref{eq.dr.affine}. Under
the assumptions of Theorem~\ref{theorem.affine.dissipative}, for any
$d_r$-Lipschitz function $f:\cZ_d\to\R$ such that
\[
f(\gamma\fz) = f(\fz)
\qquad \text{for $(\mu\otimes\sigma)$-a.e.\ $(\gamma,\fz)$,}
\]
and for every $\gamma_\star=(B_\star,A_\star)\in \Gamma_\mu$ such that $A_\star$ is proximal with
dominant eigenvalue $|\lambda_\star|\neq 1$, it holds that
\[
f(x,v) = f\circ \Phi_{\lambda_\star}(x,v)=f(x,\lambda_\star v)
\qquad \textrm{for $\sigma$-a.e. $\fz=(x,v)$.}
\]
In particular, the action of $(\Aff(\R^d), \mu)$ is
$\Phi_{\lambda_\star}$-Lipschitz minimal in the sense of
\eqref{def.lipshitz}.
\end{prop}
In  what follows, we endow $\cZ_d$ with the Euclidean topology it inherits as a subset of  $\R^d\times \R^d$.
First, observe that the almost everywhere invariance assumption on $f$ can be
promoted to pointwise invariance on a closed subset of $\cZ_d$.
Indeed, consider the probability measure
$\hat \mu:= \sum_{n\geq 1}2^{-n}\mu^{*n}$ on $\Aff(\R^d)$ whose support in
$\Aff(\R^d)$ is equal to the closed semigroup $\Gamma_\mu$ generated by
$\textrm{supp}(\mu)$ in $\Aff(\R^d)$. Since $\sigma$ is $\mu$-invariant, it is
also $\hat\mu$-invariant. It is then easy to verify that
$f(\gamma \fz)=f(\fz)$ for $\hat \mu\otimes \sigma$-a.e.~$(\gamma,\fz)$. On the
other hand, since $d_r$ is continuous (Lemma~\ref{lem-prop-dr}), the function
$f$ is continuous on $\cZ_d$.
Denote by $\supp(\sigma)$ the support of $\sigma$ in $\cZ_d$, i.e.\ the
smallest closed subset of $\cZ_d$ of full $\sigma$-measure. By
$\mu$-stationarity of $\sigma$ and continuity of the action, $\supp(\sigma)$
is $\Gamma_\mu$-invariant, and the support of $\hat\mu\otimes\sigma$ in
$\Aff(\R^d)\times \cZ_d$ is $\Gamma_\mu\times\supp(\sigma)$.
By continuity of the action, the set of $(\gamma,\fz)$ such that
$f(\gamma\fz)=f(\fz)$ is closed and has full $\hat\mu\otimes\sigma$-measure,
so it follows that $f(\gamma \fz)=f(\fz)$ for every
$(\gamma,\fz)\in \Gamma_\mu\times\supp(\sigma)$. Hence
\begin{equation}\label{eq.invariance.pointwise}
 \forall (\gamma, \fz)\in \Gamma_\mu\times \supp(\sigma), \quad f(\gamma \fz)=f(\fz).
\end{equation}
Similarly, since $f$ and $f\circ\Phi_{\lambda_\star}$ are continuous, the
conclusion of the proposition is equivalent to the equality
$f=f\circ\Phi_{\lambda_\star}$ on $\supp(\sigma)$. 

\medskip

\medskip 

\medskip

\noindent {\bf 
One-dimensional minimality}

In dimension $d=1$, the $\Phi$-Lipschitz minimality follows from
the following elementary lemma, which is the one-dimensional version of
Proposition~\ref{prop.rigid.measure}  (related arguments appear for example  in \cite[IV, Proposition 2]{Elie79} and \cite[Proposition 2.5.4]{BDM}).

\begin{lemma}\label{lemma.rigidity.1D}
Let $\tau_\star,\tau\in\Aff(\R)$, with $\tau_\star(x)=\alpha_\star x+\beta_\star$
and $\tau(x)=\alpha x+\beta$.
Suppose that they do not fix the same point and that $|\alpha_\star|\neq 1$.
Let $f:\cZ_1=\R\times\R^*\to\R$ be continuous and invariant under the diagonal
action of $\tau_\star$ and $\tau$. Then $f$ is independent of the
$b$-variable and $f(b,a)=f(b,\alpha_\star a)$.
\end{lemma}

\begin{proof}
First, observe that $f$ is invariant under the group generated by $\tau_\star$
and $\tau$. For $\beta'\in\R$, denote by $T_{\beta'}(x):=x+\beta'$ the
translation by $\beta'$.

\emph{Step 1: Find a translation.}
The commutator $\tau\tau_\star\tau^{-1}\tau_\star^{-1}$ has linear part
$1$, hence it is a translation $T_{\beta_0}$. We claim that $\beta_0 \neq 0$.
Indeed, let $c_\star = \beta_\star/(1-\alpha_\star)$ be the unique fixed point
of $\tau_\star$. If $\beta_0=0$, then
$$c_\star=T_{\beta_0}(c_\star)=\tau\tau_\star\tau^{-1}\tau_\star^{-1}(c_\star)
=\tau\tau_\star\tau^{-1}(c_\star).$$
Thus $\tau^{-1}(c_\star)$ is a fixed point of $\tau_\star$ and, by uniqueness,
$\tau^{-1}(c_\star)=c_\star$. This contradicts the assumption that $\tau_\star$
and $\tau$ do not fix the same point.

\emph{Step 2: Generate small translations.}
For each $k\in\Z$, the conjugate $\tau_\star^k T_{\beta_0} \tau_\star^{-k}$ is
the translation $T_{\alpha_\star^k \beta_0}$. Since $|\alpha_\star|\neq 1$,
there exists a sequence $(k_n)$ in $\Z$ such that
$\alpha_\star^{k_n}\to 0$.
Since $f$ is invariant under these translations, we have
$f(b + \alpha_\star^{k_n}\beta_0, a) = f(b,a)$ for every
$(b,a)\in \R\times\R^*$ and every $n$.

\emph{Step 3: Conclusion.}
The set
\[
E := \{t\in\R \mid f(b+t,a)=f(b,a)\ \forall\,(b,a)\in \R\times \R^*\}
\]
is a closed additive subgroup of $\R$ containing the sequence
$\alpha_\star^{k_n}\beta_0 \to 0$ with $\alpha_\star^{k_n}\beta_0 \neq 0$ for
all $n$. Hence $E = \R$, and $f$ is independent of the $b$-variable.
Therefore, invariance under $\tau_\star$ yields
$f(b,a) = f(\alpha_\star b + \beta_\star, \alpha_\star a) = f(b, \alpha_\star a)$.
\end{proof}
\medskip

\noindent\textbf{Reduction to dimension $d=1$.}
To extend this result to higher dimensions, one cannot directly argue
fiberwise on each $\mathcal{Z}_1(\ell)$ since $\Gamma_\mu$ may not leave
invariant a line of $\R^d$. We therefore make use of the fact that the
dynamics on $\R^d$ induced by a proximal element $A_\star$ collapse onto its
dominant direction. More precisely, recall that a $d\times d$-matrix    is said to be 
proximal if there exists a dominant eigenvalue $\lambda_{A}$, which is necessarily real. Denoting by $\R v_A$ the associated eigenline, there is a rank-one projection map $\pi_A:\R^d\to \R v_A$ such that $\ker \pi_A$ is $A$-invariant and  
$$A^n=\lambda_{A}^n(\pi_{A}+ R_n), \quad \vert R_n\vert\to 0.$$
The projection satisfies $\pi_A A=\lambda_A\pi_A$, and hence, for the
affine action of $\gamma=(B,A)\in \AffRd$,
\begin{equation}\label{eq.intertwine.proximal}
\pi_A(\gamma x)=\pi_A(Ax +B)=\lambda_A\pi_A(x)+\pi_A(B).
\end{equation}
For $\fz=(x,v)\in \cZ_d$ with $v\notin \ker\pi_A$, we write
$\pi_A(\fz) := (\pi_A(x),\pi_A(v)) \in \R v_A \times (\R v_A \setminus\{0\})$.

\begin{lemma}\label{lemma.proximal.collapse}
Let $\gamma_\star=(B_\star,A_\star)\in \Aff(\R^d)$ with $A_\star$ proximal with associated projection $\pi_\star$ .
For every $\fz=(x,v),\fz'=(x',v')\in\cZ_d$ with
$v,v'\notin \ker\pi_\star$,
\[
d_r(\gamma_\star^n\fz,\,\gamma_\star^n\fz')
\underset{n\to+\infty}{\longrightarrow}
d_r(\pi_\star(\fz),\,\pi_\star(\fz')).
\]
In particular, if in addition $\gamma_\star\in \Gamma_\mu$ and $f:\cZ_d\to\R$ is
$d_r$-Lipschitz and satisfies \eqref{eq.invariance.pointwise}, then for
every $\fz,\fz'\in \supp(\sigma)$ with
$v,v'\notin \ker\pi_\star$,
\begin{equation}\label{eq.projection.rigidity}
 \pi_\star(\fz) = \pi_\star(\fz') \implies f(\fz) = f(\fz').
\end{equation}
\end{lemma}
\begin{proof}
By proximality, $A_\star^n = \lambda_\star^n(\pi_\star + R_n)$ with $\vert R_n\vert \to 0$.
Hence, by \eqref{eq.intertwine.proximal}, setting $\gamma_\star^n=(B_{\gamma_\star^n}, A_\star^n)$,
\[
\gamma_\star^n\fz
= (B_{\gamma_\star^n},\lambda_\star^n I)\bigl(\pi_\star(\fz)+\varepsilon_n(\fz)\bigr) \ \textrm{with }\
\varepsilon_n(\fz):=(R_n x, R_n v)\to 0 \  \text{in } \R^d\times\R^d.
\]
By Lemma~\ref{lem-prop-dr}, $d_r$ is invariant under $\R^d\rtimes\R^*I$ and
continuous on $\cZ_d\times \cZ_d$. Hence, for any $\fz, \fz'$ in $\cZ_d$ with
$v,v'\not\in \ker\pi_\star$,
\[
d_r(\gamma_\star^n\fz,\,\gamma_\star^n\fz')
= d_r\bigl(\pi_\star(\fz)+\varepsilon_n(\fz),\,\pi_\star(\fz')+\varepsilon_n(\fz')\bigr)
\underset{n\to\infty}{\longrightarrow} d_r(\pi_\star(\fz),\pi_\star(\fz')).
\]
Hence if $\gamma_\star\in \Gamma_\mu$ and $\fz,\fz'\in \supp(\sigma)$ with
$v,v'\not\in \ker\pi_\star$  and $\pi_\star(\fz)=\pi_\star(\fz')$,
\eqref{eq.invariance.pointwise} yields
$$|f(\fz)-f(\fz')| = |f(\gamma_\star^n\fz)-f(\gamma_\star^n\fz')|
\leq \mathrm{Lip}(f)\,d_r(\gamma_\star^n\fz,\gamma_\star^n\fz') \underset{n
\to +\infty}{\longrightarrow} 0.$$
\end{proof}

In the following lemma, we fix $\gamma_\star=(B_\star, A_\star)\in \Gamma_\mu$
such that $A_\star$ is proximal. For short we denote
$\lambda_{A_\star}=\lambda_\star$, $\pi_{A_\star}=\pi_\star$ and
$v_{A_\star}=v_\star$. We also denote by
$\mathcal{L}_{\star}\subset \AG$ the set of lines $\ell$ such that
$\pi_\star(v_\ell)\neq 0$ and $\mathcal{Z}_1(\ell)\subset \supp(\sigma)$.
\begin{lemma}[Reduction to 1D]\label{lemma.reduction.1D}
Let $f:\cZ_d\to\R$ be a $d_r$-Lipschitz function
satisfying \eqref{eq.invariance.pointwise} and let 
$\ell=x_\ell+\R v_\ell\in \mathcal{L}_\star$. Define
$f_\ell:\cZ_1=\R\times\R^*\to\R$ by
\begin{equation}\label{def:gell}
f_\ell(b,a) := f(x_{\ell}+b\,v_{\ell},\,a\,v_{\ell})=f(\varphi_\ell(b,a)).
\end{equation}
Then, for every $\gamma=(B_\gamma,A_\gamma)\in\Gamma_\mu$ with
$\pi_\star(A_\gamma v_{\ell})\neq 0$,
\[
f_\ell(b,a)=f_\ell(\beta_\gamma+\alpha_\gamma b,\,\alpha_\gamma a)= f_\ell(\tau_\gamma\cdot(b,a)),
\]
where $\tau_\gamma=(\beta_\gamma,\alpha_\gamma)\in\AffR$ is the unique
one-dimensional affine transformation such that
$\pi_\star(\gamma \varphi_\ell(b,a))=\pi_\star(\varphi_\ell(\tau_\gamma (b,a)))$,
that is, $\alpha_\gamma,\beta_\gamma\in \R$ are defined by
$$\alpha_\gamma \pi_\star (v_\ell) = \pi_\star(A_\gamma v_{\ell})\quad\mbox{and}\quad
\beta_\gamma \pi_\star(v_\ell) = \pi_\star(\gamma x_{\ell}) - \pi_\star(x_{\ell}).$$
\end{lemma}
Let us briefly explain in words the geometric content of this lemma.
 Let $\ell$ be a line transverse to
 $\ker\pi_\star$. The lemma shows that, although $\Gamma_\mu$ does not
 stabilize $\ell$, the restriction $f_{|_{\mathcal{Z}_1(\ell)}}$ of $f$ to
 $\cZ_1(\ell)$ is still invariant under some subgroup of affinities of
 $\ell$. This group is generated by maps
 $\tilde{\tau}_\gamma:=\varphi_\ell\tau_\gamma \varphi_\ell^{-1}:\ell\to \ell$
 such that $\gamma_{|_{\ell}}=\tilde{\tau}_\gamma$ modulo
 $\ker\pi_\star$, i.e. $\pi_\star(\gamma x)=\pi_\star(\tilde{\tau}_\gamma x)$
 for every $\gamma \in \Gamma_\mu$ such that the line $\gamma \cdot \ell$ is
 still transverse to $\ker\pi_\star$.

\begin{proof}
Let $\fz=\varphi_\ell(b,a)=(x_{\ell}+b\,v_{\ell},\,a\,v_{\ell})\in \cZ_1(\ell)\subset \supp(\sigma)$.
We have
$\gamma \fz=(\gamma x_{\ell}+b A_\gamma v_{\ell},\,a\,A_\gamma v_{\ell})$.
Applying $\pi_\star$ and using the definitions of $\alpha_\gamma$ and
$\beta_\gamma$, we obtain:
\begin{align*}
   \pi_\star(\gamma  \fz)&=\pi_\star(\gamma  x_{\ell}+b\,A_\gamma v_{\ell},\,a\,A_\gamma v_{\ell})\\
   &=(\pi_\star(\gamma  x_{\ell})+ b \pi_\star(A_\gamma v_{\ell}),\,a\,\pi_\star(A_\gamma v_{\ell}))\\
&=\bigl(\pi_\star(x_{\ell})+(\beta_\gamma+\alpha_\gamma b)\pi_\star(v_\ell),\,
\alpha_\gamma a\,\pi_\star(v_\ell)\bigr)\\
&=\pi_\star\bigl(x_{\ell}+(\beta_\gamma+\alpha_\gamma b)v_{\ell},\,\alpha_\gamma a\,v_{\ell}\bigr)
=\pi_\star(\varphi_\ell(\tau_\gamma (b,a))).
\end{align*}
Since $v_{\ell}\notin \ker\pi_\star$ and
$A_\gamma v_{\ell}\notin \ker\pi_\star$, and since $\gamma\fz$ and
$\varphi_\ell(\tau_\gamma(b,a))$ belong to $\supp(\sigma)$ (by
$\Gamma_\mu$-invariance of $\supp(\sigma)$ and $\ell\in\mathcal L_\star$),
Lemma~\ref{lemma.proximal.collapse} applied to $\gamma_\star$ yields
$f(\gamma  \fz) =f\bigl(x_{\ell}+(\beta_\gamma+\alpha_\gamma b)v_{\ell},\,\alpha_\gamma a\,v_{\ell}\bigr)$.
But by \eqref{eq.invariance.pointwise}, $f(\fz)=f(\gamma \fz)$.
Hence $f_\ell(b,a)=f(\fz)=f(\gamma  \fz)=f_\ell(\beta_\gamma+\alpha_\gamma b,\alpha_\gamma a)$,
as desired.
\end{proof}

\medskip

\noindent\textbf{Proof of Proposition~\ref{prop.rigid.measure}.}
The strategy is to apply Lemma~\ref{lemma.rigidity.1D} to
$f_\ell:=f\circ \phi_\ell$, with $\tau_\star=\tau_{\gamma_\star}$ and $\tau=\tau_\gamma$ for a
suitable $\gamma\in\Gamma_\mu$, for $\eta$-a.e.\ $\ell$.
Fix a $d_r$-Lipschitz function $f$ satisfying
\eqref{eq.invariance.pointwise}, and an element
$\gamma_\star=(B_\star,A_\star)\in\Gamma_\mu$ with $A_\star$ proximal and
$|\lambda_\star|\neq 1$. Recall that $\mathcal L_\star\subset\AG$ is the set
of lines $\ell=x_\ell+\R v_\ell$ such that $\pi_\star(v_\ell)\neq0$ and
$\cZ_1(\ell)\subset\supp(\sigma)$.

\medskip

\textit{Step 1: $\eta(\mathcal{L}_\star)=1$.}
Since $\sigma$  disintegrates over $\eta$ with Haar conditional measures
$\lambda^\ell$ (Proposition~\ref{prop.invariant.sigma}) then $\supp(\lambda^\ell) \subseteq \supp(\sigma)$ for $\eta$-a.e.\ $\ell\in\AG$. But $\lambda^\ell$
has full support in $\cZ_1(\ell)$ (i.e.~$\cZ_1(\ell) = \supp (\lambda^\ell)$). Hence   $\cZ_1(\ell)\subset \supp(\sigma)$
for $\eta$-a.e.\ $\ell\in\AG$.
Moreover, since the $\overline{\mu}$-invariant probability measure $\nu$ on $\PRd$ is non-degenerate
(by strong irreducibility of $\Gamma_{\overline\mu}$), we have
$\pi_\star(v_\ell)\neq 0$ for $\nu$-a.e.\ $[v_\ell]\in\PRd$, hence for
$\eta$-a.e.\ $\ell\in\AG$. Hence $\eta(\mathcal{L}_\star)=1$.

\textit{Step 2. 
For $\eta$-a.e.\ $\ell\in\mathcal L_\star$, there exists
$\gamma\in\Gamma_\mu$ such that $\pi_\star (A_\gamma v_\ell)\neq 0$ (so that  $\tau_\gamma$ is defined),  and such that $\tau_\gamma$ and
$\tau_{\gamma_\star}$ have no common fixed point.}

Let $\ell\in\mathcal L_\star$. Since $\pi_\star\circ A_\star=\lambda_\star\pi_\star$,
$\tau_{\gamma_\star}$ is defined and $\alpha_{\gamma_\star}=\lambda_\star\neq 1$,
so that $\tau_{\gamma_\star}$ has a unique fixed point
$c_{\star,\ell}:=\frac{\beta_{\gamma_\star}}{1-\lambda_\star}$; equivalently,
$x_{\star,\ell}:=\varphi_\ell(c_{\star,\ell})=x_\ell+c_{\star,\ell}v_\ell$ is the unique point in $\ell$
such that $\pi_\star(\gamma_\star\,x_{\star,\ell})=\pi_\star(x_{\star,\ell})$.
By definition of $\tau_\gamma$, for $\gamma\in\Gamma_\mu$ with
$\pi_\star(A_\gamma v_\ell)\neq0$, we have $\tau_\gamma(c_{\star,\ell})\neq
c_{\star,\ell}$ if and only if
$\pi_\star(\gamma x_{\star,\ell})\neq\pi_\star(x_{\star,\ell})$.
Hence, it is enough to show  that for $\eta$-a.e.~$\ell$,
\begin{equation}\label{eq.mumu}
\hat{\mu}(\{ \gamma \mid \pi_\star(A_\gamma v_\ell)\neq 0,\
\pi_\star(\gamma\,x_{\star,\ell})\neq \pi_\star(x_{\star,\ell})\})>0,
\end{equation}
where $\hat \mu:= \sum_{n\geq 1}2^{-n}\mu^{*n}$.

Indeed, on the one hand, since $\nu$ is non-degenerate on $\PRd$ and
$\hat{\mu}$-invariant, for $\hat \mu\otimes \eta$-a.e.~$(\gamma,\ell)$, it
holds that $A_\gamma v_\ell \not\in \ker\pi_\star$.
On the other hand, we claim that for $\eta$-a.e.~$\ell$, the set of
$\gamma\in \Gamma_\mu$ such that
$\gamma\,x_{\star,\ell}\notin x_{\star,\ell}+\ker\pi_\star$ has positive
$\hat{\mu}$-measure. Suppose this is not the case. Then there exists a set
$\mathcal{L}\subset \AG$ such that $\eta(\mathcal{L})>0$ and, for every
$\ell\in \mathcal{L}$,
$$
\hat{\mu}(\underbrace{
\{\gamma \mid  \pi_\star(\gamma x_{\star,\ell})=\pi_\star(x_{\star,\ell})\}}_{E_\ell})=1.
$$
Since $E_{\ell}$ is closed in $\Aff(\R^d)$, one has that
$\Gamma_\mu=\supp \hat{\mu}\subset E_\ell$ for every $\ell \in \mathcal L$. In other
terms, $\gamma x_{\star,\ell}\in x_{\star,\ell}+\ker\pi_\star$ for every
$\gamma\in \Gamma_\mu$ and $\ell \in \mathcal L$.
Hence the linear span $W_\ell$ of
$\{\gamma x_{\star,\ell}-x_{\star,\ell}\mid\gamma\in\Gamma_\mu\}$ is included in
$\ker\pi_\star$, hence it is a proper subspace. Moreover it is
$\Gamma_{\overline\mu}$-invariant, since
$A_{\gamma'}(\gamma x_{\star,\ell}-x_{\star,\ell})
=(\gamma'\gamma x_{\star,\ell}-x_{\star,\ell})-(\gamma'x_{\star,\ell}-x_{\star,\ell})$
for $\gamma,\gamma'\in\Gamma_\mu$.
By irreducibility of $\Gamma_{\overline\mu}$, it holds that
$W_\ell=\{0\}$. Hence $\Gamma_\mu x_{\star,\ell}=\{x_{\star,\ell}\}$,
meaning $x_{\star,\ell}$ is a common fixed point of $\Gamma_\mu$,
contradicting our assumption. Tonelli's theorem then yields
\eqref{eq.mumu}. By the equivalence above,
this ensures that for $\eta$-a.e.\ $\ell$, there exists
$\gamma\in \Gamma_\mu$ (depending on $\ell$) such that $\tau_\gamma$ is
defined and $\tau_\gamma (c_{\star,\ell})\neq c_{\star,\ell}$.

\medskip

\textit{Step 3: Conclusion.}
Let $\ell\in\mathcal L_\star$ and $\gamma$ be as in Step 2.
The function $f_\ell:\cZ_1\to \R$ is continuous, since the $d_r$-Lipschitz
function $f$ is continuous by Lemma \ref{lem-prop-dr}. By
Lemma \ref{lemma.reduction.1D}, $f_\ell$ is invariant under both
$\tau_{\gamma_\star}$ and $\tau_\gamma$. Since $\tau_{\gamma_\star}$ and
$\tau_\gamma$ do not fix the same point and the linear part of
$\tau_{\gamma_\star}$ is $\lambda_\star$ with $|\lambda_\star|\neq 1$,
Lemma~\ref{lemma.rigidity.1D} applies with $\tau_\star=\tau_{\gamma_\star}$
and $\tau=\tau_\gamma$ and yields
$$f_\ell(b,a) = f_\ell(b,\lambda_\star a) \qquad \forall\,(b,a)\in\R\times\R^*.$$
Pushing forward by $\varphi_\ell$, for every $(x,v)\in\cZ_1(\ell)$ written as
$(x,v)=(x_\ell+bv_\ell,\,av_\ell)$, we obtain
$$f(x,v) = f_\ell(b,a) = f_\ell(b,\lambda_\star a) = f(x_\ell+bv_\ell,\,\lambda_\star av_\ell)
= f(x,\lambda_\star v).$$
Since this holds for $\eta$-a.e.~$\ell$ and $\sigma$ disintegrates over
$\eta$, we conclude that $f(x,v)=f(x,\lambda_\star v)$ holds $\sigma$-a.e., as
claimed.

\rightline{$\Box$}
\subsubsection{Proof of Theorem~\ref{theorem.affine.dissipative}}
\label{subsec.proof.dissipative}

Before concluding the proof of  Theorem
\ref{theorem.affine.dissipative}, we record a standard algebraic consequence
of assumptions {\bf I} and {\bf P}. It can be deduced from general
non-arithmeticity results for Zariski dense semigroups of reductive groups
\cite{benoist.lineaire.2, GR-large}, and \cite[Proposition 8.8]{BQ-book}.
Since in our setting a short argument is available, we include it for
completeness.

\begin{lemma}\label{lemma.proximal.modulus}
Let $\Gamma$ be a strongly irreducible and proximal semigroup of $\GL_d(\R)$,
$d\geq 2$. Then there exists a proximal element $A_\star\in \Gamma$ with
$|\lambda_{A_\star}|\neq 1$.
\end{lemma}

\begin{proof}
Let $A_\star\in \Gamma$ be proximal, with dominant eigenvalue
$\lambda_\star$, dominant eigenvector $v_\star$ and repelling hyperplane
$\ker(\tilde v_\star)$, where we normalize $\tilde v_\star$ so that
$\tilde v_\star v_\star=1$. If $|\lambda_\star|\neq 1$, we are done.
Otherwise, since $\lambda_\star$ is real, $\lambda_\star\in \{-1,1\}$, and
replacing if necessary $A_\star$ by $A_\star^2$ we may assume $\lambda_\star=1$.
Then $A_\star^n$ converges to the rank-one projection
\begin{equation}\label{eq.conv.s}
A_\star^n \underset{n\to +\infty}{\longrightarrow} \pi_\star,
\qquad \pi_{\star} x =( \tilde v_\star x)v_\star.
\end{equation}

\emph{Step 1: It suffices to find $A\in\Gamma$ such that 
$\tilde v_\star A v_\star\notin\{-1,0,1\}$.}
Indeed, for such $A$, by \eqref{eq.conv.s},
$A_\star^nA\to\pi_\star A$ as $n\to+\infty$. The endomorphism $\pi_\star A$ has
rank one, with image $\R v_\star$, and its spectrum consists of
$\tilde v_\star Av_\star\neq0$ and $0$ (with multiplicity $d-1$); hence it is
proximal. Since proximality is an open condition in $\operatorname{End}(\R^d)$
and the simple dominant eigenvalue depends continuously on the endomorphism,
$A_\star^nA\in\Gamma$ is proximal for $n$ large enough, with dominant
eigenvalue $\lambda_n\to\tilde v_\star Av_\star\notin\{-1,1\}$. Thus
$|\lambda_n|\neq1$ for $n$ large enough.

\emph{Step 2: The set $\tilde v_\star\Gamma v_\star\subset\R$ is infinite.}
This is clearly stronger than the required condition in Step 1. First, the orbit $\Gamma v_\star$ is infinite as  otherwise, the
finite union of lines $\bigcup_{x\in \Gamma v_\star}\R x$ would be
$\Gamma$-invariant, and proper since $d\geq2$, contradicting strong
irreducibility. Moreover, since $\Gamma$ acts irreducibly on the left on $\R^d$,
it acts irreducibly on the right on the dual space $(\R^d)^*$. Hence  
$\operatorname{Span}(\tilde v_\star A \mid A\in \Gamma)=(\R^d)^*$. We can 
therefore choose $A_1,\dots, A_d\in \Gamma$ such that
$\tilde v_\star A_1,\dots, \tilde v_\star A_d$ form a basis of 
$(\R^d)^*$. This implies that the linear map
$T:\R^d\to \R^d$, $x\mapsto (\tilde v_\star A_1 x, \dots, \tilde v_\star A_d x)$
is injective. But, as  $A_i\Gamma\subset\Gamma$, we have
$T(\Gamma v_\star)\subset(\tilde v_\star\Gamma v_\star)^d$. Hence, since $\Gamma v_\star$
is infinite and $T$ is injective, we deduce that $(\tilde v_\star\Gamma v_\star)^d$ is
infinite, hence so is $\tilde v_\star\Gamma v_\star$.
\end{proof}

\begin{proof}[of Theorem~\ref{theorem.affine.dissipative}]
We apply Theorem~\ref{theorem.main.dissipativity} with the companion space
$(\cZ_d^*,\sir)$ constructed in Section~\ref{subsec.companion}. Assumptions
(A1)--(A3) were verified there.
By Lemma~\ref{lemma.proximal.modulus}, choose
$\gamma_\star=(B_\star,A_\star)\in\Gamma_\mu$ such that $A_\star$ is proximal
with dominant eigenvalue $\lambda_\star$ satisfying $|\lambda_{\star}|\neq 1$.
We take
\[
\Phi=\Phi_{\lambda_\star},\qquad \Phi_{\lambda_\star}(x,v)=(x,\lambda_\star v).
\]
This map commutes with the diagonal affine action. Moreover, by
Proposition~\ref{prop.rigid.measure}, the action is
$\Phi_{\lambda_\star}$-Lipschitz minimal.

The Jacobian of $\Phi_\lambda$ with respect to $\sigma$ is the constant
$|\lambda|$: for every $\lambda\in\R^*$ and every non-negative measurable
function $F:\cZ_d\to\R_+$,
\begin{equation}\label{eq.modular.jacobian}
\int_{\cZ_d} F(\Phi_\lambda (\fz))\,d\sigma(\fz)
=|\lambda|\int_{\cZ_d} F(\fz)\,d\sigma(\fz).
\end{equation}
Indeed, we use the disintegration $\sigma=\int\lambda^\ell\,d\eta(\ell)$, where
$\lambda^\ell$ is the pushforward by $\varphi_\ell$ of the left Haar measure
$\frac{db\,da}{|a|^2}$ on $\AffR$ (see \eqref{phiell}). For each
$\ell=x_\ell+\R v_\ell$, the change of variables $a'=\lambda a$ gives
\[
\int_{\cZ_1(\ell)} F(\Phi_\lambda (\fz))\,d\lambda^\ell(\fz)
= \int_{\R\times\R^*} F(x_\ell +bv_\ell,\lambda a v_\ell)\,\frac{db\,da}{|a|^2}
= |\lambda|\int_{\cZ_1(\ell)} F(\fz)\,d\lambda^\ell(\fz),
\]
and \eqref{eq.modular.jacobian} follows by integrating over $\eta$. 

It remains only to check the integrability assumption of
Theorem~\ref{theorem.main.dissipativity}. Indeed, by
\eqref{eq.modular.jacobian}, for every $\xi_\theta\in \cZ_d^*$,
\[
\int_{\cZ_d} \xi_\theta(\Phi_{\lambda_\star}(\fz))\,d\sigma(\fz)
=
|\lambda_\star|\int_{\cZ_d} \xi_\theta(\fz)\,d\sigma(\fz)
<+\infty.
\]
Assume by contradiction that $P_\mu$ is not totally dissipative. Then
Theorem~\ref{theorem.main.dissipativity} gives, for $\sigma_r$-a.e.
$\xi_\theta$,
\[
\int_{\cZ_d} \xi_\theta(\Phi_{\lambda_\star}(\fz))\,d\sigma(\fz)
=
\int_{\cZ_d} \xi_\theta(\fz)\,d\sigma(\fz).
\]
Using again \eqref{eq.modular.jacobian}, we obtain
\[
|\lambda_\star|\int_{\cZ_d} \xi_\theta(\fz)\,d\sigma(\fz)
=
\int_{\cZ_d} \xi_\theta(\fz)\,d\sigma(\fz)
\]
with
$
\int_{\cZ_d} \xi_\theta(\fz)\,d\sigma(\fz)
=
\|\xi\|_{L^1(\Aff(\R),m)}>0
$
for $\sigma_r$-a.e.\ $\xi_\theta$, by the computation of
Section~\ref{subsec.companion}. Hence $|\lambda_\star|=1$, contradicting the
choice of $\gamma_\star$. Therefore $P_\mu$ is totally dissipative.
\end{proof}

\end{document}